\pdfoutput=1
\documentclass{amsart} 
\usepackage{amssymb,amsmath,amsfonts,graphics,pinlabel,color,verbatim}

\begin{document}  \bibliographystyle{plain}

\newtheorem{thm}{Theorem}[section] \newtheorem{lem}[thm]{Lemma} \newtheorem{cor}[thm]{Corollary}
\newtheorem{mainlem}[thm]{Main Lemma} \newtheorem{prop}[thm]{Proposition} \newtheorem{conj}[thm]{Conjecture}
\newtheorem{rmk}[thm]{Remark}

\def\square{\hfill${\vcenter{\vbox{\hrule height.4pt \hbox{\vrule width.4pt height7pt \kern7pt \vrule width.4pt} \hrule height.4pt} } }$ }

\newenvironment{pf}{{\it Proof:}\quad}{\square \vskip 12pt}
\newenvironment{pfsketch}{{\it Sketch of Proof:}\quad}{\square \vskip 12pt}
\newenvironment{pfof}[1]{\vskip 12pt {\it Proof of #1:}\quad}{\square \vskip 12pt}

\theoremstyle{definition} \newtheorem{defn}{Definition}[section]

\newcommand{\Hn}{\ensuremath{\mathbb{H}^n}}
\newcommand{\Hf}{\ensuremath{\mathbb{H}^4}}
\newcommand{\Z}{\ensuremath{\mathbb{Z} } }

\newcommand{\U}{\ensuremath{\widetilde}}
\newcommand{\Hom}{\ensuremath{\text{Hom} } }
\newcommand{\Def}{\ensuremath{\text{Def} } }
\newcommand{\Isom}{\ensuremath{\text{Isom} } }
\renewcommand{\th}{^{\text{th} } }
\newcommand{\rect}{\ensuremath{\Gamma_{\text{rect}} } }
\newcommand{\oct}{\ensuremath{\Gamma_{\text{oct}} } }
\newcommand{\Poct}{\ensuremath{P_{\text{oct}} } }
\newcommand{\cube}{\ensuremath{\Gamma_{\text{cube}} } }

\def\urltilda{\kern -.15em\lower .7ex\hbox{\~{}}\kern .04em}

\renewcommand{\th}{^{\text{th} } }
\newcommand{\p}[1]{\ensuremath{\boldsymbol{#1^+} }}
\newcommand{\m}[1]{\ensuremath{\boldsymbol{#1^-} } }
\newcommand{\pandm}[1]{\ensuremath{\boldsymbol{#1}^{\pm}} }
\renewcommand{\l}[1]{\ensuremath{\boldsymbol{#1}} }
\newcommand{\pd}[1]{\ensuremath{\boldsymbol{\dot{#1}^+} } }
\newcommand{\md}[1]{\ensuremath{\boldsymbol{\dot{#1}^-} } }
\newcommand{\ld}[1]{\ensuremath{\boldsymbol{\dot{#1} } } }

\numberwithin{equation}{section}
\numberwithin{figure}{section}
\numberwithin{table}{section}


\title{From the hyperbolic $24$-cell to the cuboctahedron}
\author{Steven P. Kerckhoff and Peter A. Storm} \date{May 2008}

\begin{abstract}
We describe a family of $4$-dimensional hyperbolic
orbifolds, constructed by deforming an infinite volume orbifold
obtained from the ideal, hyperbolic $24$-cell by removing two
walls.  This family provides an infinite number of infinitesimally
rigid, infinite covolume, geometrically finite discrete subgroups
of $\text{Isom}(\mathbb{H}^4)$.  It also leads to {\em finite}
covolume Coxeter groups which are the
homomorphic image of the group of reflections in the hyperbolic
$24$-cell.
The examples are constructed very explicitly, both from an algebraic
and a geometric point of view.  The method used
can be viewed as a $4$-dimensional, but infinite volume, analog
of $3$-dimensional hyperbolic Dehn filling.
\end{abstract}

\thanks{Kerckhoff and Storm were partially supported by NSF grants DMS-0605151 and DMS-0603711,
respectively.} \maketitle

\section{Introduction}\label{introduction}

The study of hyperbolic manifolds or, more generally, discrete
subgroups of isometries of
hyperbolic space, typically divides along dimensional lines, between dimension $2$ where there
are deformations of co-compact respresentations and higher dimensions where there are no such
deformations due to Mostow rigidity.  However, there is an equally important distinction between
dimension $3$ and higher dimensions.  While Mostow-Prasad rigidity \cite{Mos,Pr}
(or Calabi-Weil local rigidity \cite{C,We}) guarantees that there are no nontrivial
families of complete, finite volume
structures in dimension at least $3$, if either the completeness or
finite volume hypothesis is dropped, in dimension $3$ there typically are deformations.

Specifically, a noncompact, complete, finite volume hyperbolic $3$-manifold $M$ with
$k$ cusps always has a $k$--complex dimensional deformation space, through necessarily
non-complete hyperbolic structures.  Topologically, $M$ is the interior of a compact
$3$-manifold $\bar M$ with $k$ torus boundary components.  Certain of the nearby incomplete
structures, when completed, give rise to complete hyperbolic structures on
{\it closed} $3$-manifolds which are obtained topologically by attaching solid
tori to the torus boundary components of $\bar M$.  This process is called
hyperbolic Dehn filling.  There are an infinite number of topologically
distinct ways to attach a solid torus to each boundary component and a fundamental
result of Thurston's \cite[Ch.5]{Th} is that, if a finite number of fillings are excluded from each
boundary component, the remaining closed manifolds have hyperbolic structures.
In particular, the process of hyperbolic Dehn filling gives rise to an infinite
sequence of closed hyperbolic $3$-manifolds whose volumes converge to that of
the cusped hyperbolic structure $M$.

This is in marked contrast with the situation in dimension at least $4$ where
a central result of Garland and Raghunathan \cite{GR} implies that a finite covolume
discrete subgroup of  $\text{Isom}(\mathbb{H}^n), ~ n \ge 4$, is isolated, up to
conjugacy, in its representation variety.  (Their result applies, in fact, to a much more
general class of Lie groups.)  In particular, there are no nontrivial
deformations of a finite volume hyperbolic structure, even through incomplete
structures, in dimension $4$ and higher; there is no higher dimensional analog
of hyperbolic Dehn filling for finite volume structures.
This is reflected by a result of Wang \cite{Wa} which states that, for any $V$ and fixed $n\ge 4$,
there are only a finite number of hyperbolic $n$-manifolds with volume at most $V$.

There is an equally dramatic distinction between dimension $3$ and higher dimensions
for infinite volume, geometrically finite hyperbolic structures.
When there are no parabolic elements, these structures can be compactified by
adding boundary components that inherit conformally flat structures from the
sphere at infinity.  In dimension $3$ these will consist of surfaces of genus
at least $2$ and Ahlfors-Bers deformation theory \cite{And} guarantees that
the deformation space of such structures has dimension equal to the sum of the
dimensions of the Teichmuller spaces of conformal structures on the boundary surfaces.
By contrast, in higher
dimensions the $(n-1)$-dimensional boundary can have a large dimensional
space of conformally flat deformations while the $n$-dimensional hyperbolic
structure is locally rigid; i.e., none of the deformations of the boundary
representation in $\text{Isom}(\mathbb{H}^n)$ extend over the interior.
In this paper we will construct an infinite family of examples with this
property.

In dimension $3$ the situation can be summed up in both the finite volume
and the infinite volume, geometrically finite cases by saying that a
half-dimensional subset of deformations of the boundary representations
(up to conjugacy) extends over the entire $3$-manifold.  In higher dimensions,
there are no deformations of finite volume complete structures, even through
incomplete structures, and there is no general theory of deformations in
the infinite volume case.
As a result, the study of hyperbolic structures in dimension $4$ and higher
currently consists largely of constructing and analyzing examples.
This paper can be viewed as fitting into that context.  However,
we believe that the method we use to construct these examples should
be applicable in much more general situations.

We will describe an infinite collection of $4$-dimensional hyperbolic
orbifolds.  They are actually part of a continuous family of representations
into $\text{Isom}(\mathbb{H}^4)$,
each of which corresponds to a singular hyperbolic structure, called a cone manifold.
Within this continuous family there are an infinite number of geometrically finite
discrete subgroups which are infinite covolume but infinitesimally rigid.
The process we use to create this family
can be viewed as a $4$-dimensional, but infinite volume, analog of
hyperbolic Dehn filling.  As such it can be viewed as an attempt to
mimic, to as great an extent possible, the situation in dimension
$3$ and to provide some insight into what a general theory of deformations
in dimension $4$ might look like.

The examples are all reflection groups; the discrete groups are hyperbolic
Coxeter groups.  The point of departure is the ubiquitous $24$-cell in its
realization as a right-angled, ideal $4$-dimensional hyperbolic polytope.
We denote the group of reflections in the $3$-dimensional walls of this
polytope by $\Gamma_{24}$.  It is a finite volume discrete subgroup of
$\text{Isom}(\mathbb{H}^4)$.  By the Garland-Raghunathan theorem it
has no deformations.  However, when we remove the reflections in two
disjoint walls (we will be more explicit later about which walls),
the resulting reflection group $\Gamma_{22}$ with infinite volume
fundamental domain becomes flexible.  We prove that there is
a smooth $1$-dimensional family of representations of $\Gamma_{22}$
which corresponds to a family of polytopes in $\mathbb{H}^4$.  These
are all the possible deformations of $\Gamma_{22}$ near the inclusion
representation.  Geometrically, certain pairs of walls that had been
tangent in the original polyhedron pull apart a finite nonzero distance,
and other previously tangent pairs of walls intersect.  The polytopes
retain a large degree of symmetry.  All the new angles of intersection
are equal and the walls which previously intersected orthogonally continue to do
so.  The new angle of intersection, denoted by $\theta$, can be used as
a parameter for the family of representations, which we denote by $\rho_{\theta}$.

When $\theta = \pi/m$, the group
of reflections in the walls of the polytope is a discrete subgroup of
$\text{Isom}(\mathbb{H}^4)$.  It is a homomorphic image of $\Gamma_{22}$
with new relations that do not exist in $\Gamma_{22}$.  In fact, because of 
the symmetry retained throughout the family, it is possible to conclude that, 
even when 
$\theta = \frac{2 \pi}{n}$, where $n \geq 3$ is an integer, the group 
$\rho_{\theta}(\Gamma_{22})$ is discrete.

Both $\Gamma_{24}$ and $\Gamma_{22}$ are right-angled Coxeter groups.
Although $\Gamma_{22}$ is most naturally viewed as a subgroup of $\Gamma_{24}$,
a well-known property of right-angled Coxeter groups implies that it is also a
{\em quotient group} of $\Gamma_{24}$.  Since the above discrete groups are all 
quotients of $\Gamma_{22}$
we obtain surjections of $\Gamma_{24}$ onto them for any $n\geq 3$.

We can summarize the main results of this paper as:

\begin{thm} \label{summary thm}
Let $G$ be the isometry group of hyperbolic $4$-space.  There is a discrete geometrically
finite reflection group $\Gamma_{22} < G$ and a smooth family of representations
$\rho_{\theta}: \Gamma_{22} \rightarrow G$  converging algebraically, as $\theta \to 0$, to the inclusion map.  The family has the following properties:
\begin{enumerate}
\item \label{1}When $\theta = \frac{2 \pi}{n}$, where $n \geq 3$ is an integer,
the representation $\rho_{\theta}$ is not faithful.  The image group
$\Lambda_n := \rho_{\theta}(\Gamma_{22}), \theta = \frac{2 \pi}{n}$
is a discrete geometrically finite subgroup of $G$.
\item \label{2}In the representation variety $\Hom(\Lambda_n,G)$ the inclusion map is infinitesimally rigid for all $n\ge 3$.
\item \label{3}When $n = 2 m$, $m \geq 4$, $\Lambda_n$ has infinite covolume.  Its convex core does not have totally geodesic boundary. 
\item \label{4} When $n = 2 m$, $m \geq 4$, the boundary subgroups of $\Lambda_m$ are convex cocompact and have nontrivial discrete faithful convex cocompact deformations.
\item \label{5}For all $n \geq 3$, there is a surjective homomorphism from the reflection group in the regular, ideal hyperbolic $24$-cell, $\Gamma_{24}$, onto $\Lambda_n$.
\end{enumerate}
\end{thm}

Representations in this family exist for $\theta \in [0,\pi]$, and they correspond to convex $4$-dimensional polytopes for $\theta \in [0, \pi)$.
For all $\theta \in (0, \pi/3)$ they are combinatorially the same;
only certain angles vary.  However, beginning at $\frac {\pi} {3}$,
there is a fairly drastic change in the combinatorics and for 
$\theta \in [\cos^{-1}(1/3),\pi)$, the polytope has {\em finite volume}.  

The finite covolume discrete groups $\Lambda_n$ in this family, when $n = 3,4,5$,
are of particular interest; so is the case when $n=6$.
At $\theta = \pi/3$, when $n = 6$, the topology of the boundary of the convex hull
changes and the resulting boundary components become totally geodesic.  This
is reminscent of the familiar process in dimension $3$ when curves on a higher genus
boundary surface are pinched, resulting in boundary components that are all
totally geodesic triply punctured spheres.  However, unlike the triply punctured
sphere case, the boundary components in our situation have non--totally geodesic
representations.
At $\theta = \cos^{-1}(1/3) \approx 1.231$
the boundary components of the polytopes completely disappear and for all larger
values of $\theta$ the polytopes have finite volume.  In particular, when
$n = 3,4,5$ the groups $\Lambda_n$ are actually  {\em finite covolume} discrete reflection groups.  Using a criterion of Vinberg \cite{Vin2} we can decide
whether or not they are arithmetic.

Finally, as $\theta \to \pi$ the $4$-dimensional polytopes collapse to
a $3$-dimensional hyperbolic ideal polyhedron, the right-angled cuboctahedron.
This is the end of the path referred to in the title of this paper.

\begin{thm} \label{summary thm 2}
Let $\rho_{\theta} : \Gamma_{22} \rightarrow G$ be the family of representations in 
Theorem \ref{summary thm} and  let $\Lambda_n := \rho_{\theta}(\Gamma_{22}), \theta = \frac{2 \pi}{n}$ be the image groups where $n \geq 2$ is an integer.
\begin{enumerate}
\item \label{1'} When $n = 3,4,5$ the discrete groups $\Lambda_n$ have finite
covolume.  They are non-uniform lattices.
\item \label{2'} When $n = 3,4$ the lattice $\Lambda_n$ is arithmetic.  When
$n=5$ it is not arithmetic.
\item \label{3'} When $\theta = \pi$ the image group $\rho_{\pi}(\Gamma_{22})$ preserves a
$3$-dimensional hyperplane, and has a degree two subgroup conjugate into $\Isom(\mathbb{H}^3)$ as the group of reflections $\Gamma_{co}$ in the ideal right-angled cuboctahedron.
\end{enumerate}
\end{thm}

The results in Theorems \ref{summary thm} and \ref{summary thm 2} are proved throughout the
body of the paper; the theorems themselves will not reappear explicitly outside of this introductory section.  Hence, we will provide here an outline of the paper with a guide to where the various
results are proved.

Section \ref{preliminaries} provides some background material on reflection groups, Zariski 
tangent spaces and orbifolds.  In Section \ref{24-cell} we give a detailed introduction to the
geometry and combinatorics of the hyperbolic $24$-cell, emphasizing the aspects used
throughout the paper.  In Section \ref{symmetries} we describe the symmetries of the $24$-cell,
focusing on those that remain even after we remove two walls to create the flexible reflection
group $\Gamma_{22}$.  We will ultimately find that all deformations of $\Gamma_{22}$ preserve this
group of symmetries, a fact that simplifies a number of the proofs.
In Section \ref{deformation prelims} we discuss the deformation theory of some auxillary
subgroups of $\Gamma_{22}$.  

In Section \ref{defining the deformation 2} we explicitly construct the family of representations
$\rho_{\theta}$ of the reflection group $\Gamma_{22}$.  Throughout the paper we use a different parameter $t$ because it is easier to write down the deformation using this parameter instead
of the angle $\theta$.  The formula relating the two parameters appears in Proposition \ref{angle prop}.  This, together with an easy application of the Poincare fundamental domain lemma, gives a proof of the existence of the infinite family of discrete reflection groups as stated 
in part \ref{1} of Theorem \ref{summary thm}.  
The family of representations is constructed by first assuming that
the symmetries of the polytope $P_{22}$ obtained by removing two walls of the $24$-cell are
preserved by the deformation.  We show that the representations in our family are the only ones preserving this group of symmetries.  It is only much later, in Sections \ref{inf letter} and 
\ref{computing infinitesimally}, where we show that in fact \emph{all} deformations preserve
these symmetries and, hence that we have constructed all possible local deformations.  This is used
to prove the smoothness of the family as well as
the infinitesimal rigidity of the groups $\Lambda_{n}$ (which is part (\ref{2}) of Theorem \ref{summary thm}).  Because the proof of these results is a fairly complicated computation of
Zariski tangent spaces, we have postponed it until after we have analyzed the geometry of the family.

In Sections \ref{viewed in the sphere at infinity} and \ref{fundamental domain} we analyze the
geometry of the polytopes corresponding to the representations in our family.  The geometric and combinatorial aspects of the polytopes and their walls are not at all apparent from the explicit description of the family which is given simply as a family of space-like vectors corresponding to the hyperplanes that determine the walls of reflection.  This analysis is necessary in order to make further conclusions about the discrete reflection groups in the family.  Since it focuses on the geometry of the various walls, which are $3$-dimensional, it also provides a link between the $3$ and $4$ dimensional deformation theories.  The results of these sections hold for angles $\theta$ strictly between $0$ and $\pi/3$; in this range the combinatorics of the walls is constant. 
Using this analysis, in Section \ref{miscellaneous section} we prove parts (\ref{3}) and (\ref{4}) of Theorem \ref{summary thm}.

In Section \ref{disappear} we begin the description of the more intricate change in the combinatorics and geometry of the polytopes that occurs between angles $\pi/3$ and $\pi/2$.  During this period the combinatorics of the infinite volume ends of the polytopes change and then these ends completely disappear.  
We provide a series of floating point based diagrams that suggest
how this process occurs.  In Section \ref{manual lattice} we give a rigorous proof of this process.
The proof depends on passing to representations of the extended group $\U{\Gamma}_{22}$ where the finite
group of reflective symmetries have been added in.  This group has fewer generators (though more complicated walls), allowing for a manageable analysis of the Coxeter diagrams of the image
groups $\rho_{\theta}(\U{\Gamma}_{22})$.  In this latter section, part (\ref{5}) of
Theorem \ref{summary thm} and parts (\ref{1'}) and (\ref{2'}) of Theorem \ref{summary thm 2} are
proved.  In Section \ref{cubeoctahedron section} we analyze the process of converging
to the terminal group $\rho_{\pi}(\Gamma_{22})$ and its relation to the cuboctahedron reflection group.

Finally, it is perhaps worth mentioning that this paper is not a completely accurate reflection of the process by which we discovered the family of polytopes described here.  Our presentation represents an assimilation of the material after attempting to understand initial computations in direct geometric terms, both contrasting it with and drawing analogies to the general $3$-dimensional theory. Originally, we found the family of examples by a combination of computational and graphical
experimentation, both computer-aided.  After removing one and then two walls from the
$24$-cell, we computed the dimension of the Zariski tangent space to the deformation space
to be $1$-dimensional in the latter case.  We then attempted
to construct an actual deformation corresponding to this Zariski tangent space.  Rigorous
computer-aided integer polynomial computation provided the proof of the existence of a substantial subinterval of our current family.   However, after better understanding the geometry of this
family and its symmetries, the use of the computer was removed from the proof.  The only remnants
of the computer assistance are the floating point diagrams in Section \ref{disappear} (which are
rigorously confirmed in a later section) and the numerous pictures which we continue to find
useful.

The authors benefited tremendously from the help of Daniel Allcock.  In particular, his
insights into the symmetries of the right-angled Coxeter group $\Lambda_4$ 
and how they can be
made apparent by a particular choice of the space-like vectors in Minkowski space led
directly to a significant simplification of our proofs.  We would also like to thank
Igor Rivin for several helpful conversations on topics related to the original
computer-aided computations.  Also thanks to Ian Agol and Alan Reid for their help
in understanding Vinberg's work on arithmetic reflection groups.

\section{Preliminaries}\label{preliminaries}

In this section we will introduce some of the basic concepts and terminology that will
be used throughout the paper.

\subsection{Hyperbolic reflection groups} \label{reflection prelims}

In this paper we will be studying $4$-dimensional hyperbolic reflection groups.  
A {\it reflection group} is, by definition, generated by elements $ \{r_i\} $ of order $2$
with relations of the form $ \{(r_i r_j)^{m_{ij}}\}$, where $i \neq j$ and
$m_{ij} \geq 2$ are integers.  In some texts, the convention $m_{ij} = \infty$ is used to signify that the element $r_i r_j$ is of infinite order; 
however, in this case we will simply omit any relation between $r_i$ and $r_j$.

An {\it $n$-dimensional hyperbolic representation} of a reflection group $\Gamma$ is a representation of 
$\Gamma$ in the group $G_n$ of isometries of hyperbolic $n$-space.  We will further require that the generators $r_i$ of $\Gamma$ are represented by reflections in codimension-$1$ totally geodesic hyperplanes.  We do not generally require the group generated by these reflections to be discrete.

It is worth emphasizing that our requirement that the generators of $\Gamma$ are sent to reflections 
in hyperplanes is an extra assumption.  Since there are
other types of order two elements in $G_n$, such representations may not include all
representations of  $\Gamma$ into $G_n$.   However,  since any order two element in $G_n$
near a reflection in a hyperplane is again such a reflection, this restriction is not
important locally.  Any representation near one of this form still has the same property.  
In fact, because the set of reflections in hyperplanes is a connected component of the
closed subset of elements of order two in $G_n$,  the set of such hyperbolic representations is
actually a union of components of the entire representation variety of $\Gamma$.  Throughout
this paper, a hyperbolic representation will mean a representation of $\Gamma$ in $G_n$
of  this particular type.

We will be using the Minkowski model of hyperbolic space and restrict our attention to the case
$n = 4$  since that is the dimension most relevant to this paper.    The Minkowski metric 
on $\mathbb{R}^5$ is determined by the inner product
$$( v, w ) := -v_0 w_0 + v_1 w_1 + v_2 w_2 + v_3 w_3 + v_4 w_4.$$ 
The resulting Minkowski space is denoted by $\mathbb{R}^{1,4}$.  Recall that a space-like vector has positive Minkowski squared norm, a light-like vector has Minkowski squared norm zero, and a time-like vector has negative Minkowski squared norm.  Hyperbolic $4$-space $\mathbb{H}^4$ is defined as the set of points with Minkowski squared norm $-1$ and a positive $0\th$ coordinate.  

We will denote the full isometry group of $\Hf$ simply by $G$.  It  is isomorphic to the two components of $O(1,4)$ which preserve the top component of the hyperboloid of points with Minkowski squared norm $-1$.  Viewing $O(1,4)$ as a subset of $5 \times 5$ matrices, it is an affine algebraic variety in $\mathbb{R}^{25}$.  Similarly,  a finite presentation of a group $\Gamma$ gives the 
representation variety $\Hom (\Gamma,G)$ the structure of a real algebraic variety.  Specifically, if 
$\Gamma$ is generated by $k$ elements with $\ell$ relations, the representation variety $\Hom (\Gamma,G)$ is defined as the subset of $G^k \subset \mathbb{R}^{25k}$ satisfying the $25 \ell$ polynomial equations produced by the relations.  

Because of the large dimensions of the spaces involved, such varieties are quite complicated to study.
However, in the case of a hyperbolic representation of a reflection group, there is a simple
observation that significantly reduces the dimension of the space in which representation
variety lives.  (It still is generally quite large, however.)  The observation is that a reflection in a hyperplane is uniquely determined by that
hyperplane, and that a hyperplane is uniquely determined by a space-like vector in $\mathbb{R}^{1,4}$.
More specifically, a hyperplane in $\Hf$ is a totally geodesic embedded copy of $\mathbb{H}^3$.  The order two isometry of $\Hf$ defined by reflecting in a hyperplane is uniquely determined by that hyperplane and conversely,  uniquely determines that hyperplane.    In the Minkowski model, a hyperplane is given by the intersection with $\Hf$ of a linear $4$-dimensional subspace of $\mathbb{R}^{1,4}$.  This linear subspace $W \subset \mathbb{R}^{1,4}$ is in turn determined by a space-like vector $q \in \mathbb{R}^{1,4}$ satisfying $(q,W)=0$.  A reflection isometry of $\Hf$ is thus determined by such a vector $q$ and the choice of $q$ is unique up to multiplication by a nonzero scalar.  A choice of
vector $q$ can actually be interpreted as a choice of one of the two half-spaces bounded by the hyperplane or as a choice of orientation on the hyperplane.  Multiplication by a {\it positive} scalar leaves this choice unchanged but it is flipped by multiplication by $-1$.

The relations $\{ (r_i r_j)^{m_{ij}} \} $ in a hyperbolic representation of reflection group can also be described by relations between such vectors.  To see this, denote by $W_i, W_j$ the reflection hyperplanes for the representatives of a  pair of generators $r_i, r_j$. 
When the hyperplanes are a positive 
distance apart, the element $r_i r_j$ represents translation by twice the distance between the hyperplanes along their common perpendicular geodesic.  When the hyperplanes are tangent
at infinity, the product is a parabolic translation.  In either case, when $W_i$ and $W_j$ are
disjoint, the product has infinite order.  When they intersect (and are distinct), 
they do so in a codimension-$2$ geodesic subspace.  The product $r_i r_j$ 
rotates around this subspace by twice the dihedral angle between the hyperplanes.
Thus,  a relation of the form $ (r_i r_j)^{m_{ij}}  = e$ 
is equivalent to the condition that the corresponding hyperplanes intersect in angle equal to an
integral multiple of $\frac {\pi} {m_{ij}}$.   Locally, the particular integral multiple will be constant
for each $\{i,j\}$ pair; in all of our examples, the angles will be exactly equal to $\frac {\pi} {m_{ij}}$.

A pair of  hyperplanes $W_i$ and $W_j$ are determined by a pair of space-like vectors $q_i$ and $q_j$ as described above.  The geometry of the hyperplane pair can be read from the geometry of the vectors $q_i$ and $q_j$.  In particular, if $W_i$ and $W_j$ intersect, then the dihedral angle $\theta$ formed by the $q_i$ side of $W_i$ and the $q_j$ side of $W_j$ satisfies the equation
\begin{equation}\label{angle equation} 
 \cos \theta = \frac{ - ( q_i, q_j ) }{ \sqrt{ (q_i,q_i) (q_j,q_j)}}.
\end{equation}
Similarly, if $W_i$ and $W_j$ do not intersect then the length $\ell$ of the shortest hyperbolic geodesic from $W_i$ to $W_j$ satisfies
\begin{equation}\label{distance equation}
  \cosh \ell = \frac{ - ( q_i, q_j ) }{ \sqrt{ (q_i,q_i) (q_j,q_j)}}.
\end{equation}
Note that the right hand sides of these equations are invariant under multiplication of the 
individual vectors by
{\it positive} scalars.   Multiplication of a vector by $-1$ corresponds to a change of side of its hyperplane;
this is an orientation issue that is best suppressed for the moment. 

This shows that a relation of the form $ (r_i r_j)^{m_{ij}}  = e$ for a hyperbolic representation 
of a reflection group is equivalent to a quadratic polynomial relation between the
entries of the corresponding space-like vectors, $q_i, q_j$, if the latter are normalized
to have unit (or fixed) Minkowski length.  This unit normalization is
also a quadratic condition.

Thus, it is possible to view the set of $4$-dimensional hyperbolic realizations of a
reflection group with $k$ generators and $\ell$ relations (ignoring the $k$ relations stating
the generators have order $2$) as a real algebraic variety in $\mathbb{R}^{5k}$ determined 
by $k+ \ell$ polynomial equations.

\subsection{Zariski tangent space and infinitesimal rigidity} \label{def:Zariski tangent space}

It is generally quite complicated to describe the simultaneous solutions to a large number of
polynomial equations in a large number of unknowns.  As a result, one is quickly led to consider instead
the corresponding infinitesimal problem.  If one knows a point $x$ in the simultaneous solution space of a
collection $P$ of polynomials,  one can look for tangent directions at $x$ in which the system continues to be
satisfied to first order.  This is done by differentiating the polynomials at $x$ and solving the resulting
linear system of equations.  The linear solution space is called the Zariski tangent space of the set of
polynomials at $x$.   The more modern terminology is to consider $P$ as determining a scheme $SP$;
this linear space is then called the Zariski tangent space to $SP$ at $x$.  We will use the following
simple formal definition: 
 
Let $P = \{ p_\alpha \}$ be a collection of real polynomials in $n$-variables.  Let  $x \in \mathbb{R}^n$  be a point 
where all of the $\{ p_\alpha \}$ simultaneously vanish.  Then the {\it Zariski tangent space at $x$ of the scheme $SP$}
is the linear space
$$T^Z_x SP := \{ v \in T_x \mathbb{R}^n \, | \ \textstyle\frac{d}{dt} \big|_{t = 0}\, p_\alpha ( x + t v ) = 0 \ \text{for all} \ \alpha \} \subseteq T_x \mathbb{R}^n.$$

In our case we are primarily interested in the actual real algebraic variety $V$ in
 $\mathbb{R}^n$ where a collection of polynomials $P$ simultaneously vanish.  While every
 polynomial in the ideal generated by the $P$ will certainly vanish on $V$, not every polynomial
 vanishing on $V$ is necessarily in this ideal.  Consider the simple example of the single polynomial
 $p(x,y) = x^2$ in $\mathbb{R}^2$.  The variety  $V$ in this case is the $y$-axis where $x=0$.  The polynomial
 $q(x,y) = x$ vanishes on $V$ but is not in the ideal generated by $x^2$.  This can have an effect on the
 computation of the Zariski tangent space.  For the scheme determined by $x^2$ it is $2$-dimensional at any
 point where $x=0$.  The Zariski tangent space for $V$, where one uses a generating set for the ideal
 of polynomials vanishing on $V$ (in this example, $q(x,y) = x$) in the above definition, the dimension
 is $1$ at every point on  $V$.
 
 In this example, the variety $V$ is a smooth manifold and the Zariski tangent space to $V$ is isomorphic to
 that of the smooth manifold.  This is generally not the case.  A standard example is the variety $V$ (and scheme)
 in the plane determined by $p(x,y) = x^2 -y^3$.  At the origin, the Zariski tangent space 
 to $V$ (and to the scheme) is $2$-dimensional.  The variety is not smooth there.  Most importantly, there
 are tangent directions in the Zariski tangent space that do not correspond to actual curves in $V$.  An
 infinitesimal solution to the equations need not correspond to a path of actual solutions.
 
 Fortunately, there is a criterion that deals with both of these issues, the difference between the scheme and
 variety Zariski tangent spaces, as well as the existence of a curve of solutions in any tangent direction.
 This criterion uses the implicit function theorem.  If one can find a smooth manifold of the same dimension
 as the Zariski tangent space of the scheme through a point on the variety, it follows that the variety is a smooth
 manifold of that dimension in a neighborhood of the point.  The two Zariski tangent spaces coincide with the
 tangent space of the manifold and every tangent vector corresponds to a smooth path of solutions.
 
 We will provide details in Section \ref{computing infinitesimally}, where this argument is used to analyze our nontrivial
 family of representations.  However, it is worth pointing out a classical use of this argument, due to Weil,
 with respect to the concept of infinitesimal  rigidity of a representation of a finitely presented discrete group 
 $\Gamma$ into a Lie group $G$.  Here $\Gamma$ and $G$ are quite general;  it suffices, for example,
 for $G$ to be real algebraic.  In particular, the argument applies to the situation considered in this paper. 

With $\Gamma$ and $G$ as above, let $\rho$ be a representation in $\Hom(\Gamma,G)$ which is a real algebraic
variety (scheme).  As  in the previous section, it can be viewed as a subset  of $G^n$ where $n$ is the number of generators
of $\Gamma$.   For $g \in G$ we can conjugate $\rho$ by $g$ to obtain the new representation $(g.\rho) (\gamma) := g \rho(\gamma) g^{-1}.$  We say that $\rho$ is infinitesimally rigid if the following is true: for any $v$ in the Zariski tangent space of $\Hom (\Gamma,G)$ at $\rho$ there is a path $g(t)$ in $G$ such that $g(0)$ is the identity element and 
$$\textstyle\frac{d}{dt} \big|_{t=0} \, (g(t).\rho) = v.$$  This condition is commonly described in terms of the
vanishing of the cohomology group $H^1(\Gamma; \mathcal G_{\rho})$ where the coefficients are in  the Lie algebra
$\mathcal G$ of $G$, twisted by the representation $\rho$.  The cocycles correspond to the Zariski tangent space
and the coboundaries to those tangent vectors induced infinitesimally by conjugation.

The point of this definition is that infinitesimal rigidity means every infinitesimal deformation of $\rho$ is obtained as the tangent vector to a path of representations obtained via conjugation.   Conjugate representations are not usually considered to be genuinely different representations of $\Gamma$.  The representation $\rho$ is therefore considered to be rigid at the infinitesimal level.

Weil's Lemma \cite{We} implies that if the centralizer of $\rho$
in $G$ is trivial, so that conjugation by $G$ locally determines a manifold in $\Hom(\Gamma,G)$ through $\rho$ 
with dimension equal to the dimension of $G$, then infinitesimal rigidity implies that near $\rho$ this manifold coincides with 
$\Hom(\Gamma,G)$.   Near $\rho$, $\Hom(\Gamma,G)$ is a smooth manifold and all nearby representations 
are conjugate.  The latter condition is called locally rigidity.  Thus, under a mild assumption  on the conjugation
action at $\rho$, Weil's Lemma implies that
an infinitesimally rigid representation is locally rigid.  The converse is occasionally false.

\subsection{Orbifolds} \label{def:orbifold}

The main topic of this paper is the representation space of a reflection group $\Gamma_{22}$ which has a discrete 
faithful $4$-dimensional 
hyperbolic representation coming from the group of reflections in the codimension-$1$ walls of an infinite volume polyhedron in $\Hf$.   We will find  a smooth $11$-parameter family of $4$-dimensional hyperbolic representations (not necessarily discrete) of this group.  
Dividing out by the conjugation action provides a smooth $1$-dimensional family of nontrivial deformations.  As discussed in Section \ref{reflection prelims}, this
family will be described algebraically as a family of space-like vectors.  
However, it will also be described geometrically as a family of groups generated by reflections in a family of varying polyhedra.  
These polyhedra will have new intersections between walls not existing in the original polyhedron and the combinatorics of the polyhedra will change a couple of times during the family.    At a countably infinite number of times all of the walls that intersect will do so at angles of the form $\pi/k$ for various integers $k$.  In these cases, the group generated by reflections in the walls of the polyhedron will be discrete with the polyhedron as a fundamental domain.  (This follows from  Poincar{\'e}'s lemma, which is discussed in Section \ref{defining the deformation 2}.)   These special hyperbolic representations of $\Gamma_{22}$ thus give rise to new reflection groups which are quotients of $\Gamma_{22}$, new relations having been added corresponding to new pairs of walls that intersect.

For all of these discrete hyperbolic reflection groups, the quotient space of $\Hf$ by the group can be usefully viewed as an (non-orientable) orbifold.  This will provide a third way of viewing these examples, which the authors find quite informative.  In particular,
some of the language used to describe the examples is most natural in this context.  For general background on
orbifolds, the reader is referred to \cite[Ch.13]{Th} or \cite{CHK}.  We will merely stress the essential points here, particularly those that are somewhat special to the current situation and might lead to confusion.

The general definition of a smooth orbifold  is as a topological space that is locally modeled on $\mathbb{R}^n$ (with its smooth structure) modulo a
(smooth) finite group action, with compatibility conditions on the overlaps.   The set of points where the local finite group is nontrival is called the {\it singular set} of the orbifold.  In our situation,  the underlying topological space of the quotient orbifold can be identified with the polyhedron itself.  The singular set is the union of all the reflection walls.  At  an interior point of a wall the local finite group is just the order $2$ group generated by a reflection.  At a point in the interior of a face where two walls intersect the local finite group is the dihedral group of the $k$-gon where the two walls intersect at angle $\pi/k$.  A similar analysis can be applied at points where more than two walls intersect.

Of course,  in our case the orbifold has more than a smooth structure.   The local group actions and overlap maps can all be modeled on restrictions of isometries of $\Hf$.  In this case, we say that the orbifold has a $4$-dimensional hyperbolic structure.
When the orbifold comes from a polyhedron that is not bounded, it will not be compact.  If the polyhedron bounds a finite volume region in hyperbolic space, the orbifold will have finite volume.  If not, the group of reflections in the sides of the polyhedron 
will possess a nontrivial domain of discontinuity on the sphere at infinity $S^3_{\infty}$ on which it acts properly discontinuously.  The quotient of this action will be a (union of) $3$-dimensional orbifolds.  These can be naturally attached to the original orbifold, creating an {\it orbifold with boundary}, with these as the boundary components.  The boundary components inherit a conformally flat structure from the sphere at infinity $S^3_{\infty}$.

To avoid a frequent point of confusion here, we wish to emphasize here that the walls of the original polyhedron are {\it not} part of the boundary of this orbifold with boundary.  One often refers the walls  as being ``mirrored".

This is an important point in understanding the analogy between the deformation theory of groups generated by reflections
in the walls of hyperbolic polyhedra and that of hyperbolic manifolds.   The former groups have torsion, but by Selberg's lemma \cite{Sel} they have finite index torsion-free subgroups.  The lifts of the walls of polyhedron will be internal in the corresponding manifold cover.  In particular, if the original polyhedron is bounded, the manifold will be closed.   There will be a nontrivial domain of discontinuity on the sphere at infinity for the group of reflections in the walls of the polyhedron if and only if  there is one in such a finite index manifold cover.   In general, the deformation theory for reflection groups from bounded, unbounded but finite volume, and infinite volume polyhedra respectively exactly parallels that of closed, noncompact finite volume, and infinite volume complete manifolds.  This holds in all dimensions.

For any discrete group $\Gamma$ of isometries of $\mathbb{H}^n$ one can take the limit set
of the action of $\Gamma$ (which is a subset of the sphere at infinity), form its hyperbolic convex hull, and remove the limit set itself.  It is a subset of $\mathbb{H}^n$,  in fact, the smallest closed convex subset invariant under $\Gamma$.  Its quotient space is called the {\it convex core} of $\Gamma$; it is a subset of $\mathbb{H}^n/\Gamma$.  Almost all of the discrete hyperbolic reflection groups we will consider will have infinite volume quotient spaces.  However, the volume of their convex cores will always be finite.  Such groups are called {\it geometrically finite}.  

The boundary of the convex core will be homeomorphic to an orbifold.  In dimension $3$, it has the further structure of a developable surface (or orbifold) that has been much studied.  In higher dimensions, the analogous structure is not well-understood.  In special cases, a component of the boundary of the convex core is totally geodesic.  In this case, we will say the corresponding end of $\mathbb{H}^n / \Gamma$ is a \emph{Fuchsian end}.

We will not attempt to study the structure of the boundary of the convex hull of our examples in detail.  However, we will be able to show that it is never totally geodesic, except at the beginning of the deformation and in one other case.

Finally, we remark that although the language of orbifolds only applies to the special polyhedra whose dihedral angles  are all of the form $\pi/k$,  all of the polyhedra can be viewed as hyperbolic cone manifolds.  For simplicity, we will not develop the appropriate language for this generality.  The interested reader is referred to \cite{CHK} where these concepts have been developed in somewhat more restrictive contexts.

\section{The hyperbolic $24$-cell} \label{24-cell}

 In this section we will introduce the basic properties of the hyperbolic $24$-cell.  For more background information one could consult \cite[Ch.10]{Cox1} or \cite[Ch.4]{Cox2}.  (Note that in the notation of \cite[Ch.10]{Cox1} the hyperbolic $24$-cell is the $\{ 3,4,3,4 \}$ regular hyperbolic honeycomb.)  In what follows, a face refers to a $2$-cell, and a wall to a $3$-cell. 

A natural place to begin the description is with the simpler Euclidean $16$-cell, which we will denote by $EP_{16}$.  (The $E$ stands for ``Euclidean.'')  It is a regular polytope in Euclidean $4$-space with vertices at the eight points $\{ \pm 2 e_i \}_i$, where $\{e_i \}_i$ is the standard basis for $\mathbb{R}^4$.  (The factor of $2$ will be convenient later.)  A nice way to imagine $EP_{16}$ is by starting with a $3$-dimensional octahedron with vertices at $\{ \pm 2 e_1, \pm 2 e_2, \pm 2 e_3 \}$, and suspending this octahedron to the points $\pm 2 e_4$ in $4$-space.  With this description it is clear that the sixteen walls of $EP_{16}$ are tetrahedra, each corresponding to a choice of sign for each of the four coordinates.  

Consider the set of $24$ points formed by taking the midpoint of each edge of $EP_{16}$.  These are precisely the $\binom{4}{2} \cdot 4 = 24$ points $\{ \pm e_i \pm e_j \}_{i \neq j}$. The Euclidean $24$-cell, $EP_{24}$, can be defined as the convex hull of these $24$ points.  The walls of this polyhedron can be visualized as follows:  

Consider the set of six points formed by taking the midpoint of each edge of a tetrahedral wall of  $EP_{16}$.  These six points are the vertices of an octahedron sitting inside the tetrahedron.  Thus, for each tetrahedral wall of $EP_{16}$,  there is a corresponding octahedron inside it.    These sixteen octahedra form sixteen walls of $EP_{24}$.  
A tetrahedral wall of $EP_{16}$ corresponds to a choice of sign for each of the four
coordinate vectors.  The $6$ vertices of the corresponding octahedron are formed by taking all possible sums of exactly two from the chosen set of vectors.  For example, with a choice of positive
sign for all four vectors, the vertices of the octahedron would be $e_1 + e_2, e_1 + e_3, e_1 +e_4,
e_2 + e_3, e_2 + e_4, e_3 + e_4$.   Two such octahedra intersect in a face when their corresponding choices of signs differ in precisely one place.  This provides us with a way to divide these
walls into two mutually disjoint octets: those corresponding to an even number of negative signs and
those corresponding to an odd number of negative signs.

Each of the above octahedra intersects $4$ others in a triangle contained in a face of the
corresponding tetrahedral wall of $EP_{16}$.  There are $4$ other triangular faces of such an octahedron, each associated to a vertex of the tetrahedron.  These faces arise as intersections with
the remaining eight walls of $EP_{24}$, which are formed by the links of the eight vertices of $EP_{16}$.  In order to visualize this, focus on a single vertex $v$ of $EP_{16}$.  Six edges of $EP_{16}$ emanate from $v$.  Take the midpoint of each such edge.  The convex hull of these six midpoints is an octahedron which we identify with the link of $v$.  These eight links, one for each vertex of $EP_{16}$, provide the remaining eight walls of $EP_{24}$.  Such a vertex is determined by a single, unit coordinate vector (with sign); the vertices of the associated octahedron consist of all possible sums of that vector with the remaining signed coordinate vectors, excluding the negative of the vector itself. 

With this description one can count the cells of all dimensions to conclude that $EP_{24}$ has $24$ vertices, $96$ edges, $96$ faces, and $24$ walls.  It is combinatorially self-dual.    

The Euclidean polytope $EP_{24}$ has a hyperbolic analogue $P_{24} \subset \mathbb{H}^4$.  
$P_{24}$ is a finite volume ideal hyperbolic polytope.  Being ideal, its outer radius is infinite.  The hyperbolic cosine of its inner radius is $\sqrt{2}$.  Most importantly, we will see that the interior dihedral angle between two intersecting walls is always $\pi/2$ (see also  \cite[Table IV]{Cox2}). 

This hyperbolic polytope can be most simply defined using the Minkowski model of hyperbolic space (see Section \ref{reflection prelims}).  
First, we isometrically embed Euclidean $4$-space as the last four coordinates of Minkowski space.  This puts $EP_{24}$ isometrically inside Minkowski space.  The next step is to ``lift'' $EP_{24}$ into $\mathbb{H}^4$.
For each vertex $v \in \{ \pm e_i \pm e_j \}_{i \neq j}$ of $EP_{24}$ (here $i$ and $j$ take values from $1$ to $4$) define the light-like (i.e. Minkowski norm zero) vector $ \sqrt{2} e_0 + v$.  This collection of light-like vectors can be thought of as ideal points on the boundary of $\Hf$.  We then define $P_{24}$ as the hyperbolic convex hull of these $24$ ideal points.   

The walls of this polyhedron are described combinatorially exactly as in the Euclidean case.  There are $16$ that correspond to a choice of sign for each of the coordinate vectors, $\pm e_i, i = 1,2,3,4$.  The ideal vertices of a wall corresponding to such a choice are the vectors  $ \sqrt{2} e_0 + v$, where
$v$ is one of the six vectors obtained as a sum of two of the chosen coordinate vectors.  The remaining $8$ walls are each determined by choosing a single signed coordinate vector $\pm e_i$ and letting $v$ cycle through sums with the remaining $6$  signed coordinate vectors $\pm e_j, j \neq i$.

We would like to see that the walls defined in this way have the same pairwise intersection properties as in the Euclidean case and that when two walls do intersect, they do so at right angles.  To do so,
for each wall $W_i$ we find a space-like vector $q_i$ with the property that it is orthogonal to all points in $W_i$.  (Throughout this discussion, all inner products, norms, and notions of angle are defined in terms of the Minkowski quadratic form.)  As discussed in Section \ref{reflection prelims}, $q_i$ is unique up to scale
and the relative geometry between two walls $W_i, W_j$ can be read off from the inner product
between the corresponding $q_i, q_j$. 

Given our previous description of the ideal vertices that determine each of the walls, it is
easy to immediately write down the list of the corresponding space-like vectors.  Given a
collection of $6$ light-like vectors which determine a wall, we must find a $q_i$ which is
orthogonal to all $6$ of them.  Note that $4$ such vectors generically determine a hyperplane,
and hence an orthogonal space-like vector.  So the fact that we will be able to find our collection of $q_i$ reflects the fact that the light-like vectors are not at all in ``general position".

All of the ideal vertices come from light-like vectors of the form $ \sqrt{2} e_0 \pm e_i \pm e_j$,  where
$i,j$ are distinct values from $1$ to $4$.   One type of collection of $6$ such vectors comes from
choosing a sign for each of $\pm e_1, \pm e_2, \pm e_3, \pm e_4$ and then restricting $\pm e_i$
and $\pm e_j$  to come from those choices.    The corresponding $q_i$ is then just
$\sqrt{2} e_0 \pm e_1 \pm e_2 \pm e_3 \pm e_4$ where one takes the same choices of sign.
Orthogonality is easy to check since the contribution to the dot product from the $0\th$ coordinate
is always $-2$ and  since the signs agree,  the contribution from the last $4$-coordinates
is always $+2$.  The other type of vector hextet comes from fixing $\pm e_i$ and letting $\pm e_j,  j\neq i$ vary over the remaining possibilities.  It is again easy to check that the vector
$e_0 \pm \sqrt{2} e_i$, with the same choice of $\pm e_i$, is orthogonal to all six of these vectors.

\begin{figure}[ht]
\labellist
\small\hair 2pt
\pinlabel $\m{0}$ at 183 177
\pinlabel $\p{0}$ at 226 235
\pinlabel $\p{3}$ at 176 101
\pinlabel $\m{3}$ at 214 6
\pinlabel $\p{1}$ at 108 179
\pinlabel $\m{1}$ at 3 238
\pinlabel $\m{2}$ at 103 99
\pinlabel $\p{2}$ at 7 14
\pinlabel $\p{5}$ at 214 197
\pinlabel $\m{5}$ at 318 272
\pinlabel $\m{4}$ at 222 117
\pinlabel $\p{4}$ at 316 68
\pinlabel $\m{6}$ at 131 194
\pinlabel $\p{6}$ at 90 275
\pinlabel $\p{7}$ at 138 107
\pinlabel $\m{7}$ at 88 64
\pinlabel \textcolor{red}{$\l{A}$} at 130 134
\pinlabel \textcolor{red}{$\l{B}$} at 222 151
\pinlabel \textcolor{red}{$\l{C}$} at 162 209
\pinlabel \textcolor{red}{$\l{G}$} at 163 158
\pinlabel \textcolor{red}{$\l{F}$} at 196 159
\pinlabel \textcolor{red}{$\l{E}$} at 102 151
\pinlabel \textcolor{red}{$\l{D}$} at 160 77
%
\endlabellist
\centering
\includegraphics{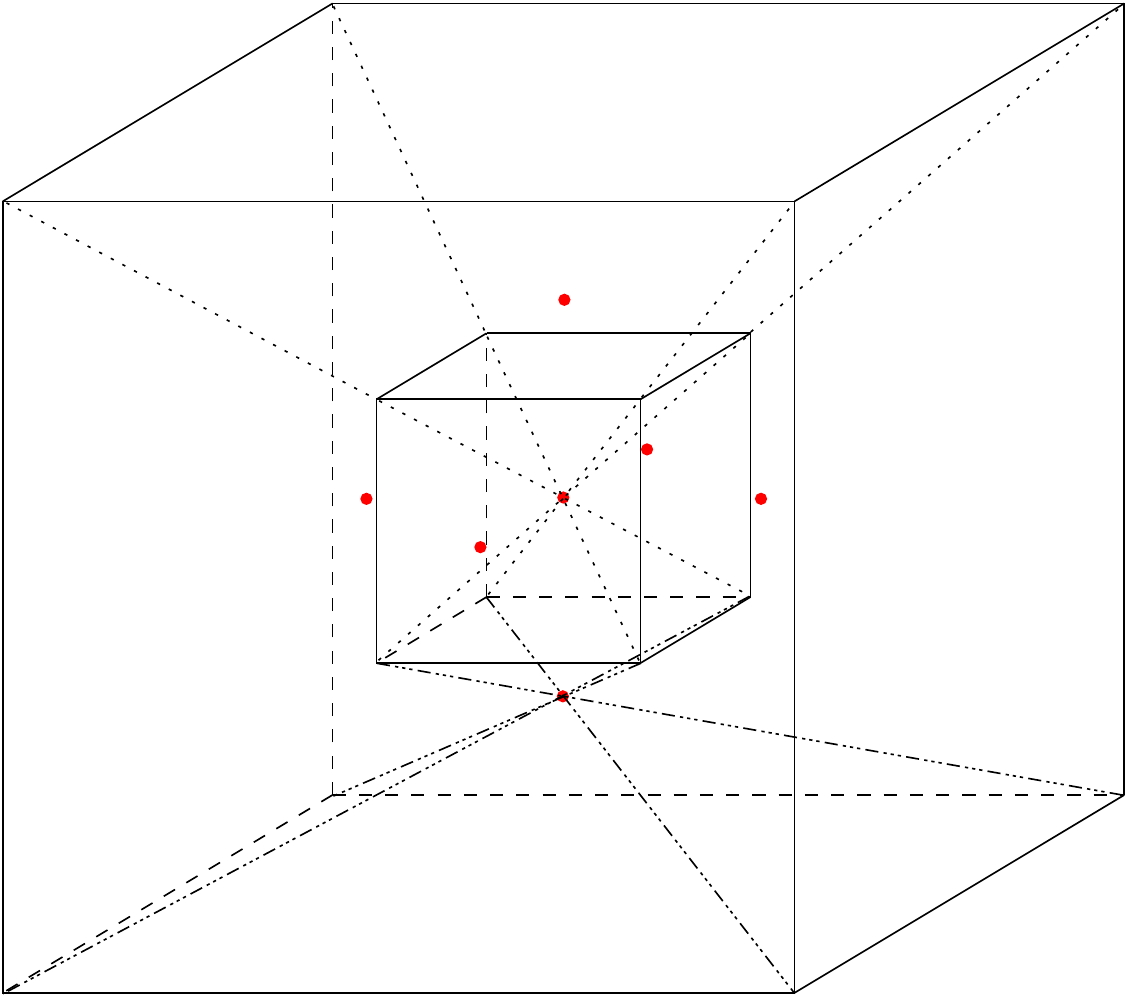}
\caption{A representation of the $24$-cell.  Vertex $\l{H}$ is at infinity.}
\label{cube figure}
\end{figure}
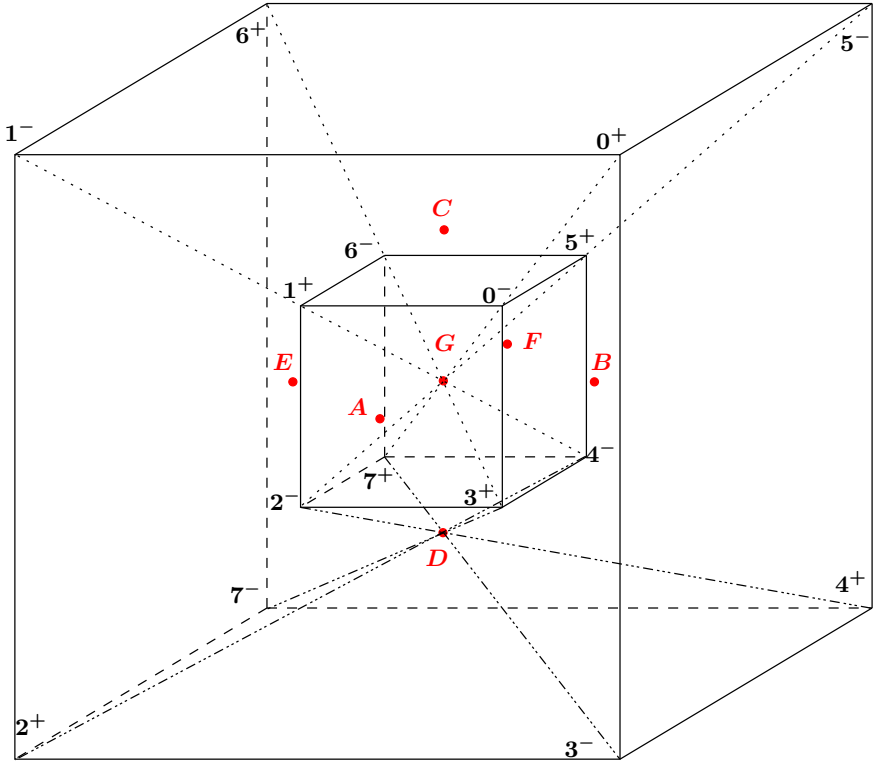

Thus our entire collection of vectors $q_i$ corresponding to the walls  of $P_{24}$ is given by:

$$\left\{ \left( \sqrt{2}, \pm 1, \pm 1, \pm 1, \pm 1\right)\right\} \bigcup \left\{ e_0 \pm \sqrt{2} e_j \right\}_{1 \le j \le 4}.$$ 
Due to their importance, we explicitly index these space-like vectors.  A list of the vectors $\{ q_i \}$ is given in table \ref{q table}, where the indices are not just integers and we have only written the index itself.   (Formally we should write $q_{\p{7}}$ rather than $\p{7}$, but we have chosen the latter simpler notation.)  The walls indexed by a number and a sign correspond to the first collection of $16$ vectors and those by a letter to the second collection of $8$ vectors.  The sign of the first collection of vectors refers to the parity of the number of negative signs in the vector.  They will be referred to as the
{\it positive} and {\it negative} walls, respectively.  The final group of $8$ will be referred to as the
{\it letter} walls.  The specific choices of numbers and letters will hopefully become more apparent once we introduce a visual device for cataloging them in figure \ref{cube figure}.

\begin{table}
\begin{eqnarray*}
\p{0} = \left( \sqrt{2},1,1,1,1 \right) , & &
\m{0} = \left( \sqrt{2},1,1,1,-1 \right),\\
\p{1} = \left( \sqrt{2},1,-1,1,-1\right),& &
\m{1} = \left( \sqrt{2},1,-1,1,1\right),\\
\p{2} = \left( \sqrt{2},1,-1,-1,1\right),& &
\m{2} = \left( \sqrt{2},1,-1,-1,-1\right),\\
\p{3} = \left( \sqrt{2},1,1,-1,-1\right),& &
\m{3} = \left( \sqrt{2},1,1,-1,1\right),\\
\p{4} =\left(\sqrt{2},-1,1,-1,1\right),& &
\m{4} =\left(\sqrt{2},-1,1,-1,-1\right),\\
\p{5} = \left( \sqrt{2},-1,1,1,-1\right),& &
\m{5} =\left(\sqrt{2},-1,1,1,1\right),\\
\p{6} =\left(\sqrt{2},-1,-1,1,1\right),& &
\m{6} =\left(\sqrt{2},-1,-1,1,-1\right),\\
\p{7} =\left(\sqrt{2},-1,-1,-1,-1\right),& &
\m{7} =\left(\sqrt{2},-1,-1,-1,1\right),\\
\l{A} =\left(1,\sqrt{2},0,0,0\right),& &
\l{B} =\left(1,0,\sqrt{2},0,0\right),\\
\l{C} =\left(1,0,0,\sqrt{2},0\right),& &
\l{D} =\left(1,0,0,-\sqrt{2},0\right),\\
\l{E} =\left(1,0,-\sqrt{2},0,0\right),& &
\l{F} =\left(1,-\sqrt{2},0,0,0\right),\\
\l{G} =\left(1,0,0,0,-\sqrt{2}\right),& &
\l{H} =\left(1,0,0,0,\sqrt{2}\right).
\end{eqnarray*}
\caption{Space-like vectors describing the $24$-cell}\label{q table}
\end{table}

This notation naturally divides the $24$ walls of the $24$-cell into three octets: the positive, negative, and letter walls.   They will play a special role at several points in this paper.

In particular, we will now check that walls in the same octet are disjoint (in $\Hf$).   (Here we are considering asymptotic walls to be disjoint.)   We will also see that walls which do intersect do so orthogonally.
This is easily done by considering the inner products between the corresponding space-like vectors
and using equations \ref{angle equation} and \ref{distance equation}.

For any pair of distinct vectors of the form $\left( \sqrt{2}, \pm 1, \pm 1, \pm 1, \pm 1\right)$, 
the contribution to the dot product from the last $4$ coordinates is $+2$, $0$, $-2$, $-4$ depending on 
whether they differ by $1$, $2$, $3$, or $4$ sign changes.  Since the contribution from the $0\th$ coordinate
is always $-2$, it is evident that two such walls are orthogonal if and only if they differ by
a single sign change and that they are otherwise disjoint.  In particular, all the positive walls are disjoint from each other as are all the negative walls.   It is worth pointing out further that
when they differ by exactly $2$ sign changes, they are tangent at infinity, and otherwise
they are a positive distance apart.  (Note that these vectors have Minkowski square norm equal to $+2$.)

Similarly, it is easy to see that the letter walls, corresponding to vectors of the form $e_0 \pm \sqrt{2} e_j$,  are mutually disjoint.  All pairs of letter walls are tangent except for pairs with the same $e_j$ and opposite sign.  The letter walls are orthogonal to those positive and negative walls with the same sign in the $j$th place.  There are $4$ positive and $4$ negative walls with this property.  All other positive
and negative walls are a positive distance away.

We note that the above pairwise intersection pattern agrees with that of the Euclidean $24$-cell.
One  can similarly check that the combinatorial pattern around lower dimensional cells also
agrees, where one needs to adapt this notion in an obvious way when vertices are ideal.
We will not do that here.

We now consider the group $\Gamma_{24} \subset G$ generated by reflections in the $24$ walls 
of the hyperbolic $24$-cell $P_{24}$.  It is a discrete subgroup of $G$ with finite covolume.  The
polyhedron $P_{24} \subset \Hf$ can be taken as a fundamental domain.  The quotient
of $\Hf$ by $\Gamma_{24}$ can be  viewed as an orbifold whose singular set is the entire collection of walls of $P_{24}$ which are all considered  mirrored.

Now that we understand the geometry of the intersections of the walls of $P_{24}$, we can
quickly describe a presentation for $\Gamma_{24}$.  It has $24$ generators $r_i, i=1,2,\dots, 24$,
each of order $2$.  The remaining relations are all of the form $(r_i r_j)^2 = e$ for any pair of
walls $W_i, W_j$ that intersect in $\Hf$; the exponent of $2$ arises from the fact that the walls intersect in angle $\frac {\pi} {2}$.  Since $r_i, r_j$ are of order two these relations are equivalent to saying that such pairs of $r_i$ and $r_j$ commute.  These pairs were seen above to correspond to a pair of walls consisting of one  positive and one negative wall whose corresponding space-like vectors differ by a single sign change or to a pair consisting of a letter wall whose associated space-like vector is of the form $e_0 \pm \sqrt{2} e_j$ and a positive or negative wall whose associated vector has the same sign in the $j$th place.  There are a total of $96$ such pairs.

To further aid our understanding of the geometry and combinatorics of the $24$-cell as well
as the structure of this reflection group,  we will now present two pictorial models.

The first model is shown in figure \ref{cube figure}.  The $24$ walls $\{ W_i \}$ of $P_{24}$ are represented by the $23$ vertices labeled in the figure, plus a point at infinity representing wall $\l{H}$.  To avoid cluttering the figure we have not drawn all the edges.  First, for each of the drawn edges, add in its orbit under the symmetries of the cube.
Then add edges connecting vertices $\p{0}$, $\m{1}$, $\p{2}$, $\m{3}$, $\p{4}$, $\m{5}$, $\p{6}$, and $\m{7}$ to the point at infinity representing wall $\l{H}$.  The resulting labeled graph $IG_{24}$ encodes how the walls of $P_{24}$ intersect: two walls intersect (orthogonally) if and only if their respective vertices are joined by an edge.

To help parse figure \ref{cube figure}, let's explain some properties of the vertex set.  If vertex $\l{G}$ is placed at the origin, then the red vertices lie on the coordinate axes.  Specifically, vertices $\l{A}$ and $\l{F}$ lie on the $x$-axis, vertices $\l{B}$ and $\l{E}$ lie on the $y$-axis, and vertices $\l{C}$ and $\l{D}$ lie on the $z$-axis.  Vertex $\l{H}$ is at infinity, and the remaining vertices all lie on a vertex of a cube.  The two cubes are both centered on vertex $\l{G}$.

The corresponding space-like vectors can also be read off from this figure.  The $0th$ coordinate is determined simply by whether it corresponds to a lettered wall, in which case it is $+1$, or a
numbered wall, in which case it is $+\sqrt{2}$.  The final coordinate, the coefficient for $e_4$,
plays a special role.  For the numbered walls it is $+1$ for those on the outer cube and
$-1$ for those in the inner cube.  For the letter walls  it is $0$ except for walls $\l{G}, \l{H}$,
at the origin and infinity respectively in the figure, which have coefficient $-\sqrt{2}$,
$+\sqrt{2}$, respectively.   These two letter walls will be removed in order to allow the
configuration to become flexible, so they play a limited role in what follows.

The remaining coordinates, the coefficients of $e_1,e_2, e_3$, closely correspond to
$x,y,z$ coordinates in the figure.  Specifically, if one rescales both cubes so they both
have vertices at the $8$ points $(\pm 1, \pm 1, \pm 1)$, then  these are the $(e_1, e_2, e_3)$ 
coordinates of the corresponding space-like vectors.  The orientation is chosen so that the
positive $x$- axis comes out of the page, the positive $y$-axis  to the right,  and the 
positive $z$-axis points up.   Similarly, the $(e_1, e_2, e_3)$ coordinates  of vectors 
associated to the remaining letter walls are determined by the positions of the red vertices along the coordinate axes if one scales them to be distance $\sqrt{2}$ from the origin.

Recall that the polytope $P_{24} \subset \Hf$ is an ideal polytope with its vertices at the boundary at infinity of $\Hf$.  Each end of $P_{24}$ is isometric to a warped product of a Euclidean cube and the half-line.  The vertices of $P_{24}$ correspond to the $3$-cells of figure \ref{cube figure}, which we now define.  For each complete graph on $3$ vertices in the graph $IG_{24} \subset \mathbb{R}^3$ add a Euclidean triangle with boundary the complete graph.  (These instructions require some obvious modifications for dealing with the point at infinity.  These are left to the reader.)  The resulting $2$-complex cuts $\mathbb{R}^3$ into $24$ regions which are combinatorial octahedra.  
For example, the vertices $\m{0}$, $\p{1}$, $\p{5}$, $\m{6}$, $\l{C}$, and $\l{G}$ ``bound'' such an octahedron.  Less obviously, $\p{5}$, $\m{5}$, $\m{6}$, $\p{6}$, $\l{C}$, and $\l{F}$ also bound an octahedron.  
With this extra structure we can now say that the vertices of $P_{24}$ correspond bijectively to these $24$ octahedra.  Moreover, for a given vertex $v$ of $P_{24}$, the six walls which are asymptotic to $v$ are given by the six vertices of $IG_{24}$ forming the respective octahedron.  
For example, the walls corresponding to $\m{1}$, $\p{2}$, $\p{6}$, $\m{7}$, $\l{E}$, and $\l{H}$ are all asymptotic to a single common vertex of $P_{24}$.  The group generated by the reflections in these $6$ walls is abstractly isomorphic to the Euclidean reflection group in the sides of cube.  Finally, two walls of $P_{24}$ are tangent at infinity if and only if the corresponding vertices of $IG_{24}$ are contained in a single octahedron and are not adjacent.  With this it is possible to determine the relative geometry of any pair of walls of $P_{24}$ using figure \ref{cube figure}.  For example, the following pairs of walls are tangent at infinity: $\{W_{\p{1}}, W_{\p{3}}\}, \{ W_{\l{C}}, W_{\l{F}} \}, \{ W_{\p{1}}, W_{\p{6}} \}$, and $ \{W_{\l{H}}, W_{\l{D}} \}$.

Consider the set of walls orthogonal to wall $W_{\l{A}}$.  Using figure \ref{cube figure} one can see they are the walls indexed by the set $J=\{\p{0},\m{0},\p{1},\m{1},\p{2},\m{2},\p{3},\m{3} \}$.  We know from the above description of $EP_{24}$ that $W_{\l{A}}$ is an ideal octahedron.  From the fact that all intersections are orthogonal, it follows that the dihedral angles of $W_{\l{A}}$ are also $\pi/2$.  Moreover, a reflection in a wall from the set $J$ will preserve the hyperplane containing $W_{\l{A}}$, which is a copy of $\mathbb{H}^3$.  This implies that reflections in the walls $\{ W_j \}_J$ generates a $3$-dimensional hyperbolic reflection group, the reflection group corresponding to a regular ideal hyperbolic octahedron with dihedral angles $\pi/2$.  Since $P_{24}$ is a regular polytope, an identical analysis applies to any other wall, they are all ideal octahedra with orthogonal dihedral angles.

Let us now present a second model for imagining $P_{24}$.  This model will take place in the $3$-sphere at infinity $\mathbb{S}^3_\infty$ of $\Hf$.  A hyperbolic hyperplane in $\Hf$ intersects $\mathbb{S}^3_\infty$ in a conformally round $2$-sphere.  Using stereographic projection, which is a conformal map, we can map the $3$-sphere minus a point to $\mathbb{R}^3$.  $2$-spheres passing through the point at infinity will appear under stereographic projection as affine planes.  Other conformally round $2$-spheres in $\mathbb{S}^3_\infty$ will appear as round spheres in $\mathbb{R}^3$.  A pair of orthogonal hyperplanes will intersect $\mathbb{S}^3_\infty$ in a pair of spheres intersecting orthogonally.  A pair of asymptotic hyperplanes will intersect $\mathbb{S}^3_\infty$ in a pair of tangent spheres.

\begin{figure}[ht]
\labellist
\small\hair 2pt
\pinlabel $\m{0}$ at 262 192
\pinlabel $\p{0}$ at 65 263
\pinlabel $\p{3}$ at 191 262
\pinlabel $\m{3}$ at 262 65
\pinlabel $\p{1}$ at 130 263
\pinlabel $\m{1}$ at 8 129
\pinlabel $\m{2}$ at 250 131
\pinlabel $\p{2}$ at 130 9
\endlabellist
\centering
\includegraphics{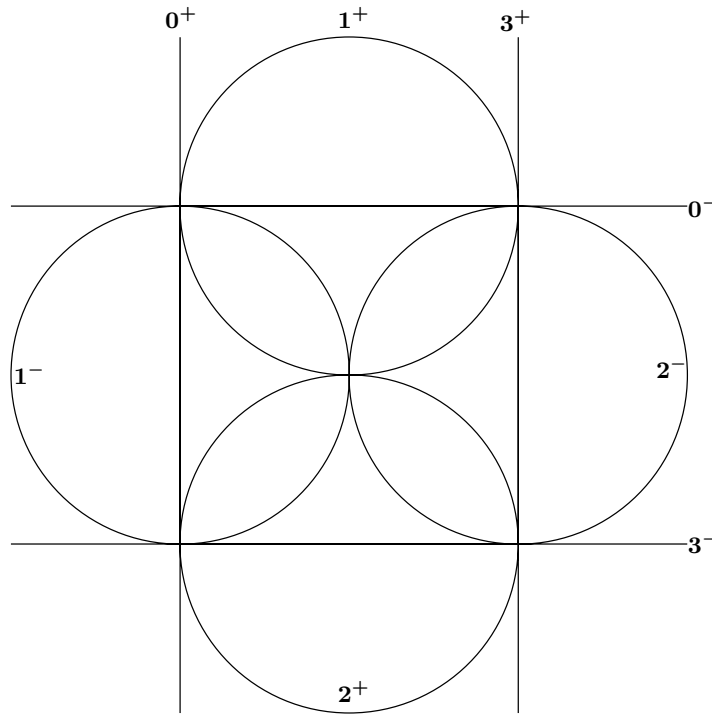}
\caption{The hyperplanes of an octahedron, viewed at infinity}
\label{circles1 figure}
\end{figure}

Each wall $W_i$ of $P_{24}$ extends to a hyperplane $H_i$ of $\Hf$.  Using stereographic projection, we can imagine how these hyperplanes will intersect $\mathbb{S}^3_\infty$.  Before considering all $24$ walls of $P_{24}$, let us begin with an easier example one dimension down.  Consider again the walls $\{W_j \}_{  J}$ intersecting $W_{\l{A}}$.  Each $H_i$ intersects $\mathbb{S}^3_\infty$ in a $2$-sphere $S_i$.  We will envision these spheres under stereographic projection in $\mathbb{R}^3$.  We may choose the coordinates of the stereographic projection so that $S_{\l{A}}$ is the $xy$-plane.  A moment's thought will reveal that the only possible arrangement of $\{ S_j \}_J$, up to conformal automorphism of $\mathbb{S}^3_\infty$ preserving $S_{\l{A}}$, is as shown in figure \ref{circles1 figure}.  The figure shows only the curves of intersection of the spheres of $\{ S_j \}_J$ with $S_{\l{A}}$.  The actual spheres of $\{S_j \}_J$ are determined uniquely by figure \ref{circles1 figure} because they all intersect $S_{\l{A}}$ orthogonally.

\begin{figure}[p!]
\includegraphics[scale=0.4,angle=-90]{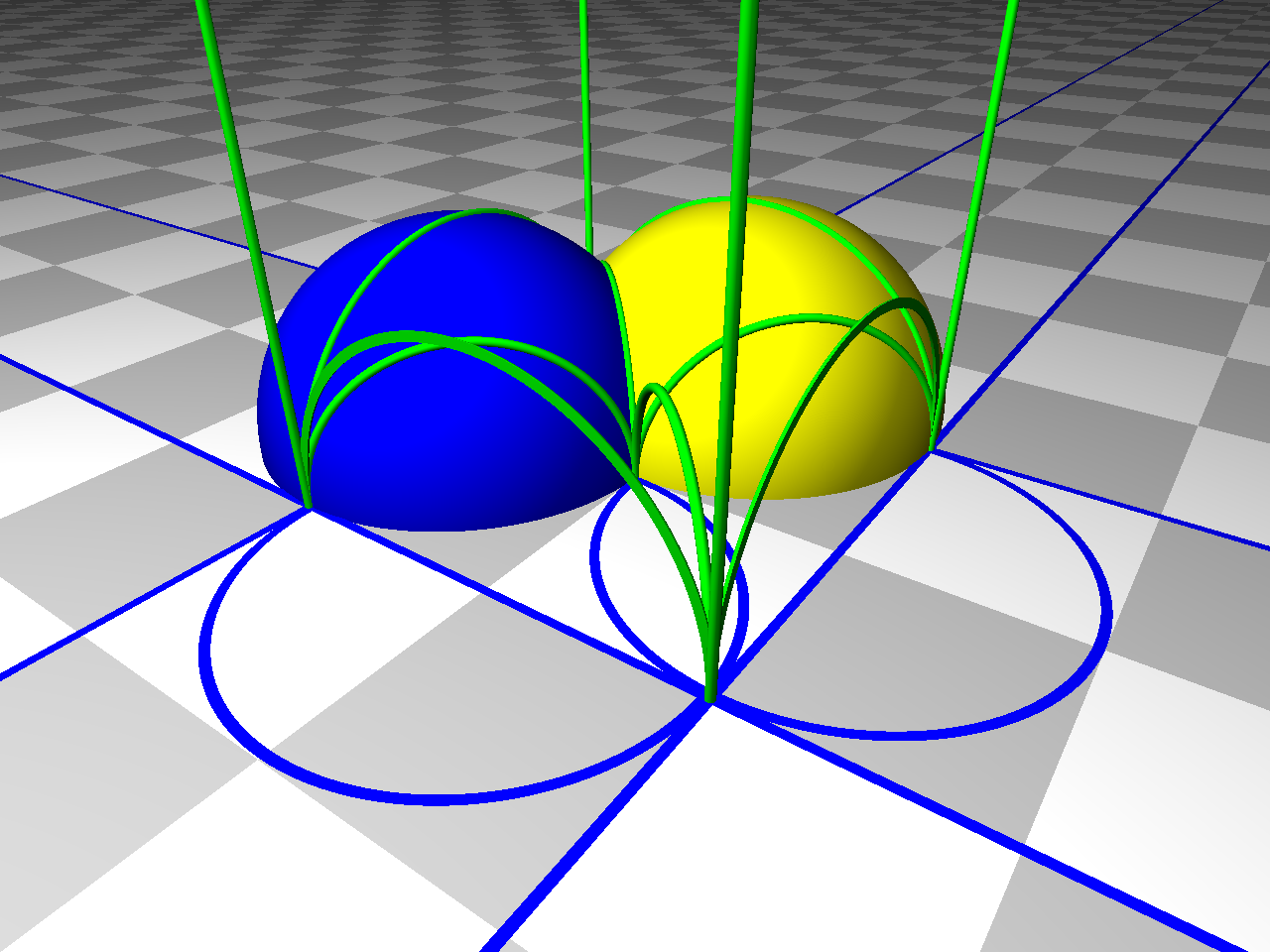}
\caption{The ideal octahedron $W_{\l{A}}$}
\label{octahedron2 figure}
\end{figure}

The hyperplane $H_{\l{A}}$ is a copy of hyperbolic $3$-space which we can imagine as the upper half-space lying over the $xy$-plane of figure \ref{circles1 figure}.  The intersections $H_{\l{A}} \cap H_{j}$ for $j \in J$ will be $2$-dimensional hyperplanes in this upper half-space.  As described above, these $2$-dimensional hyperplanes bound an ideal octahedron which is $W_{\l{A}}$.  To aid the readers imagination, we have included figure \ref{octahedron2 figure}, which shows the edges of the ideal octahedron $W_{\l{A}}$ and a pair of $2$-dimensional hyperplanes.  The checkerboard is the $xy$-plane of figure \ref{circles1 figure}, which represents the sphere at infinity of $H_{\l{A}}$.  The blue circles and lines form a copy of figure \ref{circles1 figure}.  The green edges are the edges of the ideal octahedron.  The green edge furthest from the viewer passes through the origin.  The blue hemisphere is the intersection $H_{\p{2}} \cap H_{\l{A}}$.  The yellow hemisphere is the intersection $H_{\m{1}} \cap H_{\l{A}}$.  We encourage the reader to examine figure \ref{octahedron2 figure} until its meaning is clear.  Several variations will appear later.

\begin{figure}[ht]
\centering
\includegraphics[scale=0.20]{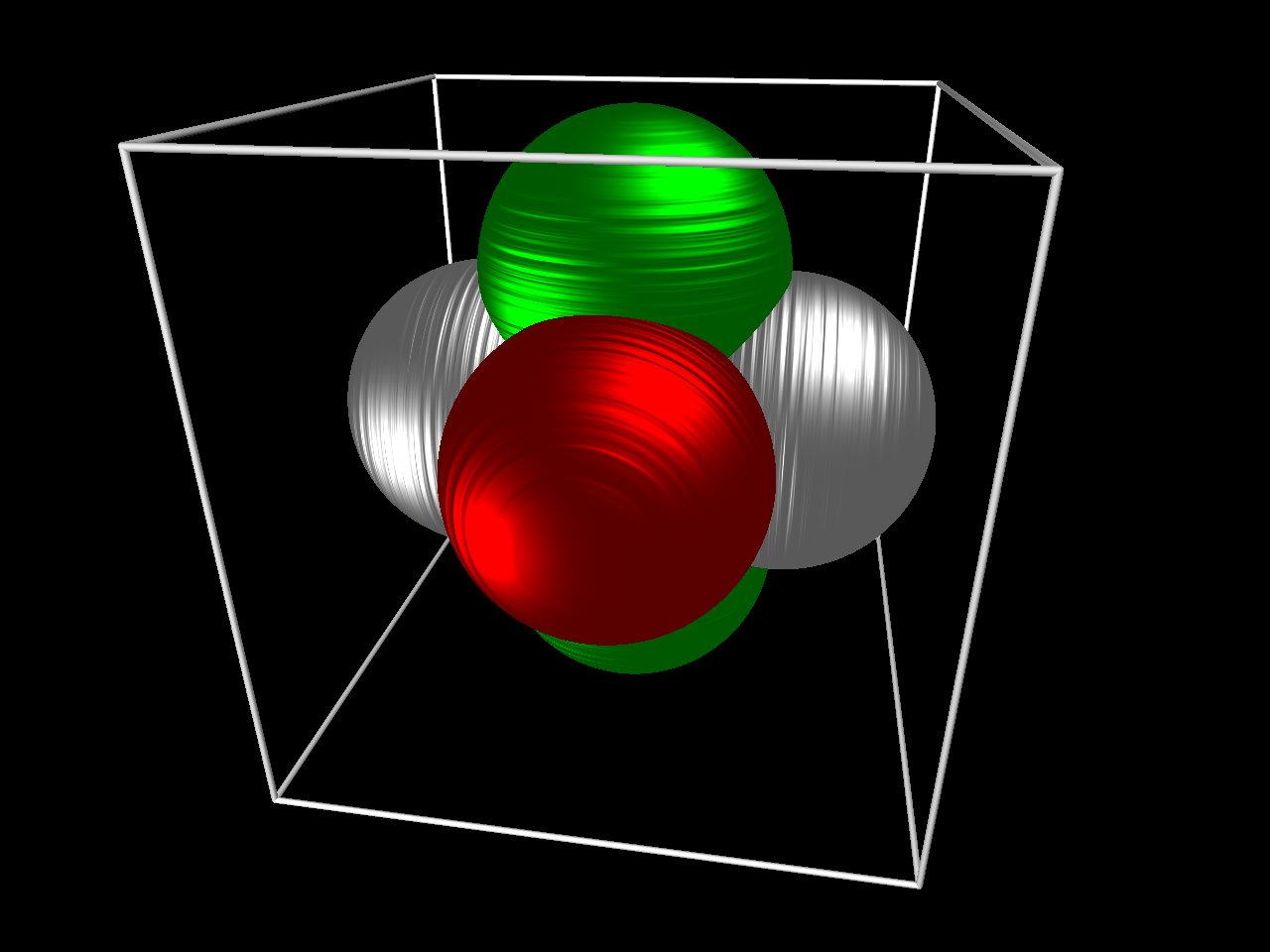}
\caption{The inner hyperplanes of the $24$-cell, viewed at infinity}
\label{24_cell figure 1}
\end{figure}

\begin{figure}[ht]
\centering
\includegraphics[scale=0.20]{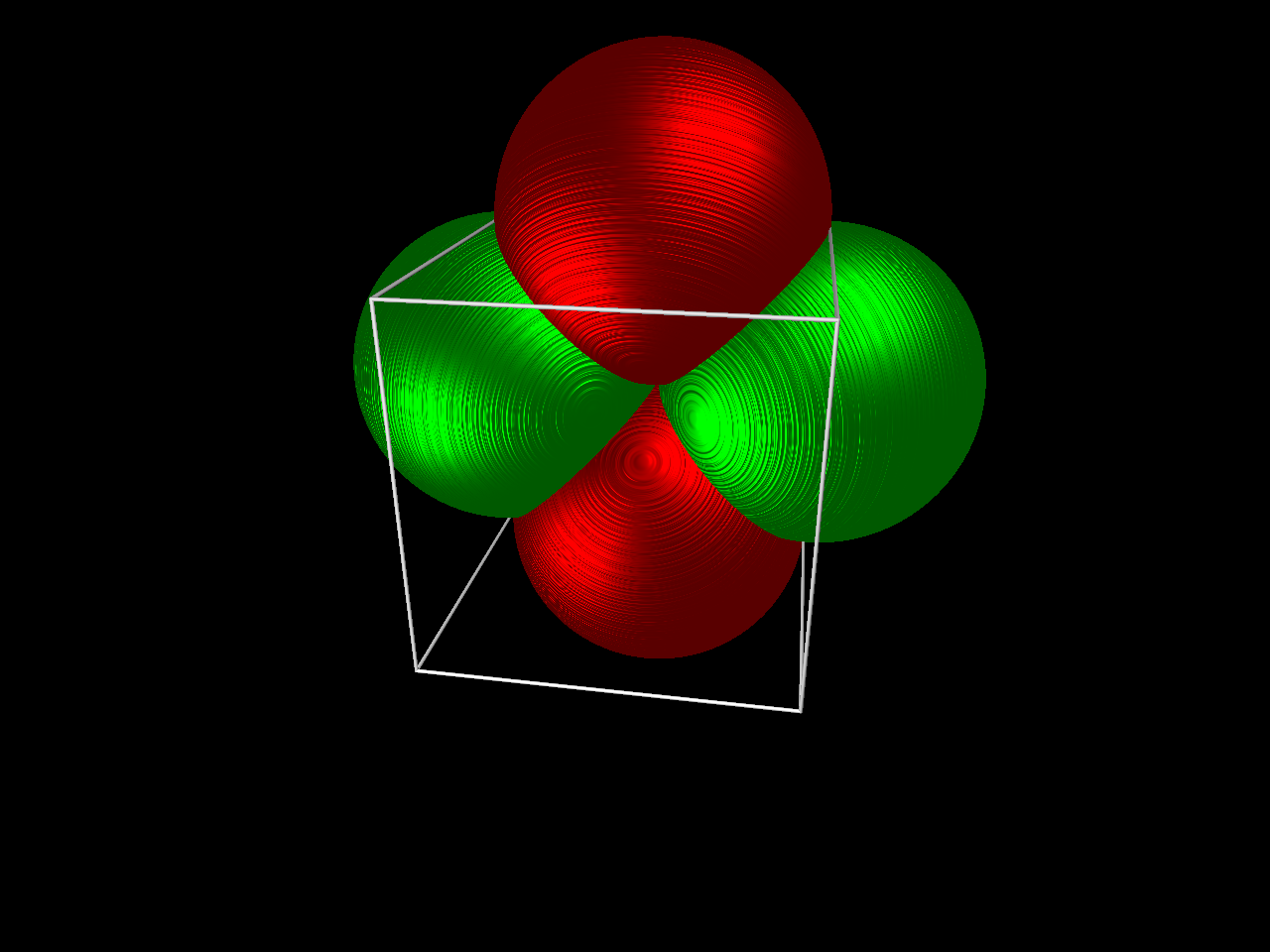}
\caption{The hyperplanes of the back face of the $24$-cell, viewed at infinity}
\label{24_cell figure 2}
\end{figure}

\begin{figure}[ht]
\centering
\includegraphics[scale=0.20]{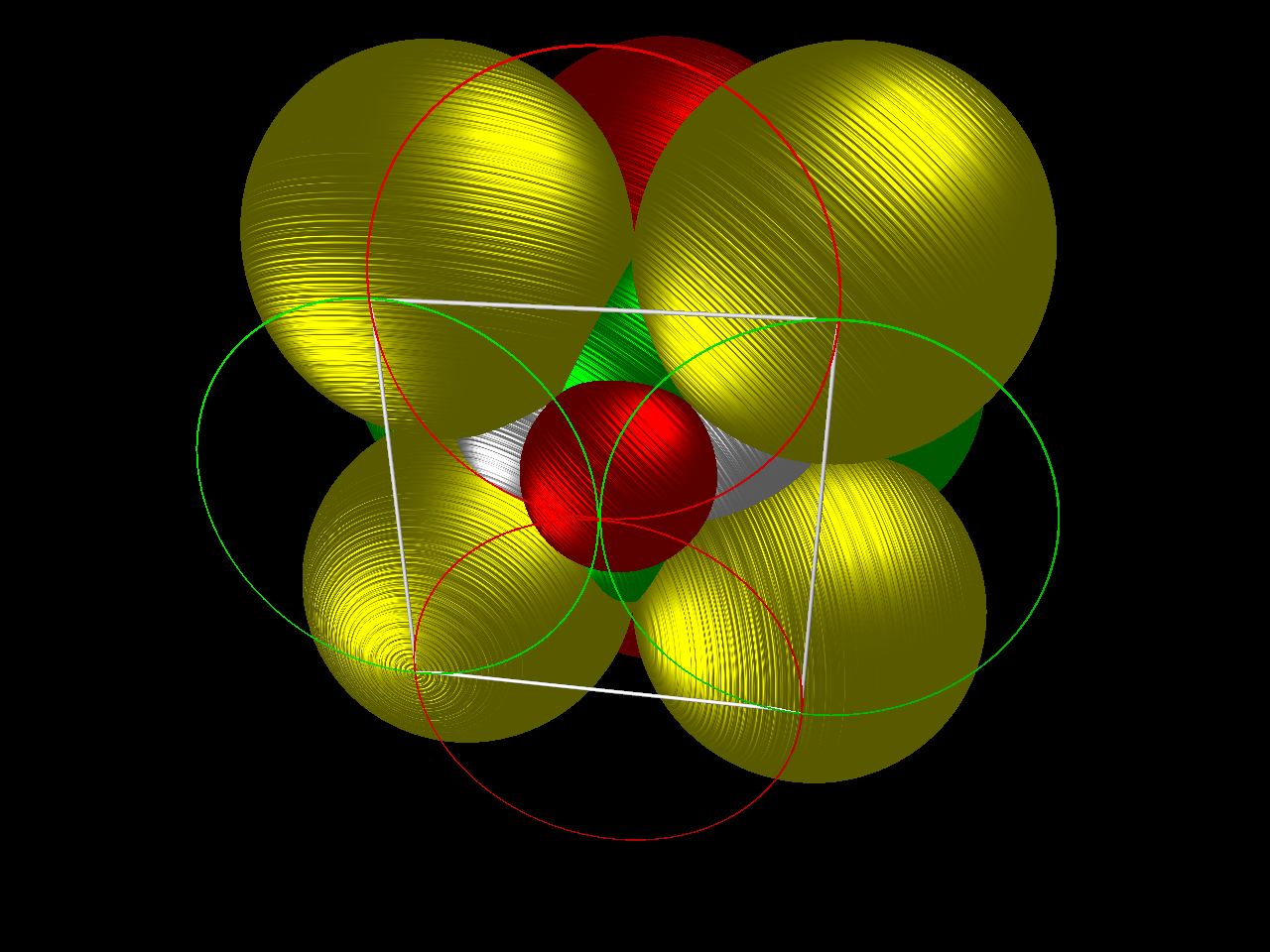}
\caption{The hyperplanes of the $24$-cell, viewed at infinity}
\label{24_cell figure 3}
\end{figure}

With these initial normalizations, the centers and radii for the spheres $S_i$ are listed in table \ref{24-cell center radii}.  We will later explain one method for deriving the information of table \ref{24-cell center radii}.  We encourage the reader to imagine each new sphere, and sketch a few pictures to get a feel for the arrangement.  The most difficult spheres to imagine are $S_{\m{6}}$, $S_{\p{6}}$, $S_{\p{7}}$, and $S_{\m{7}}$.  It is perhaps easier to mentally add them to the arrangement last.  To aid the readers imagination we have included figures \ref{24_cell figure 1}, \ref{24_cell figure 2}, and \ref{24_cell figure 3}, which we now explain.  

\begin{table}$
\begin{array}{ccc}
\text{sphere} & \text{       center} & \text{radius} \\
S_{\m{0}} & \text{ the } $xz$-\text{plane shifted by the vector }  (0,2,0)\\
S_{\p{0}} & \text{ the } $yz$-\text{plane}\\
S_{\p{3}} & \text{ the } $yz$-\text{plane shifted by the vector }  (2,0,0)\\ 
S_{\m{3}} & \text{ the } $xz$-\text{plane}\\ 
S_{\p{1}} & (1,2,0) & 1\\
S_{\m{1}} & (0,1,0) & 1\\
S_{\m{2}} & (2,1,0) & 1\\
S_{\p{2}} & (1,0,0) & 1\\
S_{\p{5}} & (1,2,2) & 1\\
S_{\m{5}} & (0,1,2) & 1\\
S_{\m{4}} & (2,1,2) & 1\\
S_{\p{4}} & (1,0,2) & 1\\
S_{\m{6}} & (1,3/2,1) & 1/2 \\
S_{\p{6}} & (1/2,1,1) & 1/2 \\
S_{\p{7}} & (3/2,1,1) & 1/2 \\
S_{\m{7}} & (1,1/2,1) & 1/2 \\
S_{\l{A}} & \text{ the } $xy$-\text{plane}\\
S_{\l{B}} & \text{ the } $xy$-\text{plane shifted by the vector }  (0,0,2)\\
S_{\l{C}} & (0,2,1) & 1\\
S_{\l{G}} & (2,2,1) & 1\\
S_{\l{F}} & (1,1,3/2) & 1/2\\
S_{\l{E}} & (1,1,1/2) & 1/2\\
S_{\l{D}} & (2,0,1) & 1\\
S_{\l{H}} & (0,0,1) & 1
\end{array}$
\caption{Centers and radii for the $24$ spheres of the $24$-cell}\label{24-cell center radii}
\end{table}

Note that not all $24$ hyperplanes are shown in these figures.  We have omitted several walls to make the images understandable.  To begin, figure \ref{24_cell figure 1} shows the edges of a cube and four spheres.  Each face of the implied cube is a sphere of table \ref{24-cell center radii}.  (Keep in mind that a plane counts as a sphere in this setting.)  Specifically, the spheres $S_{\m{0}}$, $S_{\p{0}}$, $S_{\p{3}}$, $S_{\m{3}}$, $S_{\l{A}}$, and $S_{\l{B}}$ form the faces of the cube.  The round spheres of figure \ref{24_cell figure 1} are the inner spheres, and they correspond to $S_{\m{6}}$, $S_{\p{6}}$, $S_{\p{7}}$, $S_{\m{7}}$, $S_{\l{F}}$, and $S_{\l{E}}$.  Notice that if we perform conformal inversion in the cube's inscribed sphere, then the faces of the cube are taken to the inner spheres and vice versa.

In figure \ref{24_cell figure 2} we have removed the inner spheres to show clearly the spheres of the back face.  The four spheres shown correspond to $S_{\m{1}}$, $S_{\m{5}}$, $S_{\l{C}}$, and $S_{\l{H}}$.  In the chosen coordinate system, the back face of the cube is the $yz$-plane.  The colored spheres are all cut by the $yz$-plane along equators.

Finally, in figure \ref{24_cell figure 3} we have attempted to show all $24$ spheres in a single picture.  To make the inner spheres visible, the spheres of the front face are not shown in their entirety.  The spheres of the front face are $S_{\m{2}}$, $S_{\m{4}}$, $S_{\l{G}}$, and $S_{\l{D}}$.  Figure \ref{24_cell figure 3} shows only an equator for each of these.  The yellow spheres correspond to $S_{\p{1}}$, $S_{\p{2}}$, $S_{\p{5}}$, and $S_{\p{4}}$.  Each yellow sphere is tangent to the planes $\{ x =0 \}$ and $\{ x=2\}$.  In this figure the many conformal symmetries of the $24$-cell are visible.

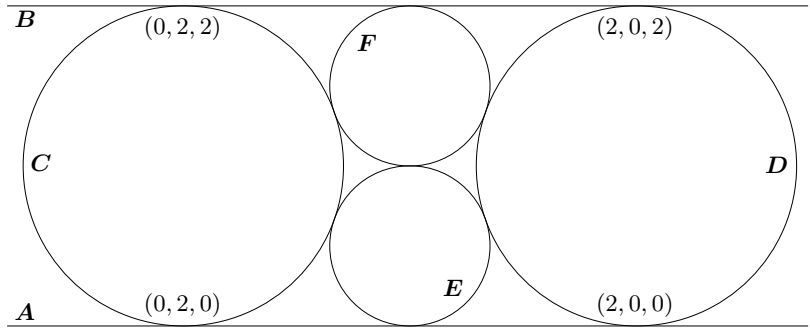
\begin{figure}[ht]
\labellist
\small\hair 2pt
\pinlabel $\l{A}$ [bl] at 3 2
\pinlabel $\l{B}$ [tl] at 3 201
\pinlabel $\l{C}$ [l] at 12 103
\pinlabel $\l{F}$ [tl] at 218 186
\pinlabel $\l{E}$ [br] at 290 17
\pinlabel $\l{D}$ [r] at 495 102
\pinlabel $(0,2,0)$ [b] at 111 3
\pinlabel $(0,2,2)$ [t] at 111 199
\pinlabel $(2,0,0)$ [b] at 395 3
\pinlabel $(2,0,2)$ [t] at 395 199
\endlabellist
\centering
\includegraphics[scale=0.6]{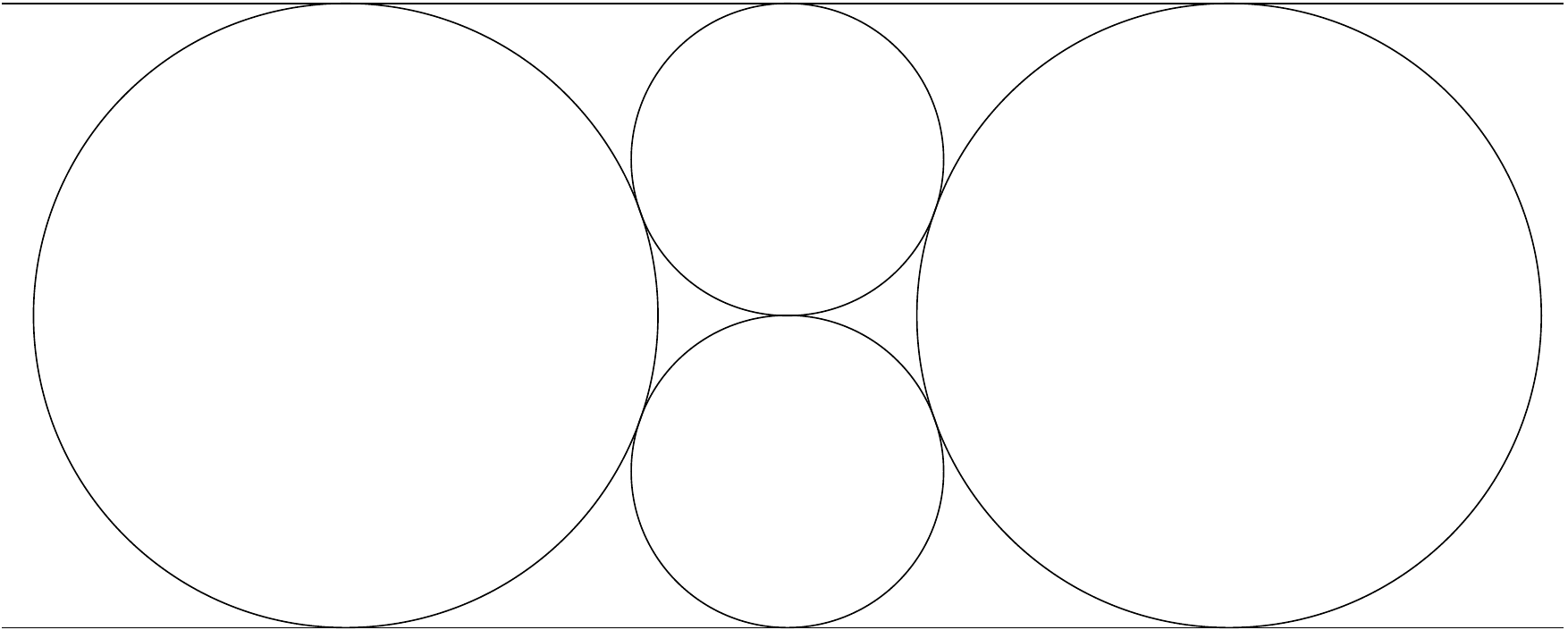}
\caption{A diagonal slice of the sphere arrangment}
\label{diag_slice}
\end{figure}

To clarify the arrangment of the inner spheres of figure \ref{24_cell figure 1}, we have also included figure \ref{diag_slice}.  Define the vertical plane $P:=\{x+y=2\}$.  Figure \ref{diag_slice} represents the intersections of $P$ with the spheres $S_{\l{A}}, S_{\l{B}}, S_{\l{C}}, S_{\l{F}}, S_{\l{E}},$ and $S_{\l{D}}$.

How could one determine this arrangement of $24$ spheres?  The first and best method is to begin with the familiar octahedral arrangement of $S_{\m{0}}$, $S_{\p{0}}$, $S_{\p{3}}$, $S_{\m{3}}$, $S_{\p{1}}$, $S_{\m{1}}$, $S_{\m{2}}$, $S_{\p{2}}$, and $S_{\l{A}}$, and one by one add the new spheres.  The combinatorics and geometry of the arrangement described by figure \ref{cube figure} will dictate the placement of each new sphere.  When doing this it is perhaps easier to finish by determining the placement of spheres $S_{\m{6}}$, $S_{\p{6}}$, $S_{\p{7}}$, and $S_{\m{7}}$.  

However, it is important to also have a more directly computation method.  This is particularly useful for determining and  drawing the arrangement of spheres once we deform things and there is less 
symmetry.  Let us sketch a way in which this can be done.

Let $W_i$ be a wall of $P_{24}$.  It is contained in a hyperplane of $\Hf$, which we denote by $H_i$.  We are interested in the intersection $\partial H_i = \overline{H_i} \cap \mathbb{S}^3_\infty$. 
$H_i$ is determined uniquely by the six ideal vertices of $P_{24}$ contained in $\partial H_i$.  Recall that the ideal vertices of $P_{24}$ are given by the $24$ light-like vectors $\sqrt{2} e_0 \pm e_i \pm e_j$ ($i,j \in \{1,2,3,4\}$).   The sphere at infinity can be identified with the set of geodesic rays based at any chosen point in $\Hf$ and this set of rays is naturally identified with the unit tangent space of $\Hf$ at that point.  If we choose the point to equal $e_0 = (1,0,0,0,0)$, then  the unit tangent space, $S_{e_0} \Hf$, is identified with the set of unit vectors in $\mathbb{R}^{1,4}$ whose first coordinate is zero.   This determines a map $\pi' : \mathbb{S}^3_\infty \rightarrow S_{e_0}\Hf$.  There is also a one-to-one correspondence that takes a geodesic ray based at $e_0$  to the unique light-like line in 
$\mathbb{R}^{1,4}$ to which it is asymptotic.  This provides an identification between the set of nonzero 
light-like vectors, up to scale, and the sphere at infinity.  With this identification the map $\pi'$ can be
described explicitly in coordinates.  It simply takes a light-like vector representing a point in $\mathbb{S}^3_\infty$, scales it to have $0\th$ coordinate $+1$, and then changes the $0\th$ coordinate to zero.

The six ideal vertices in $\partial H_i$ are exactly the six ideal vertices of $P_{24}$ orthogonal to the space-like vector $q_i \in \mathbb{R}^{1,4}$.  This is how we computed the
$q_i$ for table \ref{q table}.  For each $q_i$,  the conditions that a vector  be perpendicular to $q_i$, have Minkowski norm $0$ and be normalized to have $0\th$ coordinate $1$ defines a $2$-sphere.  Setting
the $0\th$ coordinate to be $0$ defines a sphere in $S_{e_0}\Hf$.  It is the image of 
$\overline{H_i} \cap \mathbb{S}^3_\infty$ under the map $\pi'$.

In order to place the arrangement of spheres in $\mathbb R^{3} \cup \infty$,  we use 
stereographic projection.  First we make the obvious identification of 
$S_{e_0} \Hf \subset \mathbb{R}^{1,4}$ with $\mathbb S^3 \subset \mathbb R^4$ by dropping the
$0\th$ coordinate.  Then we choose a point from which to define the stereographic projection.
This is an arbitrary choice.  To compute the explicit set of radii and centers of the spheres,
as well as to draw the pictures we have done the following:

First we apply to the entire configuration an element of $ \phi \in O(1,4)$ which is the identity in the last $3$ coordinates
and takes the space-like vector associated to the wall $W_{\l{A}}$ to the vector $(0,1,0,0,0)$.
In particular the image hyperplane $\phi(H_{\l{A}})$ consists of all points in $\Hf$ whose second  Minkowski coordinate is $0$.
(Keep in mind the coordinates are indexed from $0$ to $4$.)  We now identify the sphere at infinity with $\mathbb{S}^3$ as described above.  Under this identification the intersection of the closure of $\phi(H_{\l{A}})$ with the sphere at infinity is the $2$-sphere consisting of those vectors whose first coordinate is $0$.  We then take the standard
stereographic projection to $\mathbb{R}^3$ from the point $(0,1,0,0) \in \mathbb{S}^3$.  After a signed
permutation of the coordinates of $\mathbb{R}^3$,  the image of the $2$-sphere $S_{\l{A}}$ becomes the $x,y$-plane.  Finally, we use an affine map of $\mathbb{R}^3$ preserving the $x,y$-plane
to place everything in the form presented above.

 For completeness we record that this map  $\pi: \mathbb{S}^3_\infty \to  \mathbb R^{3} \cup \infty$ is expressed in coordinates by:
 
$$(0,x_1,x_2,x_3,x_4) \mapsto \left(
   \frac{\sqrt{2} - x_1 - x_2 - x_3 - x_4}{\sqrt{2} - x_1 - x_2} ,
  \frac{\sqrt{2} - x_1 - x_2 + x_3 - x_4}{\sqrt{2} - x_1 - x_2},
  \frac{\sqrt{2} - 2 x_1}{\sqrt{2} - x_1 - x_2}   \right).$$

One can explicitly determine the final configuration of $2$-spheres by applying the first map
$\phi$ above to all the space-like vectors $q_i$, using the image vectors to find spheres
in the sphere at infinity, and then mapping them to $\mathbb{R}^3$ using the above sequence of 
conformal maps.   However, it is simpler to apply $\pi$ to
the light-like vectors corresponding to the cusps of $P_{24}$ and then determine the spheres
in $\mathbb{R}^3$ from these.


\section{Describing the symmetries of the $24$-cell} \label{symmetries}

Recall that the hyperbolic $24$-cell, $P_{24}$, is the convex hull of the $24$ ideal points 
corresponding to the light-like vectors $\{ \sqrt{2} e_0 \pm e_i \pm e_j \}_{1 \le i < j \le 4}$ in $\mathbb{R}^{1,4}$.  It is well known that $P_{24}$ is a regular polytope \cite{Cox2}.  However, we will rederive this fact during the process of describing its symmetry group.

Begin with the obvious fact that the symmetry group of an ideal triangle has order $6$.  From this we deduce that the symmetry group of a regular, ideal octahedron has order $8 \cdot 6 = 48$.  Each wall of $P_{24}$ is an ideal octahedron and any symmetry that is the identity on one wall is the identity on all of $P_{24}$.  This shows that the symmetry group $S$ of the hyperbolic $24$-cell has order at most $24 \cdot 48 = 2^7 \cdot 3^2$.

For each $i = 1, \ldots, 4$, the map sending $e_i$ to $-e_i$ is an isometry of $P_{24}$.  The only point fixed by all these symmetries is $e_0$.  The barycenter of $P_{24}$ is also fixed by any symmetry, so it must equal $e_0$, implying $e_0$ is fixed by any symmetry.  Therefore the symmetry group of $P_{24}$ lies inside $\text{O}(4)$, the orthogonal group of the space-like coordinates.
The symmetries of $P_{24}$ can be identified with the isometries of $\mathbb{R}^4$ that permute the set of points $\{ \pm e_i \pm e_j \}_{1 \le i < j \le 4}$  in 
$\mathbb{R}^4$.  We will suppress the $0\th$ coordinate while freely identifying  $\mathbb{R}^4$ 
with the last $4$ coordinates of $\mathbb{R}^{1,4}$.  Notice that the set $\{ \pm e_i \pm e_j \}_{1 \le i < j \le 4}$ is exactly the set of integer points in $\mathbb{R}^4$ of norm $\sqrt{2}$.  Therefore the group $\text{O}(4;\mathbb{Z})$ is a subgroup of $S$.  An easy count shows that 
$$| \text{O}(4;\mathbb{Z}) | = 8 \cdot 6 \cdot 4 \cdot 2 = 2^7 \cdot 3,$$
implying $\text{O}(4;\mathbb{Z})$ has index at most $3$ inside $S$.

To prove $|S| = 2^7 \cdot 3^2$, and consequently that the $P_{24}$ is regular, it suffices to find an element in $S$ of order $3$ which is not in $\text{O}(4;\mathbb{Z})$.  Some experimenting produces the desired order $3$ matrix
\[ M = \frac{1}{2} \cdot \left( \begin{array}{rrrr}
  1 & 1 & 1 & -1 \\
  1 & 1 & -1 & 1 \\
  1 & -1 & 1 & 1 \\
  1 & -1 & -1 & -1 \end{array} \right) \]
It is easy to check that $M$ acts on the (last $4$ coordinates of the) vertex set of $P_{24}$ and is thus a symmetry.  Therefore $|S| \ge 2^7 \cdot 3^2$.  Combined with the previous upper bound we have proved that the hyperbolic $24$-cell is regular with symmetry group $S$ of order $2^7 \cdot 3^2$.

In search of a cleaner description of $S$, consider the group $\text{O}(4; \mathbb{Z}[{\textstyle \frac{1}{2}}])$.  Clearly $S \le \text{O}(4; \mathbb{Z}[{\textstyle\frac{1}{2}}])$.  It's not hard to see that an element of $\text{O}(4; \mathbb{Z}[{\textstyle\frac{1}{2}}]) - \text{O}(4;\mathbb{Z})$ must have all its matrix entries equal to $\pm \frac{1}{2}$.  From this observation, and the definition of the orthogonal group, one can count 
\[ \left| \text{O}(4; \mathbb{Z}[{\textstyle \frac{1}{2}}]) - \text{O}(4;\mathbb{Z}) \right| = 2^8 \cdot 3. \]
By recalling $| \text{O}(4;\mathbb{Z})|= 2^7 \cdot 3$ and summing we see that
\[ \left| \text{O}(4; \mathbb{Z}[{ \textstyle \frac{1}{2} }])\right| = 2^7 \cdot 3^2 = |S|,\]
yielding the desired clean description of the action on the last $4$ coordinates
\[ S = \text{O}(4; \mathbb{Z}[{ \textstyle \frac{1}{2}}]). \]
 
Recall that the walls of $P_{24}$ decompose naturally into three octets: the letter octet $\l{A}$ through $\l{H}$, the positive octet $\p{0}$ through $\p{7}$, and the negative octet $\m{0}$ through $\m{7}$.  The walls of the letter octet are associated to space-like vectors of the form  $\{ e_0 \pm \sqrt{2} e_i \}_{1\le i\le 4}$.  The walls of the positive and negative octets have associated vectors of the form 
$(\sqrt{2} e_0 \pm e_1 \pm e_2 \pm e_3 \pm e_4)$.  The positive octet consists of those with an \emph{even} number of negative coordinates, while the negative octet  consists of those with an \emph{odd} number of negative coordinates.  The octet to which a wall belongs consists of the six walls to which it is tangent and its opposite wall.  (Here, the opposite wall is obtained from the given wall by negating all its space-like coordinates.)  Since these properties are preserved by symmetries, the isometry group acts on the set of octets.  It is easy to check that this induces a surjective homomorphism onto $\Sigma_3$, the permutation group of order $6$.  The kernel $\U{K}$ of this homomorphism is $\text{SO}(4; \mathbb{Z}) < S$ acting on the last $4$ coordinates.  Thus this action produces the short exact sequence
\[ 1 \rightarrow \text{SO}(4; \mathbb{Z}) \rightarrow S \rightarrow \Sigma_3 \rightarrow 1. \]

Consider the group $\left( \mathbb{Z} /2 \mathbb{Z}\right)^4$ of sign changes to the coordinates of $\mathbb{R}^4$.  To avoid transposing the positive and negative octets, we restrict attention to the order $8$ subgroup $E$ which changes sign in an even number of coordinates.  Let $\Sigma_4$ act by permutation on the coordinates.  Then
$ E \rtimes \Sigma_4 \le \U{K}$.  As 
\[ | E \rtimes \Sigma_4 | = 2^6 \cdot 3 = | \U{K} |,\]
we know that $\U{K} \cong E \rtimes \Sigma_4$.

Soon we will be considering the $24$-cell with its walls $\l{G}$ and $\l{H}$ removed.  In preparation for this let us consider the subgroup $K$ of $\U{K}$ which either fixes or transposes walls $\l{G}$ and $\l{H}$.  As $\l{G}=(1,0,0,0,-\sqrt{2})$ and $\l{H}=(1,0,0,0,\sqrt{2})$,  $K$ is given by the elements of $\U{K}$ which fix or change the sign of the last coordinate of $\mathbb{R}^4$.  Therefore $K \cong E \rtimes \Sigma_3$, implying $|K| = 2^4 \cdot 3 = 48$.

For later use let us describe a generating set for $K$.  Let $\l{L}$ permute coordinates $1$ and $2$ 
of $\mathbb{R}^{1,4}$.  $\l{L}$ is reflection in the hyperplane orthogonal to $(0,-1, 1,0,0)$.  Let $\l{M}$ permute coordinates $2$ and $3$; it is reflection in the hyperplane orthogonal to $(0, 0,-1,1,0)$.  Let $\l{N}$ permute coordinates $2$ and $3$ while reversing their sign, reflection in the hyperplane orthogonal to $(0, 0,-1,-1,0)$.  Finally let $R$ change sign in the last $2$ coordinates.  $R$ is an angle $\pi$ rotation in the plane of the final $2$ coordinates.  Clearly $\langle \l{L},\l{M} \rangle = \Sigma_3$.  A little thought shows that $\langle \l{L},\l{M},\l{N},R \rangle = K$.  We will refer to $R$ as the roll symmetry.

Now, beginning with the hyperbolic $24$-cell, $P_{24}$, we remove walls $\l{G}$ and $\l{H}$.  Consider the group $\Gamma_{22}$ generated by reflections in the remaining $22$ walls.  To this group we may add the reflective generators $\l{L}, \l{M}, \l{N}$ because they are symmetries.  Call the resulting finite extension $\U{\Gamma}_{22}$.  The group $\langle \l{L},\l{M},\l{N}\rangle$ acts by conjugation on the $22$ generators of $\Gamma_{22}$.  Using this action one can conclude the set
\[ \left\{ \p{0},\m{0},\p{3},\m{3},\l{A},\l{L},\l{M},\l{N} \right\} \]
generates all of $\Gamma_{22} \cdot \langle \l{L},\l{M},\l{N} \rangle = \U{\Gamma}_{22}$.  With this relatively small generating set one can produce the following Coxeter diagram for $\U{\Gamma}_{22}$.  The roll symmetry $R$, which is not a reflection, is visible in this diagram as reflection in the horizontal line through $\l{L}$ and $\l{A}$.  

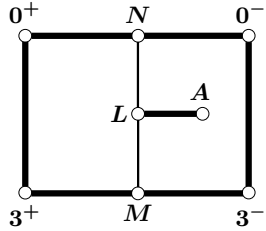
\begin{figure}[ht]
\labellist
\small\hair 2pt
\pinlabel $\p{0}$ [b] at 25 650
\pinlabel $\p{3}$ [t] at 25 4
\pinlabel $\l{M}$ [t] at 451 4
\pinlabel $\l{L}$ [r] at 435 327
\pinlabel $\l{A}$ [b] at 685 360
\pinlabel $\l{N}$ [b] at 451 650
\pinlabel $\m{0}$ [b] at 880 650
\pinlabel $\m{3}$ [t] at 880 4
\endlabellist
\centering
\includegraphics[scale=0.1]{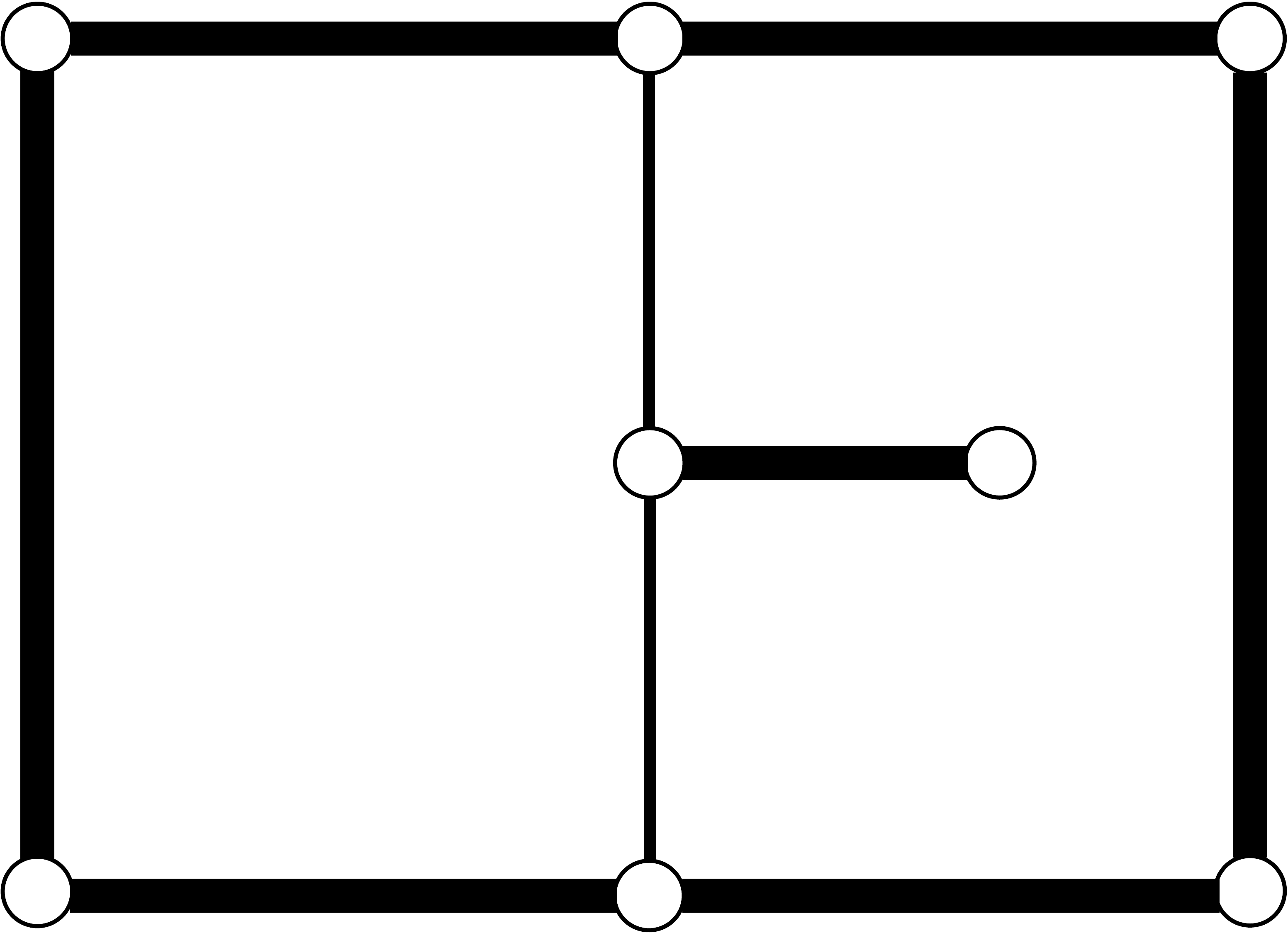}
\caption{A Coxeter diagram for the extension ${\widetilde{\Gamma}_{22}}$}
\label{coxeter_diagram1}
\end{figure}

Let us recall the meaning of the edges of figure \ref{coxeter_diagram1}.  The thick edges indicate walls which are tangent at infinity.  
The thin edges indicate walls which intersect at angle $\pi/3$.  Vertices not joined by an edge indicate walls which intersect orthogonally.



\section{Deformation preliminaries} \label{deformation prelims}

Having analyzed the geometry of the hyperbolic $24$-cell, $P_{24}$ in some detail, we are now in a position to begin the description of our family of  hyperbolic reflection groups.
As before, let $\Gamma_{24} < \text{Isom}(\mathbb{H}^4)$ be the hyperbolic reflection group in the walls of $P_{24}$.  Since the $24$-cell has finite volume, the group $\Gamma_{24}$ is a lattice, and is therefore locally rigid as a subgroup of $\text{Isom}(\mathbb{H}^4)$ by Garland-Raghunathan rigidity \cite{GR}.  In order to have a nontrivial deformation theory we will have to pass to an infinite index subgroup $\Gamma_{22}$ of $\Gamma_{24}$ obtained by removing two walls from $P_{24}$.  Before describing this process, we will first prove some preliminary results about the deformations of some simpler subgroups of $\Gamma_{24}$.  These will allow us to understand the local rigidity of
$\Gamma_{24}$ directly, and,  most importantly, understand which principles underlying local rigidity still persist even when there are nontrivial deformations.

First, consider the group $\Gamma_{\text{oct}} < \text{Isom}(\mathbb{H}^3)$ of reflections in a regular ideal right-angled hyperbolic octahedron $\Poct$.  This group is isomorphic to any of the $24$ subgroups of $\Gamma_{24}$ preserving a wall of the $P_{24}$.  Each ideal vertex of this octahedron 
is preserved by the subgroup generated by reflections in the $4$ faces which share that vertex.  The corresponding point at infinity is the unique point fixed by this entire subgroup;  we call such a group a {\em cusp group}.   It follows from Calabi-Weil rigidity that there are no nontrivial deformations of $\oct$ under which all of the vertex subgroups remain cusp groups.  Alternatively,  one can note that  the $8$ faces of $\Poct$ intersect $\partial \mathbb{H}^3 = \hat{\mathbb{C}}$ as shown in figure \ref{circles1 figure}, and then prove that  a continuous deformation in $\hat{\mathbb{C}}$ of the circles in that figure preserving orthogonal intersections and tangencies is necessarily a conformal motion of $\hat{\mathbb{C}}$.  

However, if the vertex subgroups are not required to remain cusp groups, there are nontrivial deformations.  In order to see the geometry of these deformations, we first analyze the deformations of a vertex subgroup.  It is isomorphic to the group of reflections in the sides of a Euclidean square, which equals the intersection of a neighborhood of a cusp of $\Poct$ with a horosphere based at the corresponding ideal vertex.  It can be presented in the following way: 
$$\Gamma_{\text{rect}} = \left\{r_i, ~ i=1,2,3,4 ~|~ r_i^2 = e = (r_i r_j)^2 ~~~ \text{for} ~~~ |i-j| \equiv 1\  (\text{mod }4) \right\}.$$ 
There is a $1$-parameter family of representations of this group as cusp groups, parametrized by the modulus of the rectangular, horospherical cross-section.  None of the deformations in this family extend to cusp preserving deformations of $\oct$.  However, there are other nearby representations of $\oct$
that are not cusp groups.  To see this, note that for nearby representations the order $2$ generators will be represented  by reflections in hyperplanes in $\mathbb{H}^3$ and the planes corresponding to neighboring edges of the rectangle must intersect in right angles.  But there is no requirement that the other pairs of planes, corresponding to opposite edges of the rectangle, remain tangent.

\begin{figure}[ht]
\labellist
\small\hair 2pt
\pinlabel $P_3$ [tr] at 64 703
\pinlabel $P_1$ [tl] at 541 696
\pinlabel $P_4$ [r] at 5 399
\pinlabel $P_2$ [r] at 156 399
\endlabellist
\centering
\includegraphics[scale=0.15]{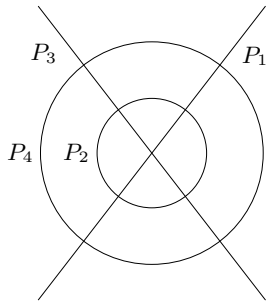}
\caption{A deformed cusp group}
\label{bullseye fig}
\end{figure}

Suppose one pair of opposite planes, $P_1$, $P_3$ intersect transversely.  The planes intersect the sphere at infinity in circles which intersect  in two points.  By applying an isometry of $\mathbb{H}^3$ we can take these points to be $0$ and $\infty$; i.e., the circles are two lines in $ \hat{\mathbb{C}}$ intersecting at the origin at some angle.  Any circle orthogonal to these two lines is centered at the origin and any pair of such circles are disjoint (or equal).  Thus, the remaining planes, $P_2$, $P_4$ must be a nonzero distance apart.  Similarly, if $P_1$, $P_3$ are disjoint and not tangent, their corresponding circles at infinity can be taken to be concentric circles centered at the origin, implying that $P_2$, $P_4$ intersect $\hat{\mathbb{C}}$ in two lines through the origin.  See figure \ref{bullseye fig} for a picture of this arrangement.  Such configurations are easily seen to be parametrized locally by the angle  of intersection between the two intersecting planes and the distance between the non-intersecting planes.  If $P_1$, $P_3$ intersect and $P_2$, $P_4$ are disjoint, then the products of the corresponding reflections, $r_1 r_3$ and $r_2 r_4$ are, respectively, a rotation by twice the angle of intersection and translation along the geodesic of intersection by twice the distance between $P_2$ and $P_4$.  From this it is not difficult to conclude the following.

\begin{lem} \label{rectangle deformations}
Let $\rho: \Gamma_{\text{rect}} \to \text{Isom}(\mathbb{H}^3)$ be a discrete faithful representation whose image is a cusp group.  Then the space of representations 
$\text{Hom}(\rect, \text{Isom}(\mathbb{H}^3))$ near $\rho$, up to conjugation,  
is $2$-dimensional.  In particular, there are nearby representations that are not cusp groups.  For such representations, one pair of opposite planes must intersect and the other must be a positive distance apart.   The angle of intersection and the distance are local parameters for such representations.
\end{lem}

\begin{figure}[ht]
\centering
\includegraphics[scale=0.1]{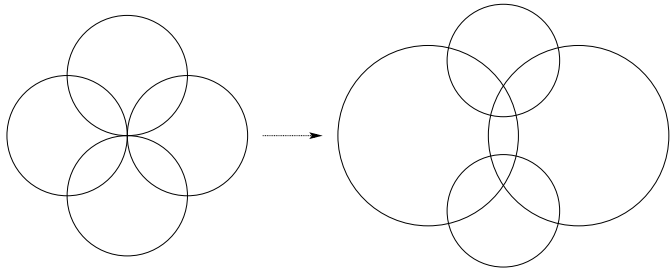}
\caption{A deformation of a cusp group}
\label{clover fig}
\end{figure}

\begin{proof}
The only thing that isn't clear from the previous discussion is how the non-cusp representations approximate the cusp representations.  This is best understood by examining figure \ref{clover fig}, where the picture on the left shows a cusp group with its parabolic fixed point at the origin, and the picture on the right shows a deformation of the cusp group.  
\end{proof}

The non-cusped deformations of the representations of the vertex subgroups of $\oct$ can give rise to 
deformations of the representation of $\oct$ itself.  These latter representations can be realized as
groups of reflections in the faces of families of polyhedra.   If a vertex group deforms
so that one pair of previously tangent planes intersect and the other pair move a positive distance
apart, this is realized in terms of polyhedra by a new edge appearing where there had been
an ideal vertex.  The edge is formed by the newly intersecting faces; it has two finite vertices
where these faces intersect with the two faces that are now a finite distance apart.  We will
see explicit examples of this phenomenon in later sections.  Figure \ref{deformedoctahedron  figure} 
shows a combinatorial
model of such a new polyhedron where two opposite vertices of the octahedron have been
deformed in this way.  The combinatorial choice of which pair of faces intersect can be interpreted as a choice of sign for the angle of intersection, where a positive sign is arbitrarily assigned to a given pair of faces at each vertex.
In this way, the parametrization can be seen to be smooth.


\begin{prop} \label{octahedral deformations}
Let $\rho: \oct \to \text{Isom}(\mathbb{H}^3)$ be the discrete faithful representation given by reflections in the faces of the ideal, regular right angled octahedron.  The space of representations
$\Hom(\oct, \text{Isom}(\mathbb{H}^3))$, up to conjugation,  near $\rho$ is $6$-dimensional. 
There is one parameter for each vertex of the octahedron, corresponding to the signed angle of intersection.
\end{prop}

\begin{proof}  The existence and parametrization by angle follows from the corresponding theorem for polyhedra due to Andreev \cite{Andreev} when all of the new dihedral angles are less than $\pi/2$.  It requires slightly more work to see how these bounded polyhedra approximate ones with ideal vertices, or, equivalently, how the corresponding groups of reflections approximate those with some cusp group representations for the vertex subgroups.

Note that these hyperbolic representations of $\oct$ need not be discrete or faithful.  \end{proof}

Although we have used the polyhedral theory of Andreev in explaining Proposition \ref{octahedral deformations},  there is another viewpoint which connects it more closely to the theory of hyperbolic $3$-manifolds.  We note that the vertex groups have an index $4$ subgroup isomorphic to $\Z \oplus \Z$.
A cusp representation of $\rect$ restricted to this subgroup is a parabolic representation fixing a
point at infinity; the quotient of any horosphere based at that point is a torus.  When the representation
of $\rect$ is deformed into a non-cusp representation, the corresponding representation of 
$\Z \oplus \Z$ preserves a geodesic in $\mathbb{H}^3$; one generator rotates around the
geodesic, the other translates along it.    

There is an index $4$ subgroup of $\oct$ in which all of the vertex subgroups lift to
torus subgroups;  the hyperbolic structure on the ideal octahedron lifts to a complete, finite volume hyperbolic structure on a $6$-component link \cite[Ch.6]{Th}.  Thurston's {\it hyperbolic Dehn filling theorem}
\cite[Ch.5]{Th} implies that, for any sufficiently small angles, one can always deform such a
hyperbolic structure so that each of the representations of the torus subgroups corresponding
to an end of the manifold is of the form above.  Topologically, this corresponds to truncating
the link complement so that it has $6$ torus boundary components and then attaching a solid torus so that the element corresponding to the rotation bounds a disk.  Geometrically there will be a singular hyperbolic structure on the resulting closed manifold with singularities along a core geodesic for each attached solid torus.

Proposition \ref{octahedral deformations} can be viewed as a special case of the hyperbolic
Dehn filling theorem, where the filling is done equivariantly with respect to this order
$4$ group.  Because of this we will describe such deformations of $\oct$ as {\it reflective hyperbolic Dehn filling}.

Secondly, we note that Proposition \ref{octahedral deformations} provides a $1$-dimensional family of deformations for each ideal vertex while Lemma \ref{rectangle deformations} gives a
$2$-dimensional family for  each vertex.  This is an explicit example of the half-dimensional phenomenon in the $3$-dimensional hyperbolic deformation theory that was mentioned in the introduction.  Precisely a half-dimensional subset of the deformations of the ($2$-dimensional) ``boundary groups" (vertex groups in this case) extend over the whole $3$-dimensional space.

The process of reflective hyperbolic Dehn filling does not have a precise analog in dimension $4$.
 We can explicitly see the reason for this in our particular situation by analyzing the deformation theory 
 of a vertex subgroup  of $\Gamma_{24}$.   Such a subgroup is isomorphic to the group of reflections in the faces of a Euclidean cube, $\cube$.  It has a presentation with $6$ generators, $\{r_i\}$, $i = 1,\ldots 6$, of order $2$
corresponding to the $6$ faces of the cube;  there are $12$ relations of the form $(r_i r_j)^2 = e$
corresponding to pairs of orthogonal faces of the cube.  There is no relation between the reflections corresponding to opposite faces.

In the discrete faithful representation of $\Gamma_{24}$ in $\text{Isom} (\Hf)$, the vertex subgroups
$\cube$ are represented by reflections in $3$-dimensional hypersurfaces that have exactly one point at infinity in common.  In other words they are cusp groups.  The key fact is that, unlike the situation in
dimension $3$ and the representations of $\rect$,
{\em all} nearby representations of $\cube$ are cusp groups.

\begin{lem} \label{cube deformations}  Let $\rho: \cube \to \text{Isom}(\Hf)$ be a discrete faithful representation as a cusp group.  Then the space $\Hom(\cube, \text{Isom}(\Hf))$, up to conjugation,
is $2$-dimensional near $\rho$.  All nearby representations are cusp groups.
\end{lem}

\begin{proof}
All representations of $\cube$ near $\rho$ will be generated by reflections in $6$  distinct 
$3$-dimensional hyperplanes (walls)  which we label by $X_i, Y_i, Z_i, i=1,2$.  Each wall
is required to intersect orthogonally those labelled by a distinct letter.  At $\rho$, or any other representation by a cusp group  the pairs labelled by the same letter are tangent; however,
this is not required by the representation.

There is a $2$-dimensional family of cusp representations, corresponding 
to the shape, up to scale, of the rectangular parallelepiped resulting from the intersection of reflection walls with a horosphere based at the fixed point of the cusp group.  We will show that these are the only representations of $\cube$ near $\rho$.

To see this, note that the walls $X_i, Y_i$ all intersect the $3$-dimensional hyperplane $Z_1$ orthogonally, and hence that the subgroup generated by reflections in these $4$ walls is conjugate to a representation of $\rect$ in $\Isom(\mathbb{H}^3)$.  Suppose that there is a representation near 
$\rho$ for which this subgroup is not a cusp group.  Then Lemma \ref{rectangle deformations} implies that, after possibly relabelling, the pair $X_1,X_2$ of walls will intersect while $Y_1, Y_2$
pull apart a nonzero distance.  Applying the same argument to $Y_i, Z_i$ implies that the
walls $Z_1, Z_2$ intersect.  However,  Lemma \ref{rectangle deformations}, applied to
the group generated by reflections in $X_i, Z_i$, says that if the $Z_i$ intersect, the
$X_i$ must pull apart, a contradiction.  Thus, all nearby representations of $\rho$ are
cusp groups.
\end{proof}

We can summarize Lemma \ref{cube deformations} as saying that locally ``cusp groups stay cusp groups".   This is the key step in the general local rigidity theorem of Garland-Raghunathan.
The subgroup of a cofinite volume discrete representation into $\Isom(\mathbb{H}^n)$,  
corresponding to an end of the quotient orbifold will be a cusp representation of an
$(n-1)$-dimensional Euclidean group.  They show that, when $n\geq 4$, all nearby representations 
of this Euclidean group into $\Isom(\mathbb{H}^n)$ will still be cusp groups.  It then
follows from Calabi-Weil local rigidity (or Mostow-Prasad global  rigidity) that there are
no nontrivial deformations.

In our situation, we can even avoid appealing to Calabi-Weil local rigidity to prove $\Gamma_{24}$ is locally rigid.  Lemma \ref{cube deformations} implies that  any continuous deformation of
$\Gamma_{24}$ takes its $24$ cusp subgroups to cusp subgroups of $\text{Isom}(\mathbb{H}^4)$.  Then Proposition \ref{octahedral deformations} implies that the $24$ subgroups of $\Gamma_{24}$, each of which preserves a wall of the $24$-cell,  are rigid under a continuous deformation.  From this the desired local rigidity of $\Gamma_{24}$ follows quickly.

\section{Defining the deformation} \label{defining the deformation 2}

With these preliminary observations concluded, we are now able to  better motivate our process
of finding a family of nontrivial deformations related to the inflexible group $\Gamma_{24}$.

In search of more flexibility it is natural to first consider a subgroup generated by reflections in $23$ of the $24$ walls of the hyperbolic $24$-cell.  For example we could omit the reflection in the 
letter wall $W_{\l{H}}$ and denote the resulting subgroup by $\Gamma_{23}$.  If we consider the 
subgroup of $\Gamma_{23}$ preserving the wall $W_{\l{G}}$, we see that it is isomorphic to 
$\oct$.  Furthermore, referring to figure \ref{cube figure} we see that all the vertex subgroups of this octahedral subgroup remain part of  subgroups of $\Gamma_{23}$ isomorphic to $\cube$, and thus,
by Lemma \ref{cube deformations}, they remain cusp groups.  It follows that the octahedral  subgroup
preserving wall $W_{\l{G}}$ is locally rigid.  This implies that, after conjugation, the walls corresponding to $\l{G}, \m{0}, \p{1}, \m{2}, \p{3}, \m{4}, \p{5}, \m{6}, \p{7}$ remain unchanged.  Repeating the above
argument about vertex subgroups for the eight numbered walls in this list, we see that the
subgroups preserving them are also locally rigid.  It is then not difficult to conclude that
$\Gamma_{23}$  is also infinitesimally rigid despite the fact that its covolume is clearly infinite.  (A variant of this minimally computational approach to the infinitesimal rigidity of $\Gamma_{23}$ was explained to us by D. Allcock after we first observed it through a brute force computation.)

Since the above local rigidity argument relied on the local rigidity of the octahedral subgroup preserving to wall $W_{\l{G}}$, we are led to remove that wall as well.  Thus, we continue on to the subgroup $\Gamma_{22}< \Gamma_{24}$ generated by reflections in all the walls other than $W_{\l{G}}$ and $W_{\l{H}}$.  The hyperbolic orbifold $\mathbb{H}^4 / \Gamma_{22}$ has infinite volume.  It is obtained by taking the hyperbolic polytope $P_{24}$ and attaching a $4$-dimensional Fuchsian end to its $\l{G}$ and $\l{H}$ walls (see Section \ref{def:orbifold}).  The convex core of $\Gamma_{22}$ has totally geodesic boundary.  It is the orbifold with boundary obtained from $P_{24}$ by mirroring all its walls other than $W_{\l{G}}$ and $W_{\l{H}}$.  The $\l{G}$ and $\l{H}$ walls remain unmirrored, forming the boundary of the convex core.

The $24$-cell is self-dual, and therefore has $24$ cusps.  Passing from $\Gamma_{24}$ to $\Gamma_{22}$ opens up exactly $12$ of these cusps.  Specifically, all $24$ cusps of $\Gamma_{24}$ have rank $3$ (i.e. they are virtually a rank $3$ free abelian group).  After throwing away the wall $W_{\l{G}}$, the six cusps involving the wall $W_{\l{G}}$ become reflections in $5$ sides of a cube, which is virtually a rank $2$ abelian group.  Similarly, after throwing away the wall $W_{\l{H}}$, six more cusps become rank $2$.  It's easy to check that no cusp involves both $W_{\l{G}}$ and $W_{\l{H}}$.  $\Gamma_{22}$ is left with $12$ cusps of rank $3$, and $12$ cusps of rank $2$.  For example, the cusp of $\Gamma_{24}$ formed by the walls corresponding to $\l{A}$, $\l{B}$, $\p{0}$, $\m{0}$, $\p{3}$, and $\m{3}$ remains rank $3$ in $\Gamma_{22}$.  The cusp of $\Gamma_{24}$ formed by the walls corresponding to $\l{A}$, $\l{G}$, $\m{0}$, $\p{1}$, $\m{2}$, and $\p{3}$ opens up to a rank $2$ cusp in $\Gamma_{22}$.  By Lemma \ref{cube deformations}, under any deformation of $\Gamma_{22}$ the rank $3$ cusps will remain rank $3$ cusps, at least near the initial group $\Gamma_{22}$.  As we will soon see, the same is not true for the rank $2$ cusps.

Although we can no longer repeat the argument used for the local rigidity of $\Gamma_{23}$,  it is not clear whether or not  $\Gamma_{22}$ is again infinitesimally rigid as a representation in $\Hom (\Gamma_{22}, G)$.  We can quickly do a dimension count.  The space of arrangements of $22$ hyperplanes in $\mathbb{H}^4$ has dimension $88$.  There are $80$ orthogonality conditions, each assumed to kill $1$ dimension.  This leaves $8$ dimensions.  After subtracting $10$ for the action of $G$ we are left with a $(-2)$-dimensional deformation space.  From this we may naively expect $\Gamma_{22}$ to be infinitesimally rigid.  This simple dimension count yields the wrong answer in this case;  the group $\Gamma_{22}$ will turn out to have nontrivial deformations.  

A natural way to start the search is by looking for deformations of $\Gamma_{22}$ which preserve its symmetries.  In Section \ref{computing infinitesimally} we will prove that all nearby representations of $\Gamma_{22}$ retain these symmetries.  However, in the interest of getting more quickly to the family of examples, we will assume that the symmetries are preserved.  By narrowing the possibilities in this way it is much easier to find a nice deformation.

We first consider the possible configurations of the letter walls.  The first important observation is
that any deformation of $\Gamma_{22}$ will preserve the tangencies between the remaining $6$
letter walls.  This is because the cusp groups to which they belong in $\Gamma_{24}$ remain rank 
$3$ even after the removal of walls $W_{\l{G}}$ and $W_{\l{H}}$.  By Lemma \ref{cube deformations} these subgroups will remain cusp groups and walls that were tangent remain so.  It turns out that, if we further assume that the deformation preserves the symmetries of $\Gamma_{22}$,  the entire configuration of the letter walls is locally rigid.  Thus, it is quite remarkable that, even with this rigid configuration 
of walls, we will still be able to deform our polyhedron and its corresponding reflection group.

From the initial list of space-like vectors (Table \ref{q  table}), we consider those corresponding to the six remaining walls from the letter octet.  
\[
\begin{array}{lll}
\l{A} = \left( 1, \sqrt{2},0,0,0 \right) & \qquad & \l{B} = \left(1,0,\sqrt{2},0,0\right) \\
\l{C} = \left( 1,0,0,\sqrt{2},0 \right) & \qquad & \l{D} = \left(1,0,0,-\sqrt{2},0\right) \\
\l{E} = \left( 1,0,-\sqrt{2},0,0\right) & \qquad & \l{F} = \left(1,-\sqrt{2},0,0,0\right)
\end{array}
\]

Suppose that $\l{A}(t)$ through $\l{F}(t)$ are paths of space-like vectors defining one-parameter deformations of walls $W_{\l{A}}$ through $W_{\l{F}}$, beginning at $t=1$.  Let us assume that $\l{A}(1)=\l{A}$, $\l{B}(1) = \l{B}$, etc..  We assume that the deformation preserves the tangencies of the original walls.  
Furthermore, we restrict attention to symmetric deformations of the letter walls.  In fact, we will only require that the deformation respect  the roll symmetry $R$, defined near the end of Section \ref{symmetries}.  The action of $R$ on the letters is given by preserving $\l{A}$, $\l{B}$, $\l{E}$, and $\l{F}$, while transposing $\l{C}$ and $\l{D}$.  (Note that $R$ does not fix any wall pointwise.)  Hence, we will assume the existence of a path $R(t)$ of hyperbolic involutions satisfying
\[
\begin{array}{lll} 
R(t) \cdot \l{A}(t) = \l{A}(t) & \qquad & R(t) \cdot \l{B}(t) = \l{B}(t) \\
R(t) \cdot \l{C}(t) = \l{D}(t) & \qquad & R(t) \cdot \l{D}(t) = \l{C}(t) \\
R(t) \cdot \l{E}(t) = \l{E}(t) & \qquad & R(t) \cdot \l{F}(t) = \l{F}(t) 
\end{array}
\]
for all $t$.

It is easiest to describe this in $\partial \mathbb{H}^4$, where the $1$-parameter families of hyperplanes, $H_{\l{A}(t)}$-$H_{\l{F}(t)}$, become the $1$-parameter families of spheres $S_{\l{A}(t)}$-$S_{\l{F}(t)}$, and Minkowski isometries become conformal automorphisms of $\partial \mathbb{H}^4$.  Recall that, at time $t=1$,  the roll isometry $R(1)$ changes the sign of the final two coordinates of $\mathbb{R}^{1,4}$.   At  the sphere at infinity $\partial \mathbb{H}^4$ it  is rotation by $\pi$ with axis $\alpha$ orthogonally intersecting the spheres $S_{\l{A}(1)}$, $S_{\l{B}(1)}$, $S_{\l{E}(1)}$, and $S_{\l{F}(1)}$, passing through their points of tangency $S_{\l{A}(1)} \cap S_{\l{E}(1)}$, $S_{\l{E}(1)} \cap S_{\l{F}(1)}$, $S_{\l{F}(1)} \cap S_{\l{B}}(1)$, and $S_{\l{B}(1)} \cap S_{\l{A}(1)}$.  On these four spheres the roll acts by rotation through angle $\pi$.  By assumption, $R(t)$ preserves each of the four spheres, $S_{\l{A}(t)}$, $S_{\l{B}(t)}$, $S_{\l{E}(t)}$, and $S_{\l{F}(t)}$, for all $t$.  Since the spheres are tangent in the same combinatorial pattern for all $t$, the rotational axis for $R(t)$ must always be orthogonal to the spheres and pass through their various points of tangency.  

At this point we apply  a path of Minkowski isometries to put our paths $\l{A}(t)$ through $ \l{F}(t)$  into a simple form.  Identifying  $\partial \mathbb{H}^4$ with $\mathbb{R}^3 \cup \infty$,  we
can assume that $S_{\l{A}(t)}$ and $S_{\l{B}(t)}$ are parallel to the $xy$-plane at heights
$z = \pm 1$, respectively.  They are tangent at infinity so the axis of $R(t)$ is always a vertical line.  It can be assumed to be the $z$-axis for all $t$ so $R(t) = R$ will be constant.  The points of tangency  $S_{\l{A}(t)} \cap S_{\l{E}(t)}$ and $S_{\l{F}(t)} \cap S_{\l{B}}(t)$  will then be $(0,0, \pm 1)$, respectively, for all $t$.
Since $S_{\l{C}(t)}$ and $S_{\l{D}(t)}$ are tangent to both planes $z = \pm 1$ they must be spheres of radius $1$;  since $R$ transposes $S_{\l{C}(t)}$ and $S_{\l{D}(t)}$ they must be symmetric with respect to the $z$-axis.  After a final isometry, their points of tangency with the planes
can be put in the form $(-x(t), 0, \pm 1)$ and $(x(t), 0, \pm 1)$.   This implies that the $yz$-plane
is perpendicular to $S_{\l{C}(t)}$, $S_{\l{D}(t)}$, $S_{\l{E}(t)}$, and $S_{\l{F}(t)}$, for all $t$.  It is
then easy to see that the only way that these spheres can be tangent in the required way is
for $x(t) = 2$ and for $S_{\l{E}(t)}$ and $S_{\l{F}(t)}$ to be tangent at the origin for all $t$.  The pattern of circles given by the intersection with the $yz$-plane is shown in figure \ref{diag_slice} (where a slightly different normalization was used).

We have therefore proved the following proposition.

\begin{prop} \label{letter walls prop}
Suppose that $\l{A}(t)$ through $\l{F}(t)$ are paths of Minkowski-unit vectors defining one-parameter deformations of walls $W_{\l{A}}$ through $W_{\l{F}}$ such that $\l{A}(1)=\l{A}$, $\l{B}(1) = \l{B}$, etc..  Assume the deformation preserves the tangencies between the walls present when $t=1$.  Moreover, assume there exists a path of hyperbolic involutions $R(t)$ such that $R(1)$ is the roll symmetry, and $R(t)$ induces the same permutation on the letters as $R$.  Then the deformation is trivial, i.e. it is conjugate to the constant deformation.
\end{prop}

Thus we see that any symmetry preserving deformation of $\Gamma_{22}$ will, up to conjugation,  fix the space-like vectors corresponding to the letter walls.  So, we will assume that these vectors are fixed exactly and try to find space-like vectors for the remaining walls that satisfy the required orthogonality conditions.

Begin with the vector $\p{0} = \left( \sqrt{2},1,1,1,1\right)$.  $\p{0}$ is orthogonal to $\l{A}$, $\l{B}$, and $\l{C}$.  If we extend the vector to a deformation $\p{0}(t)$, it must remain orthogonal to $\l{A}$, $\l{B}$, and $\l{C}$ for all $t$.  This implies that, up to scale,
\[ \p{0}(t) =  \left( \sqrt{2}, 1,1,1, a^+ (t) \right).\]
A similar analysis for $\m{0}$ shows that
\[ \m{0}(t) = \left( \sqrt{2}, 1,1,1, a^- (t) \right).\]
Finally, the vectors $\p{0}(t)$ and $\m{0}(t)$ must remain orthogonal for all $t$.  Therefore $a^+ (t) \cdot a^-(t) = -1$.  From this we can conclude that
\[ \p{0}(t)=\left(\sqrt{2},1,1,1,1/t\right) \quad \text{and  } 
    \m{0}(t)=\left( \sqrt{2},1,1,1,-t \right) \]
are the unique paths satisfying the required orthogonality conditions, up to scaling and reparametrization.

Repeating this analysis for the $14$ other numbered space-like vectors forces us to the list of paths of space-like vectors found in table \ref{qt table}.  These paths are similarly unique up to scaling and reparametrization.  

\begin{table}
\begin{eqnarray*}
\p{0} = \left( \sqrt{2},1,1,1,1/t \right) , & &
\m{0} = \left( \sqrt{2},1,1,1,-t \right),\\
\p{1} = \left( \sqrt{2},1,-1,1,-1/t\right),& &
\m{1} = \left( \sqrt{2},1,-1,1,t\right),\\
\p{2} = \left( \sqrt{2},1,-1,-1,1/t\right),& &
\m{2} = \left( \sqrt{2},1,-1,-1,-t\right),\\
\p{3} = \left( \sqrt{2},1,1,-1,-1/t\right),& &
\m{3} = \left( \sqrt{2},1,1,-1,t\right),\\
\p{4} =\left(\sqrt{2},-1,1,-1,1/t\right),& &
\m{4} =\left(\sqrt{2},-1,1,-1,-t\right),\\
\p{5} =\left(\sqrt{2},-1,1,1,-1/t\right),& &
\m{5} = \left( \sqrt{2},-1,1,1,t\right),\\
\p{6} =\left(\sqrt{2},-1,-1,1,1/t\right),& &
\m{6} =\left(\sqrt{2},-1,-1,1,-t\right),\\
\p{7} =\left(\sqrt{2},-1,-1,-1,-1/t\right),& &
\m{7} =\left(\sqrt{2},-1,-1,-1,t\right),\\
\l{A} =\left(1,\sqrt{2},0,0,0\right),& &
\l{B} =\left(1,0,\sqrt{2},0,0\right),\\
\l{C} =\left(1,0,0,\sqrt{2},0\right),& &
\l{D} =\left(1,0,0,-\sqrt{2},0\right),\\
\l{E} =\left(1,0,-\sqrt{2},0,0\right),& &
\l{F} =\left(1,-\sqrt{2},0,0,0\right)
\end{eqnarray*}
\caption{A deformation of the group $\Gamma_{22}$}\label{qt table}
\end{table}

The computation showing that the original $80$ orthogonality relations between the $22$ generators of $\Gamma_{22}$ are valid for all $t$ is essentially the same as that done in Section \ref{24-cell} for $\Gamma_{24}$, the reflection group of the hyperbolic $24$-cell.  The
coefficients of the numbered space-like vectors only change in the final coordinate;  that coordinate
is $\pm 1/t$ for all the positive vectors and $\pm t$ for all the negative vectors with the same
signs as before.  Thus the contribution to the dot product of a positive and a negative vector will
still be $\pm 1$, the signs remaining the same.  Those that were orthogonal will remain so.
The letter vectors remain exactly the same and have coefficient $0$ in the final coordinate
so their dot products with the numbered vectors are exactly the same.  

Note, however, that
the norms of the numbered vectors, the dot products among the positive 
vectors, and the dot products among the negative vectors are all changing.  So the relationships between the non-orthogonal walls will vary.   In particular, some pairs of previously tangent positive walls intersect and some previously tangent negative walls pull apart.  The pairs behaving in this way are those 
which are diagonal across a face of the inner or the outer cube in figure \ref{cube figure}.  The remaining pairs of previously tangent numbered walls and all the previously tangent letter walls remain tangent.  Those walls that are a positive distance apart at $t=1$ remain so for nearby values of $t$, indeed up to $t= \sqrt{1/7}$.
(There is no particular significance to this value of $t$ as the new intersections of walls occurring at that time take place
away from the polytope itself.  This will become clearer once we have discussed the geometry of the extended group $\U{\Gamma}_{22}$ in more detail.)

These facts can be easily derived using table \ref{qt table} and formulae \ref{angle equation} and 
\ref{distance equation},  using symmetry to reduce the number of computations required.  For later reference, we record the results of these computations here:

\begin{prop} \label{angle prop}
Suppose $\sqrt{1/7} < t < 1$ and consider any pair of the $22$ hyperplanes of $\Hf$ determined by the space-like vectors $\{ \p{0}, \m{0}, \ldots, \m{7}, \l{A}, \ldots , \l{F} \}$.  
\begin{enumerate}
\item If they are separated by a strictly positive distance in $P_{24}$ then they are still separated by a strictly positive distance.
\item If they intersect orthogonally in $P_{24}$ then they still intersect orthogonally.
\item If they are tangent at infinity in $P_{24}$ then they are either tangent at infinity, separated by a distance $\ell$, or intersect at an angle $\theta$, where
$$\cos (\theta) =  \frac{ 3t^2 - 1}{1+t^2} \quad \text{and} \quad 
      \cosh (\ell) = \frac{ 3 - t^2 }{1+ t^2} .$$
Explicitly, the following pairs intersect at angle $\theta$: 
$$\{\p{1},\p{3}\},\{\p{3},\p{5}\},\{\p{5},\p{7}\},\{\p{7},\p{1}\},\{\p{1},\p{5}\},\{\p{3},\p{7}\},$$
$$\{\p{2},\p{0}\},\{\p{0},\p{4}\}, \{\p{4},\p{6}\},\{\p{6},\p{2}\},\{\p{2},\p{4}\},\{\p{0},\p{6}\}.$$  
The following pairs are separated by a distance $\ell$: 
$$\{\m{2},\m{0}\},\{\m{0},\m{4}\},\{\m{4},\m{6}\},\{\m{6},\m{2}\},\{\m{2},\m{4}\},\{\m{0},\m{6}\},$$
$$\{\m{1},\m{3}\},\{\m{3},\m{5}\},\{\m{5},\m{7}\},\{\m{7},\m{1}\},\{\m{1},\m{5}\},\{\m{3},\m{7}\}.$$  
Finally, any pair which is tangent at infinity in $P_{24}$ and is not contained in either of the previous two lists, remains tangent at infinity.
\end{enumerate}
\end{prop}

With table \ref{qt table} we can define the path of homomorphisms
\[ \rho_t : \Gamma_{22} \rightarrow \text{Isom}(\mathbb{H}^4) \]
in the obvious way:  $\rho_t$ takes reflection in the wall $W_{\l{V}}$ to reflection in the wall $W_{\l{V}}(t)$ orthogonal to $\l{V} (t)$, where $\l{V}$ is a symbol from the list $\left\{ \p{0},\m{0}, \ldots, \m{7}, \l{A}, \ldots, \l{F} \right\}$.  The representation $\rho_t$ is defined for all $t>0$, and $\rho_1$ is the inclusion map.  Let $\sigma : \Gamma_{22} \rightarrow \Gamma_{22}$ be the automorphism described by changing the sign of the third spatial coordinate of the vectors $\left\{ \p{0},\m{0}, \ldots, \m{7}, \l{A}, \ldots, \l{F} \right\}$.  The deformation $\rho_t$ has the symmetry $\rho_t \circ \sigma = \rho_{1/t}$.  For this reason we will restrict our attention to values of $t \le 1$.  

Recall the symmetries $\l{L}$, $\l{M}$, and $\l{N}$ of $\Gamma_{22}$: $\l{L}$ permutes the first two spatial coordinates, $\l{M}$ permutes the middle pair of spatial coordinates, and $\l{N}$ permutes the middle pair of spatial coordinates while changing their signs.  A quick inspection of table \ref{qt table} reveals that these reflections remain symmetries of $\rho_t (\Gamma_{22})$ for all $t$.  This means our deformation extends to a deformation $\U{\rho}_t$ of the extended reflection group 
\[ \U{\Gamma}_{22} := \langle \Gamma_{22}, \l{L}, \l{M}, \l{N} \rangle, \]
where $\U{\rho}_t$ fixes $\l{L}$, $\l{M}$, and $\l{N}$, and equals $\rho_t$ on $\Gamma_{22}$.  For the sake of notational simplicity, let us denote this extended representation $\U{\rho}$ by simply $\rho$.

\begin{figure}[ht]
\labellist
\small\hair 2pt
\pinlabel $\p{0}$ [b] at 25 650
\pinlabel $\p{3}$ [t] at 25 4
\pinlabel $\l{M}$ [t] at 451 4
\pinlabel $\l{L}$ [r] at 435 327
\pinlabel $\l{A}$ [b] at 685 360
\pinlabel $\l{N}$ [b] at 451 650
\pinlabel $\m{0}$ [b] at 880 650
\pinlabel $\m{3}$ [t] at 880 4
\pinlabel $\nu$ [b] at 235 632
\pinlabel $\nu$ [b] at 235 35
\endlabellist
\centering
\includegraphics[scale=0.1]{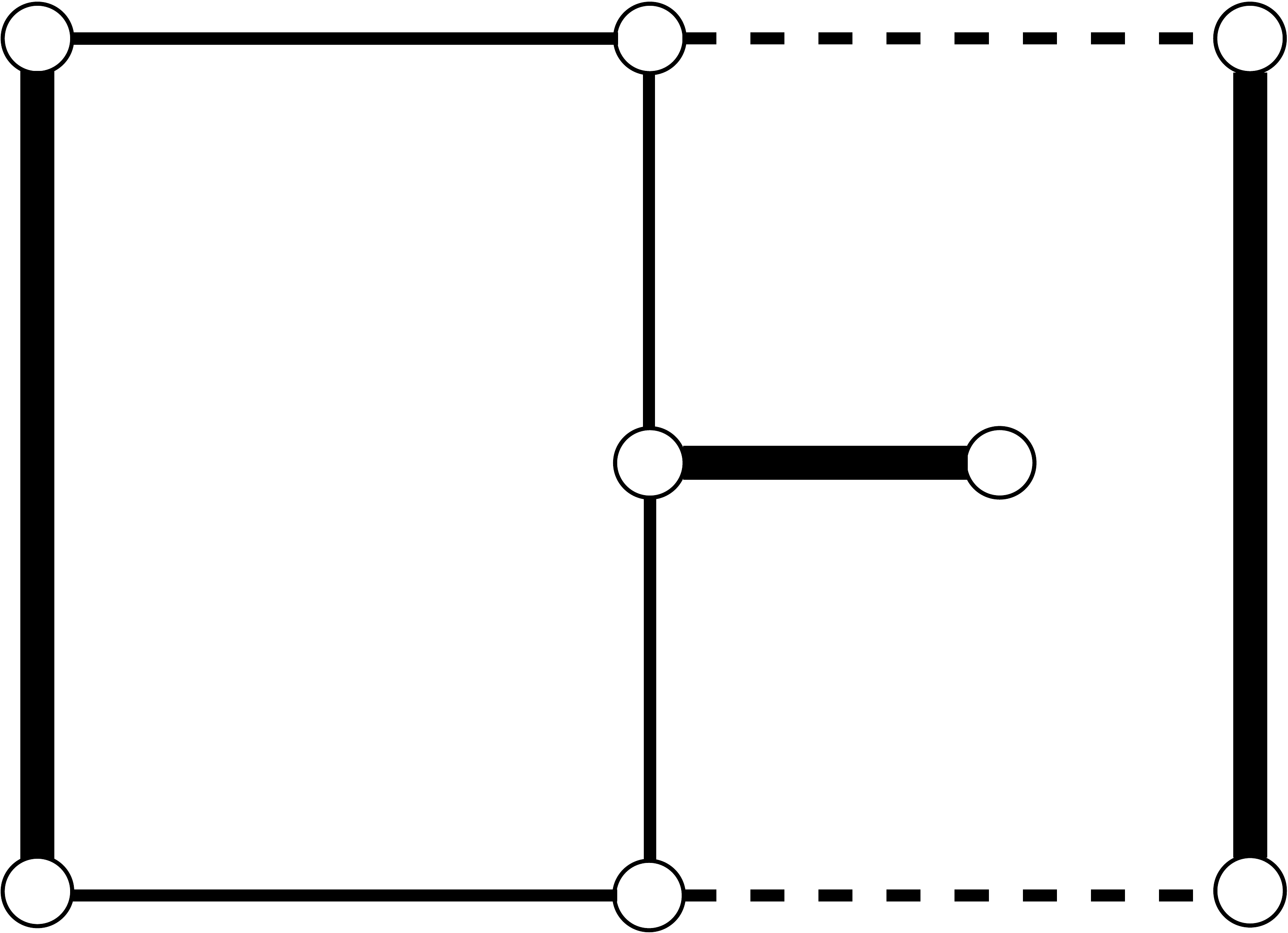}
\caption{A Coxeter diagram for the deformation $\widetilde{\rho}_t ( \widetilde{\Gamma}_{22} ) $ }
\label{coxeter_diagram2}
\end{figure}
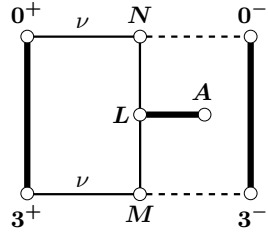

As shown in Section \ref{symmetries}, $\U{\Gamma}_{22}$ is generated by the elements of the set 
\[ \left\{ \p{0}, \m{0}, \p{3}, \m{3}, \l{A},\l{L},\l{M},\l{N} \right\}.\]
By abusing the standard notation slightly we can draw the Coxeter diagram of $\rho_t (\U{\Gamma}_{22})$ for $0<t<1$ as in figure \ref{coxeter_diagram2}.  Let us recall the meaning of the edges of figure \ref{coxeter_diagram1}.  The thick solid edges indicate walls which are tangent at infinity.  The dotted edges indicate walls which are separated by a positive distance.  The thin unlabeled solid edges indicate walls which intersect at angle $\pi/3$.  Vertices not joined by an edge indicate walls intersecting orthogonally.  Finally, edges labeled by $\nu$ indicate walls intersecting at angle 
\[ \frac{\pi}{\nu} = \arccos \sqrt{  \frac{ 2t^2}{1 + t^2} } . \]   One can check that this angle is $\theta/2$, where $\theta$ is the angle from Proposition \ref{angle prop}.

Note that $\nu$ is not necessarily an integer.  However, the values of $t$ for which $\nu$ is an integer will be of particular interest.

Recall the following special case of the Poincar{\'e} lemma.

\begin{thm}\cite{dlH}  \label{poincare lem v2}
Suppose $\Gamma< \text{Isom}(\mathbb{H}^n)$ is a group generated by reflections in a finite set of hyperplanes.  Assume these hyperplanes bound a convex region $P$ and every pair is either disjoint, tangent at infinity, or intersect at an angle of the form $\pi / n$ for some $n \in \{2,3,4,\ldots\}$.  Then $\Gamma$ is a discrete subgroup of $\text{Isom}(\mathbb{H}^n)$ with a fundamental domain given by $P$.
\end{thm}

Let's apply this theorem to the groups $\rho_t (\U{\Gamma}_{22})$.  First we must show that the hyperplanes of $\rho_t (\U{\Gamma}_{22})$ bound a convex domain in $\mathbb{H}^4$.  For this it suffices to find a time-like vector $v \in \mathbb{R}^{1,4}$ such that $(v, w(t)) <0$ for all $t>0$, where $w$ is any symbol from the list $\{ \p{0}, \m{0}, \p{3}, \m{3}, \l{A},\l{L},\l{M},\l{N} \}$.  For such a $v$, $v/\|v\|$ will be a point of $\mathbb{H}^4$ lying on the negative side of all $8$ generating walls of $\rho_t (\U{\Gamma}_{22})$.  A good first guess at $v$ is $(1,0,0,0,0)$, except this lies in the walls of $\l{L}$, $\l{M}$, and $\l{N}$.  Perturbing this to $v = (1, \frac{1}{10},\frac{1}{100},0,0)$ yields the desired vector for all $t$.  

With $v$ in hand we can apply Theorem \ref{poincare lem v2} to conclude that $\rho_t ( \U{\Gamma}_{22})$ is a discrete subgroup of $\text{Isom}(\mathbb{H}^4)$ if
\[ \nu(t) = \frac{ \pi }{ \arccos \sqrt{  \frac{ 2t^2}{1 + t^2} } } \]
is an integer.

Let $t_n <1$ be parameter values chosen such that $\nu(t_n) = n$, i.e.
\begin{equation} \label{discrete values}
t_n = \sqrt{ \frac{ \cos^2 (\frac{\pi}{n} ) }{ 2 - \cos^2 (\frac{\pi}{n})} }, 
\end{equation}
where $n \in \{ 2,3,4, \ldots\} $.  Note that $t_2 =0 < t_3 < t_4 < \ldots$, and $t_n \rightarrow 1$.  For each $n>2$ we see that $\rho_{t_n} (\U{\Gamma}_{22})$ is a discrete subgroup.  Therefore $\rho_{t_n} (\Gamma_{22})$ is also discrete.  We will write $\rho_{t_n}$ as simply $\rho_n$.  

Let us summarize these results in the following proposition.

\begin{prop} \label{global rigidity prop}
Consider the subgroup $\Gamma_{22}$ of $\text{Isom}(\mathbb{H}^4)$.  There exists a $1$-parameter family of representations $\rho_t : \Gamma_{22} \rightarrow \text{Isom}(\mathbb{H}^4)$ defined for $t>0$ such that $\rho_1$ is the inclusion map, and $\rho_t$ fixes the letter generators $\{\l{A}, \l{B} ,\ldots, \l{F} \}$ for all $t$.  Moreover, $\rho_t$ is the unique such family up to reparametrization of the $t$-coordinate.  There exists a sequence of parameter values $\{ t_n \}$ increasing to $1$ such that $\rho_n = \rho_{t_n}$ has discrete image for all $n$.  This sequence $\{ t_n \}$ is given explicitly by formula \ref{discrete values}.
\end{prop}

\begin{figure}[ht]
\labellist
\small\hair 2pt
\pinlabel $\l{A}$ [b] at 7 85
\pinlabel $\l{B}$ [t] at 7 287
\pinlabel $\l{C}$ [l] at 10 188
\pinlabel $\l{D}$ [r] at 496 188
\pinlabel $\l{E}$ [br] at 298 104 
\pinlabel $\l{F}$ [tl] at 210 271
\pinlabel $\pandm{0}$ [r] at 114 188
\pinlabel $\pandm{1}$ [b] at 186 14
\pinlabel $\pandm{2}$ [b] at 328 14
\pinlabel $\pandm{3}$ [l] at 394 188
\pinlabel $\pandm{4}$ [t] at 324 359
\pinlabel $\pandm{5}$ [t] at 182 359
\pinlabel $\pandm{6}$ [r] at 185 188
\pinlabel $\pandm{7}$ [l] at 324 188
\endlabellist
\centering
\includegraphics[scale=0.5]{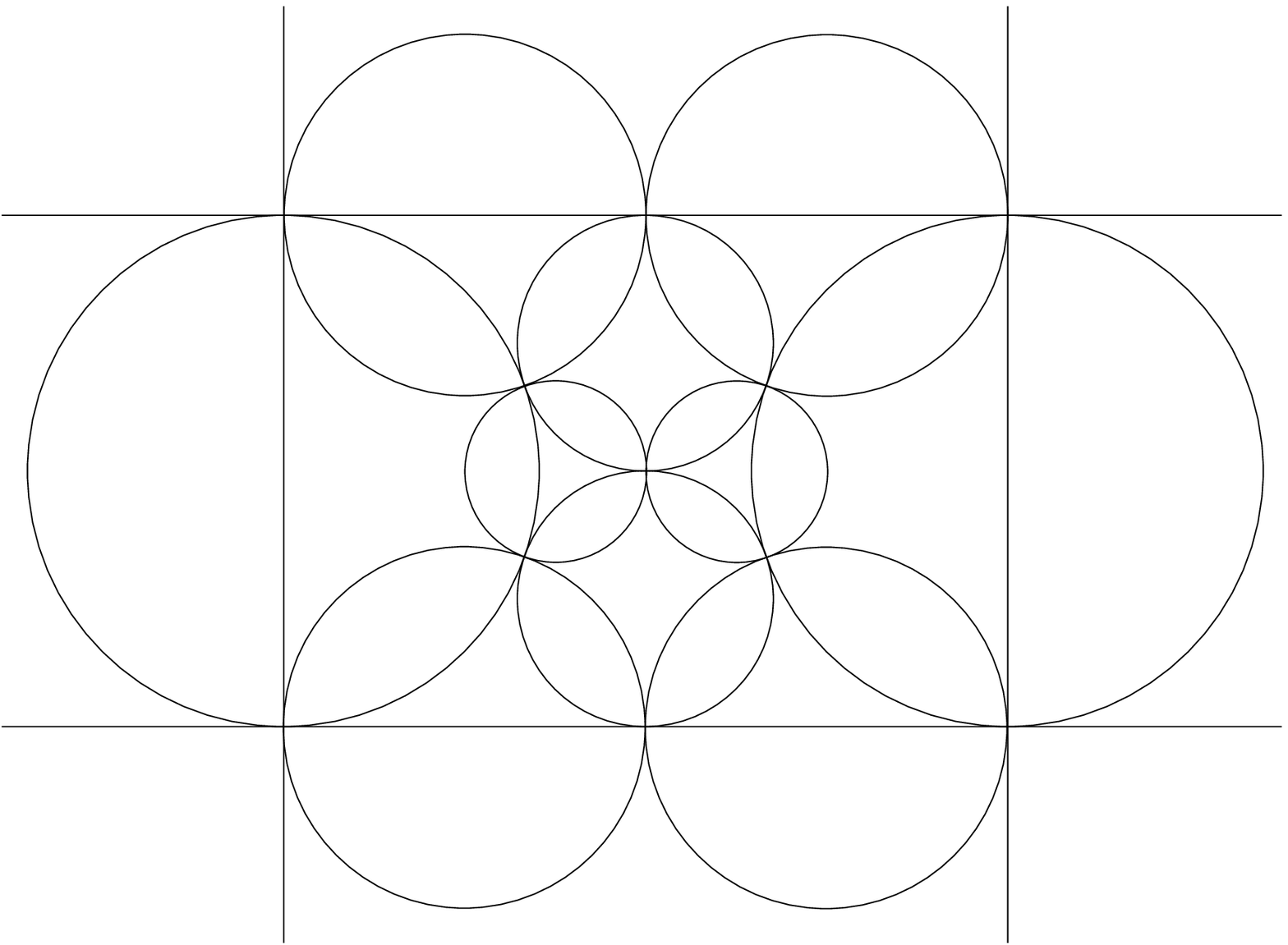}
\caption{The intersection pattern with $\partial V$}
\label{constant core intersection pattern}
\end{figure}

Notice that the letter space-like vectors $\l{A}$ through $\l{F}$ are all orthogonal to the vector $e_4 = (0,0,0,0,1)$.  Let $V \subset \Hf$ be the codimension $1$ hyperplane orthogonal to $e_4$,  $\partial V \subset \partial \Hf$ its boundary at infinity.  A remarkable feature of the deformation described in table \ref{qt table} is the pattern in $\partial V$ formed by its intersections with the walls corresponding to $\p{0}$, $\m{0}$, $\ldots$, $\m{7}$, $\l{A}$, $\ldots$, $\l{F}$.  Figure \ref{constant core intersection pattern} shows this pattern.  The notation $\pandm{0}$, $\pandm{1}$, etc., indicates that the triple intersections $\partial V \cap \partial W_{\p{0}} \cap \partial W_{\m{0}}$, $\partial V \cap \partial W_{\p{1}} \cap \partial W_{\m{1}}$, etc. are each $1$-dimensional, rather than the expected dimension of $0$.  

In figure \ref{constant core intersection pattern} there are $12$ points given by the common intersection of $4$ distinct circles, for example the intersection of $\l{B}$, $\l{C}$, $\pandm{0}$, and $\pandm{5}$.  These points are cusps of the initial group $\Gamma_{22}$ and they remain the cusps of the deformed group for all $t>0$.  To see this recall that the letter walls remain constant throughout the deformation.  In particular, they will always be orthogonal to $\partial V$ at infinity and intersect it in precisely the same way.  Their $12$ points of tangency at infinity correspond to the fixed points of the $12$ vertex cusp groups that must remain cusped by Lemma \ref{cube deformations}.  Since the reflections in the numbered walls 
are also part of these cusp groups in the same combinatorial pattern,  they must continue to pass through these same points at infinity.  Since each numbered wall passes through $3$ cusps, they must intersect $\partial V$ in the same circles (since $3$ points determine a circle).

What {\it does} change in this configuration is the angle at which the numbered walls intersect 
$\partial V$ (hence $V$).    Unlike the letter walls, which are orthogonal to $V$ throughout the deformation,  the numbered walls all begin at $t=1$ intersecting $V$ at angle $\frac {\pi} {4}$.  The angles of intersection of $V$ with the {\it negative walls} are all the same;  they go towards 
$\frac {\pi}{2}$ as $t \to 0$.  Similarly,  all the angles of intersection with the {\it positive walls} are equal and go towards $0$ as  $t \to 0$.   In fact,  all the positive walls converge to $V$ itself.

Thus there is a $3$-dimensional, totally geodesic subset that exists in all of the polyhedra throughout the deformation.  As $t \to 0$, the polyhedra collapse onto this subset.  Algebraically, the representations of $\Gamma_{22}$ converge to one preserving the $3$-dimensional geodesic subspace $V$.  The image of this limiting representation, viewed as being in $\Isom(\mathbb{H}^3)$, equals the group of reflections in the ideal {\it cubeoctahedron}.  This limiting process will be discussed further in Section \ref{cubeoctahedron section}.

\section{Viewed in the sphere at infinity} \label{viewed in the sphere at infinity}

The explicit list of vectors in table \ref{qt table} does not give much insight into the geometry of the deformation it determines.  For this we first return to the conformal model of spheres in $\mathbb{S}^3$, viewed in $\mathbb{R}^3$ using stereographic projection.  
We will then describe a family of polytopes  $\mathcal{F}_t$ in $\Hf$ which are deformations of the fundamental domain for $\Gamma_{22}$ and that, for certain values of $t$, serve as fundamental 
domains for the the discrete subgroup $\rho_t( \Gamma_{22})$ of $G$ described in Proposition \ref{global rigidity prop}.  The descriptions given in this section are valid for values of $t$ strictly between
$1$ and $\sqrt{3/5} \approx 0.775$, which correspond to the values of $\theta$ (the angle of intersection
between newly intersecting walls) strictly between $0$ and $ \pi/ 3 $.  We will worry later about other values of the deformation parameter.  

In figure \ref{24_cell figure 3} the arrangement of the $2$-spheres at infinity for the walls of $P_{24}$, the hyperbolic $24$-cell, is described  in terms of a cube in $\mathbb{R}^3$.  The planes determining the $6$ faces of the cube are $S_{\l{A}}, S_{\l{B}}, S_{\p{0}}, S_{\m{0}}, S_{\p{3}}$, and 
$S_{\m{3}}$. They meet in the point at infinity of $\mathbb{R}^3 \cup \infty$ which corresponds to one of the cusps of $P_{24}$.  
The edges of the cube are diameters of $12$ more of the spheres and, finally, there is the internal collection of  $6$ spheres pictured in figure \ref{24_cell figure 1} that meet in the center of the cube and are each tangent to one face of the cube.  The arrangement of spheres corresponding to the polygon 
$P_{22}$, defined by removing two walls from $P_{24}$, is obtained by removing the two spheres 
$S_{\l{G}}$ and $S_{\l{H}}$.   Their  diameters are opposite vertical edges of the cube (the front right and back left ones in figure \ref{24_cell figure 3}).
Note that $S_{\l{G}}$ and $S_{\l{H}}$ are disjoint, not even tangent.   

After we remove these two spheres there are still $80$ orthogonal pairs among the remaining $22$ spheres.  These correspond to the $80$ relations among the generating reflections of the group
$\Gamma_{22}$.   According to Proposition \ref{global rigidity prop}, there is a $1$-parameter family of nontrivial deformations of $\Gamma_{22}$, hence, of nontrivial deformations of these remaining spheres preserving the $80$ orthogonality conditions.
We will first try to understand this family by describing how the spheres at
infinity move.

By Lemma \ref{cube deformations} the walls $W_{\l{A}}, W_{\l{B}}, W_{\p{0}}, W_{\m{0}}, W_{\p{3}}$, and 
$W_{\m{3}}$ continue to represent a cusp group so their spheres at infinity continue to meet in
a single common point which we can keep as the point at infinity.  Then these spheres will continue to be planes.  They still determine a rectangular parallelpiped but we will see below that it will no longer be a cube. 

Because of the symmetries preserved by the deformation, the walls corresponding to the members of each of the original three octets -- the letter walls (of which there remain only $6$), the positive, and the negative walls --
all behave in the same manner.  But the three types of walls behave very differently.  We can see the $3$ types of behavior by looking at the faces of the cube and their intersections with the other
spheres.

\begin{figure}[ht]
\labellist
\small\hair 2pt
\pinlabel $\m{0}$ [b] at 451 359
\pinlabel $\p{0}$ [r] at 142 6
\pinlabel $\p{3}$ [l] at 324 6
\pinlabel $\m{3}$ [b] at 451 129
\pinlabel $\p{1}$ [bl] at 350 422
\pinlabel $\m{1}$ [tr] at 49 251
\pinlabel $\m{2}$ [bl] at 401 278
\pinlabel $\p{2}$ [tr] at 95 95
\pinlabel $(0,0)$ [tl] at 141 126
\pinlabel $(2t,0)$ [tl] at 324 126
\pinlabel $(0,2)$ [br] at 136 357
\pinlabel $\theta$ [tl] at 325 440
\endlabellist
\centering
\includegraphics[scale=0.4]{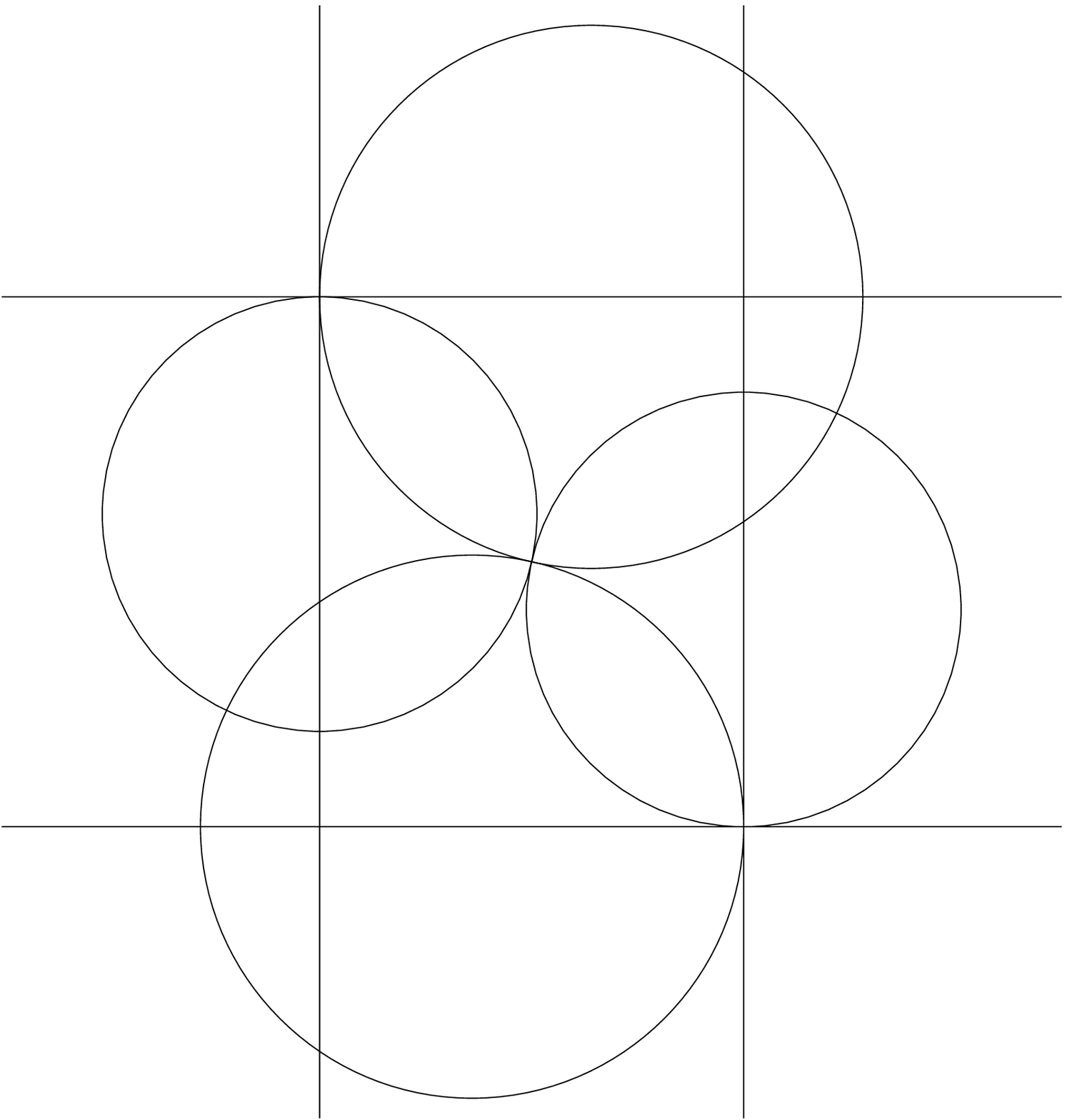}
\caption{Slicing along wall $\l{A}$ at $t = 0.8$.}
\label{circles2 figure}
\end{figure}

First, consider $S_{\l{A}}$ which is the $xy$-plane in figure \ref{24_cell figure 3}.  It is intersected
orthogonally by the $\pandm{0}$, $\pandm{1}$, $\pandm{2}$, $\pandm{3}$ spheres, forming the ideal octahedral pattern
seen in figure \ref{circles1 figure}.  The ideal vertices of this octahedron correspond to the
four vertices of the square, the point at infinity, and the center of the square.   These are each part of a cusp vertex of the original hyperbolic $24$-cell.  However, when the two walls are removed, two of the vertex subgroups of this octahedral reflection group (corresponding to the lower left and upper right vertices of the square in figure \ref{circles1 figure})  are no longer part of a rank $3$ cusp group and thus are no longer constrained to be represented by cusp groups.  The configuration of circles in 
$S_{\l{A}}$ determined by orthogonal intersection with the other eight spheres is no longer
rigid.  Figure \ref{circles2 figure} provides a picture of how they move.

Four of the original six vertices of the ideal octahedron are still part of rank $3$ cusps so they will remain cusped.  Keeping one cusped vertex at the point at infinity, four of the circles remain straight lines.  However,  the rectangle they determine is no longer a square.  At the two non-cusped vertices of the rectangle one circle pulls away from the line to which it was previously tangent, and one circle intersects a previously tangent line in an angle labeled $\theta$.  In other words, as required by Proposition \ref{octahedral deformations} for any nontrivial deformation, reflective hyperbolic Dehn filling has occurred at these two vertices.

The same process occurs on the parallel face of the cube, determined by $S_{\l{B}}$,  where
circles of intersection with the $\pandm{1}$ and $\pandm{2}$ spheres are replaced by intersections
with $\pandm{4}$ and $\pandm{5}$ spheres.  Specifically, on $S_{\l{A}}$ the circles of intersection with
$S_{\m{1}}$ and $S_{\m{2}}$ pull away from one of the sides of the rectangle to which they were
tangent while the circles of intersection with $S_{\p{1}}$ and $S_{\p{2}}$ push across one of
the sides.  On $S_{\l{B}}$ the circles of intersection with $S_{\m{4}}$ and $S_{\m{5}}$ pull away
while $S_{\p{4}}$ and $S_{\p{5}}$ push across.  This describes the behavior of all the spheres
whose diameters had been the horizontal edges of the cube.  The two remaining spheres
whose diameters had been vertical edges, $S_{\l{C}}$ and $S_{\l{D}}$, continue to be tangent
to the horizontal planes $S_{\l{A}}$ and $S_{\l{B}}$.

Finally, we note that the $6$ spheres which had been inside the cube continue to intersect
at a single point corresponding to a rank $3$ cusp of the deformed group.  Of these 
$S_{\l{E}}$ and $S_{\l{F}}$ continue to be tangent to the planes $S_{\l{A}}$ and $S_{\l{B}}$, respectively,
while $S_{\m{6}}$ and $S_{\m{7}}$ pull away from the planes $S_{\m{0}}$ and $S_{\m{3}}$
and spheres $S_{\p{6}}$ and $S_{\p{7}}$ poke through the planes $S_{\p{0}}$ and $S_{\p{3}}$.

With this qualitative description of the movement of the spheres we can describe how the intersection
pattern of the spheres with the vertical faces of the rectangular parallelpiped  change under
the deformation.
On the vertical faces  the initial configuration  of circles in figure \ref{24_cell figure 3} is 
conformally equivalent to that of the horizontal faces because of the symmetries of $P_{24}$.   However,  the effect of removing the walls $W_{\l{G}}$ and $W_{\l{H}}$
is more dramatic on the vertical faces because now one of the spheres that had provided a circle of intersection has been completely removed.  In $P_{22}$ the walls that had intersected these removed walls in $P_{24}$
now intersect the sphere at infinity in a conformal ideal triangle corresponding to the
face of intersection with $W_{\l{G}}$ or $W_{\l{H}}$ that has been removed.  Every numbered
wall has exactly one such removed face.

Consider the vertical face corresponding to $W_{\m{0}}$, a negative wall.  On the sphere at infinity the effect of removing the wall $W_{\l{G}}$ is that an ideal triangle that has been uncovered, corresponding to the intersection of $W_{\m{0}}$ with the sphere at infinity.  This ideal triangle  is bounded by the circles (line) of intersection of 
$S_{\m{0}}$ with $S_{\p{1}}$, $S_{\p{5}}$, and $S_{\p{3}}$.
During the deformation, these sides, which were tangent, begin to intersect in a positive angle
that will increase as $t$ decreases from $1$ towards $0$.  This can be seen in figure \ref{circles2 figure} where one sees $S_{\p{1}}$ and $S_{\p{3}}$ begin to intersect.  From the parallel figure 
corresponding to the top of the cube ($S_{\l{B}}$),  one sees the same thing occurring between 
$S_{\p{5}}$ and $S_{\p{3}}$.   Because their radii are increasing, the spheres $S_{\p{3}}$ and $S_{\p{5}}$ 
also begin to intersect.

\begin{figure}[ht]
\labellist
\small\hair 2pt
\pinlabel $\l{A}$ [bl] at 1 64
\pinlabel $\l{B}$ [br] at 252 191
\pinlabel $\p{0}$ [r] at 67 253
\pinlabel $\p{3}$ [l] at 193 253
\pinlabel $\p{1}$ [b] at 128 1
\pinlabel $\p{5}$ [t] at 128 254
\pinlabel $\l{C}$ [l] at 2 129
\pinlabel $(0,2,0)$ [tr] at 63 64
\pinlabel $(0,2,2)$ [br] at 63 194
\pinlabel $(2,2,0)$ [tl] at 194 63
\pinlabel {$W_{\m{0}} \cap \partial H_{\m{0}}$} [l] at 190 132
\endlabellist
\centering
\includegraphics[scale=0.66]{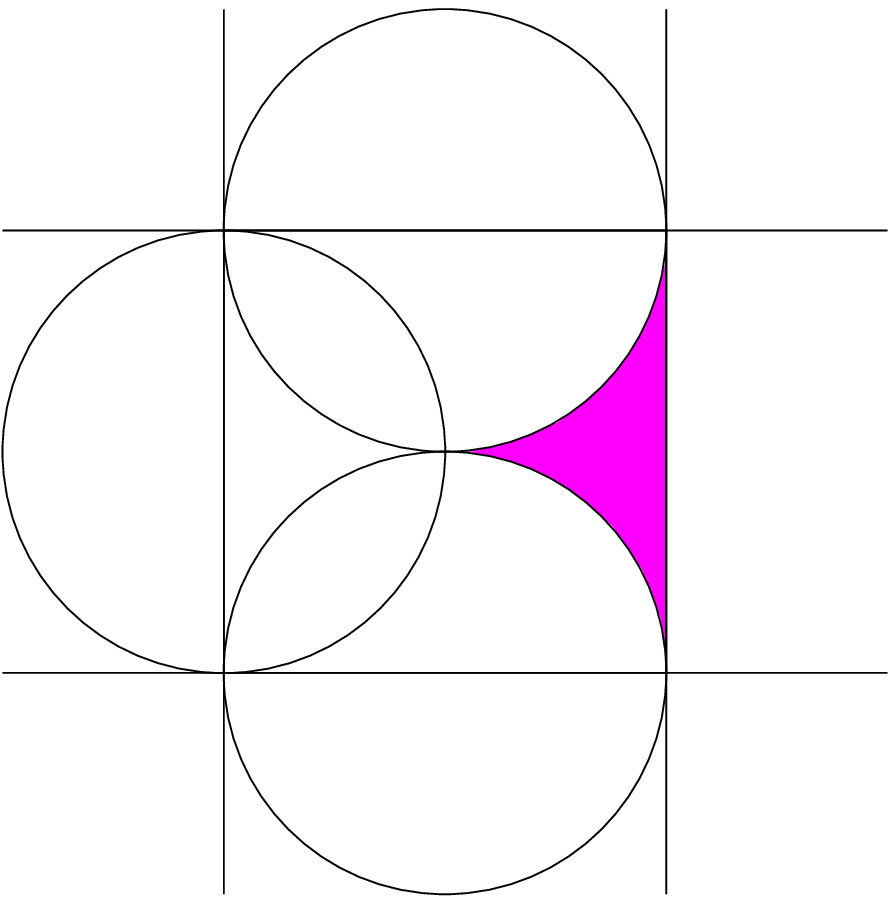}
\caption{Slicing along wall $\m{0}$ at $t=1$} 
\label{beginslicewith1.fig}
\end{figure}

\begin{figure}[ht]
\labellist
\small\hair 2pt
\pinlabel $\l{A}$ [bl] at 6 344
\pinlabel $\l{B}$ [bl] at 945 858
\pinlabel $\p{0}$ [tr] at 309 1187
\pinlabel $\p{3}$ [tl] at 769 1187
\pinlabel $\p{1}$ [b] at 601 54
\pinlabel $\p{5}$ [t] at 601 1143
\pinlabel $\l{C}$ [l] at 58 601
\pinlabel $(0,2,0)$ [tr] at 309 338
\pinlabel $(0,2,\sqrt{2}s)$ [br] at 309 858
\pinlabel $(2t,2,0)$ [br] at 763 344
\pinlabel $(\frac{1+t^2}{t},2,0)$ [bl] at 895 344
\pinlabel {$W_{\m{0}} \cap \partial H_{\m{0}}$} [l] at 814 596
\pinlabel $\theta$ at 786 1040
\endlabellist
\centering
\includegraphics[scale=0.16]{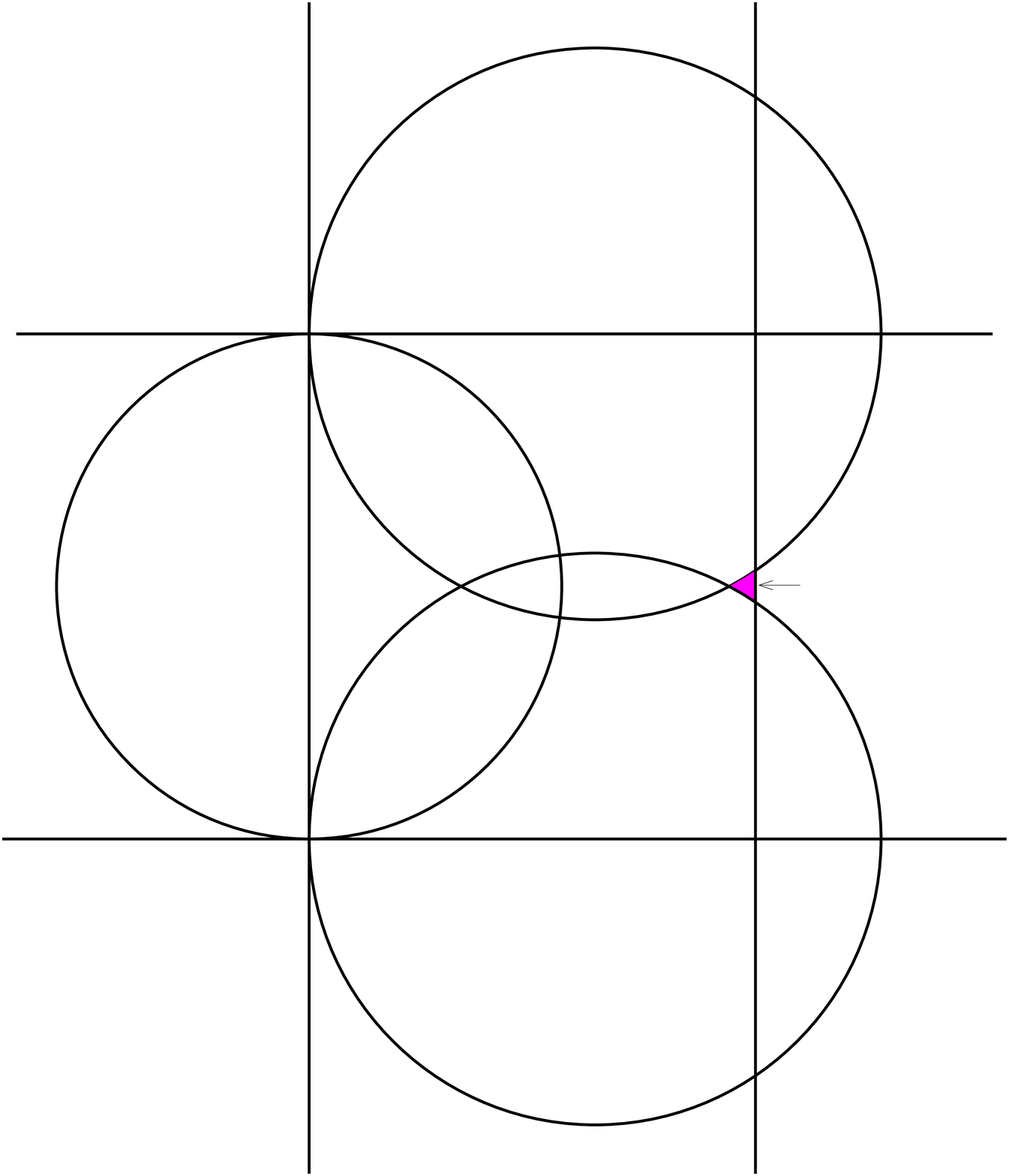}
\caption{Slicing along wall $\m{0}$ at $t=0.8$}
\label{slicewith1.fig}
\end{figure}

At $t=1$, $W_{\m{0}}$ has infinite volume and intersects  $\partial H_{\m{0}}$ in a conformal triangle with interior angles zero.  For the values of $t$ we are considering,  $W_{\m{0}}$ continues to have infinite volume and intersect  $\partial H_{\m{0}}$ in a conformal triangle.  However, in this case the triangle will have interior angles $\theta$.  This transition is described in figures \ref{beginslicewith1.fig} and \ref{slicewith1.fig}, where we have labeled the magenta triangle where $W_{\m{0}}$ intersects $\partial H_{\m{0}}$.

\begin{figure}[ht]
\labellist
\small\hair 2pt
\pinlabel $\l{A}$ [bl] at 1 64
\pinlabel $\l{B}$ [br] at 252 191
\pinlabel $\m{3}$ [r] at 67 253
\pinlabel $\m{0}$ [l] at 193 253
\pinlabel $\m{2}$ [b] at 128 1
\pinlabel $\m{4}$ [t] at 128 254
\pinlabel $\l{D}$ [l] at 2 129
\pinlabel $(2,0,0)$ [tr] at 63 64
\pinlabel $(2,0,2)$ [br] at 63 194
\pinlabel $(2,2,0)$ [tl] at 194 63
\pinlabel {$W_{\p{3}} \cap \partial H_{\p{3}}$} [l] at 190 132
\endlabellist
\centering
\includegraphics[scale=0.66]{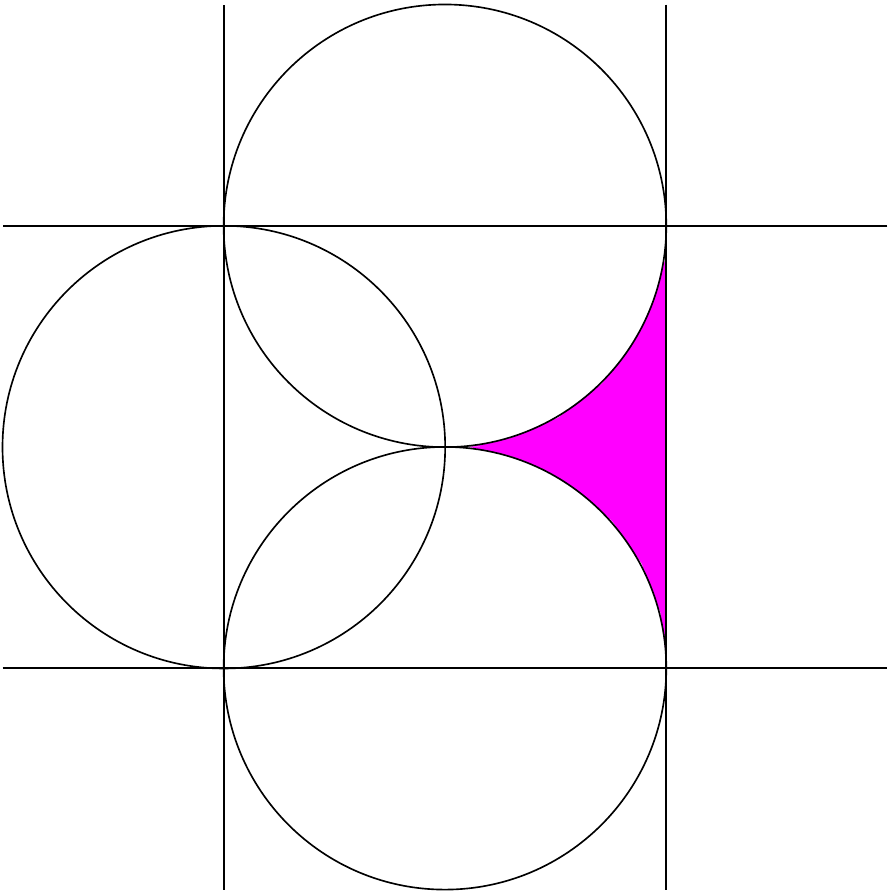}
\caption{Slicing along wall $\p{3}$ at $t=1$} 
\label{beginslicewith3.fig}
\end{figure}

\begin{figure}[ht]
\labellist
\small\hair 2pt
\pinlabel $\l{A}$ [b] at 27 312
\pinlabel $\l{B}$ [b] at 1220 822
\pinlabel $\m{3}$ [tr] at 339 1130
\pinlabel $\m{0}$ [tl] at 913 1135
\pinlabel $\p{1}$ [br] at 1089 149
\pinlabel $\p{5}$ [tr] at 1089 987
\pinlabel $\m{2}$ [bl] at 427 130
\pinlabel $\m{4}$ [tl] at 427 1007
\pinlabel $\l{D}$ [l] at 86 570
\pinlabel $\p{7}$ [r] at 746 576
\pinlabel $(2t,0,0)$ [tr] at 335 310
\pinlabel $(2t,0,\sqrt{2}s)$ [br] at 337 828
\pinlabel $(2t,2,0)$ [bl] at 909 314
\pinlabel {$W_{\p{3}} \cap \partial H_{\p{3}}$} [l] at 1051 562
\endlabellist
\centering
\includegraphics[scale=0.16]{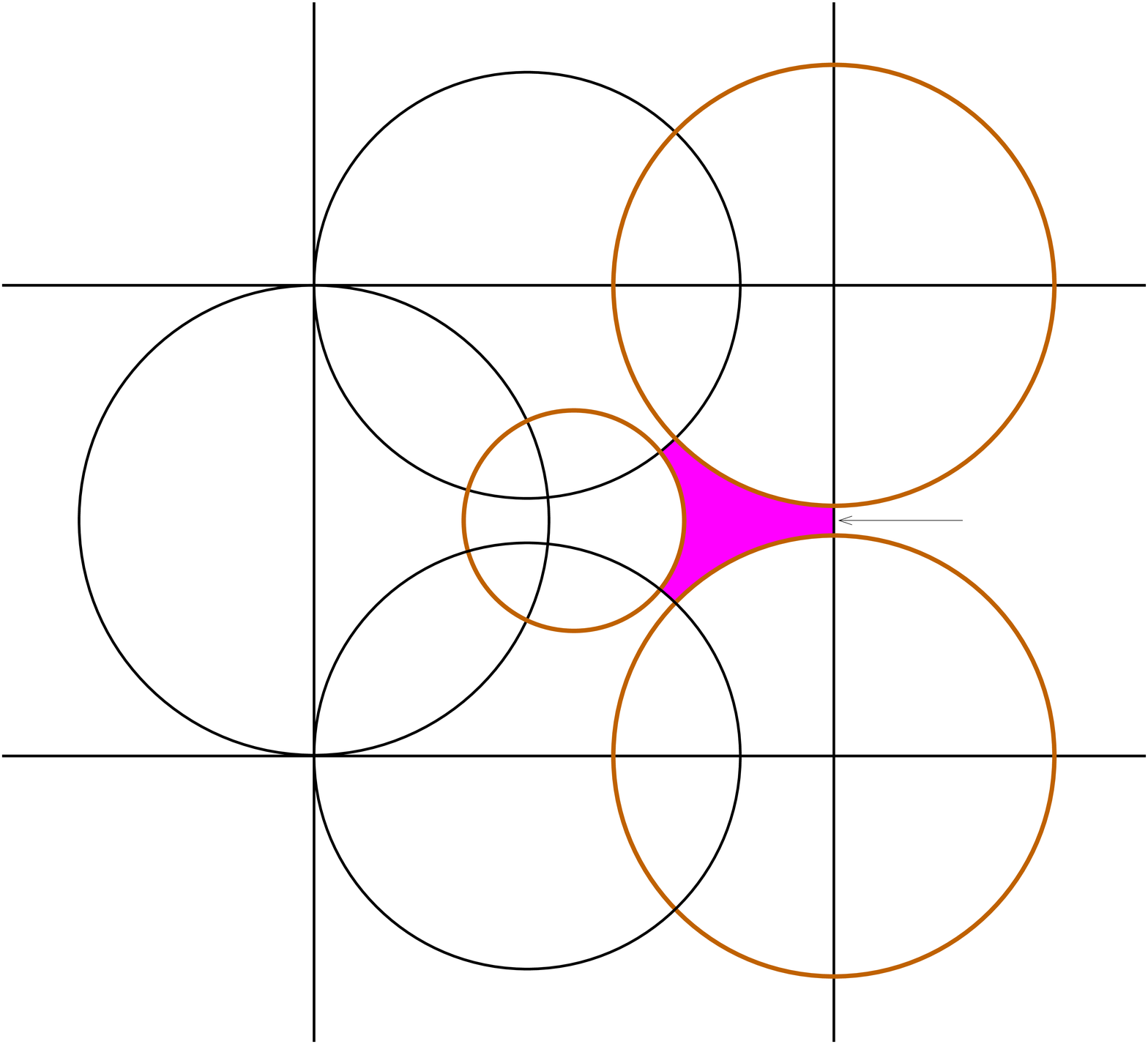}
\caption{Slicing along wall $\p{3}$} 
\label{slicewith3.fig}
\end{figure}

Similarly, we consider the vertical face corresponding to $W_{\p{3}}$, a positive wall.  In this case, the uncovered ideal triangle is bounded by the circles (line) coming from
intersection with $S_{\m{2}}$, $S_{\m{4}}$, and $S_{\m{0}}$.
During the deformation, these sides begin to pull apart and the uncovered region becomes
a right angled hexagon.  There is a further complication as the spheres $S_{\p{1}}, S_{\p{5}}$, 
and $S_{\p{7}}$ that were tangent to $S_{\p{3}}$ push through it, causing new intersections.  These 
create the $3$ new edges in the transition from a triangle to a hexagon.  The new intersections with
$S_{\p{1}}$ and $S_{\p{5}}$ can be seen in figure \ref{circles2 figure} and its counterpart figure for the top of the cube while $S_{\p{7}}$ is an internal sphere that begins to poke through.
This transition is described in  figures \ref{beginslicewith3.fig} and \ref{slicewith3.fig},  where we have labeled the magenta hexagon where $W_{\p{3}}$ intersects $\partial H_{\p{3}}$.  

There is a new phenomenon in this case coming from the new spheres that intersect  $S_{\p{3}}$.  In the previous two cases, all of the transverse intersections between spheres were at angle $\pi/2$.  In this case, whereas the spheres $S_{\l{V}}$ for  $\l{V} \in \{\m{0},\m{2},\m{3},\m{4},\l{A},\l{B},\l{D} \}$ intersect $S_{\p{3}}$ orthogonally, 
for $\l{V} \in \{\p{1},\p{5},\p{7} \}$ the intersection is at an angle $\theta$.  Therefore the dihedral angle associated to the faces $\p{1}$, $\p{5}$, and $\p{7}$ of $W_{\p{3}}$ will be $\theta$.  In figure \ref{slicewith3.fig}  orthogonal intersections are indicated by black circles and lines, while  intersections at angle $\theta$ are indicated by the three brown circles.

In order to support the previous qualitative description of the behavior of the $2$-spheres in the sphere at infinity,  we present  table \ref{centers and radii 2}, which is a list of the centers and radii  for the $22$ spheres as a function of $t$.   We have introduced a parameter $s = \sqrt{1+t^2}$.  

\begin{table}
$
\begin{array}{ccc}
\text{sphere} & \text{       center} & \text{radius} \\
S_{\m{0}} & \text{ the } $xz$-\text{plane shifted by the vector }  (0,2,0)\\
S_{\p{0}} & \text{ the } $yz$-\text{plane}\\
S_{\p{3}} & \text{ the } $yz$-\text{plane shifted by the vector }  (2t,0,0)\\ 
S_{\m{3}} & \text{ the } $xz$-\text{plane}\\ 
S_{\p{1}} & \left( \frac{s^2}{2t},2,0 \right) & \frac{s^2}{2t}\\
S_{\m{1}} & \left( 0,\frac{3 - t^2}{2},0 \right) & \frac{s^2}{2}\\
S_{\m{2}} & \left( 2t , \frac{s^2}{2},0\right) & \frac{s^2}{2}\\
S_{\p{2}} & \left( \frac{3 t^2 - 1}{2 t},0,0\right) & \frac{s^2}{2t}\\
S_{\p{5}} & \left( \frac{s^2}{2t},2,\sqrt{2} s\right) & \frac{s^2}{2t}\\
S_{\m{5}} & \left( 0,\frac{ 3 - t^2}{2} , \sqrt{2} s \right) & \frac{s^2}{2}\\
S_{\m{4}} & \left( 2t,\frac{s^2}{2},\sqrt{2} s \right) & \frac{s^2}{2}\\
S_{\p{4}} & \left( \frac{3 t^2 - 1}{2t},0,\sqrt{2} s \right)  & \frac{s^2}{2t}\\
S_{\m{6}} & \left( t, \frac{4 + s^2}{4}, \frac{\sqrt{2} s}{2}\right) & \frac{s^2}{4}\\
S_{\p{6}} & \left( t - \frac{s^2}{4t} , 1, \frac{\sqrt{2} s}{2} \right) & \frac{s^2}{4t}\\
S_{\p{7}} & \left( t +\frac{s^2}{4t},1,\frac{\sqrt{2} s}{2} \right) & \frac{s^2}{4t}\\
S_{\m{7}} & \left( t, \frac{3-t^2}{4}, \frac{\sqrt{2} s}{2} \right) & \frac{s^2}{4}\\
S_{\l{A}} & \text{ the } $xy$-\text{plane}\\
S_{\l{B}} & \text{ the } $xy$-\text{plane shifted by the vector }  (0,0,\sqrt{2} s)\\
S_{\l{C}} & \left( 0,2,\frac{\sqrt{2} s}{2} \right) & \frac{ \sqrt{2} s}{2} \\
S_{\l{F}} & \left( t,1,\frac{3 \sqrt{2} s}{4} \right) & \frac{\sqrt{2} s}{4}\\
S_{\l{E}} & \left( t,1,\frac{\sqrt{2} s}{4} \right) & \frac{ \sqrt{2} s}{4}\\
S_{\l{D}} & \left( 2t,0,\frac{\sqrt{2}s}{2}\right) & \frac{ \sqrt{2} s}{2}\\
\end{array}$
\caption{Centers and radii describing the deformation, where $s = \sqrt{1+t^2}$}
\label{centers and radii 2}
\end{table}

This list can be obtained from the list of space-like vectors in table \ref{qt table}, and stereographic projection.  Our choice of stereographic projection is a variant of the one used at the end of section 
\ref{24-cell}.  

Using this list we can see,  in a more computational form,  some of the qualitative features of the deformation as described above.  First, from the equations for the planes, one can see how the shape of the rectangular parallelpiped is changing.  Also, the process of deforming figure \ref{circles1 figure} into figure \ref{circles2 figure}  can be seen explicitly.  We see that the spheres $S_{\m{0}}$, $S_{\p{0}}$, $S_{\m{3}}$, and $S_{\l{A}}$ do not move.  The plane $S_{\p{3}}$ moves to the left with speed $2$.  The spheres $S_{\p{1}}$, $S_{\m{1}}$, $S_{\m{2}}$, and $S_{\p{2}}$ ``rotate'' about a central axis.  
One could also, using formulae from Euclidean geometry, compute the angle of intersection between
various pairs of spheres.  However, this is more easily done using the space-like vectors in 
table \ref{qt table}.
 


\section{The geometry of the polytopes}\label{fundamental domain}

 We will now build a family of convex infinite volume polytopes $\mathcal{F}_t$ in $\Hf$ that correspond to the hyperbolic reflection groups  $\rho_t( \Gamma_{22})$ and are determined by the hyperplanes corresponding to the space-like vectors in table \ref{qt table}.   These polytopes will be deformations of the fundamental region of $\Gamma_{22} < G$.  The hyperplanes intersect the sphere at infinity in 
 $2$-spheres whose intersection patterns we have just discussed.  This information will allow us to describe each of the walls.  The polytopes can then be described simply by giving the dihedral angles between the walls.  
 
A space-like vector determines not only a hyperplane but a half-space bounded by that hyperplane, 
consisting of all points in $\Hf$ whose dot product with that vector is non-positive.   Let $I$ be the set of symbols $\{\p{0}, \m{0}, \ldots , \m{7}, \l{A}, \ldots , \l{F} \}$.  Consider the list of
space-like vectors $\{\l{V},t \}_{\l{V} \in I}$ in table \ref{qt table}.  Each $\l{V}$ determines a one parameter family of hyperplanes $H_{\l{V},t}$ and halfspaces  $\{ \text{Half}_{\l{V},t} \}$.  With this notation we can define $\mathcal{F}_t$ to be the intersection
$$\mathcal{F}_t := \bigcap_{\l{V} \in I} \text{Half}_{\l{V} ,t}.$$
$\mathcal{F}_t$ will have $22$ walls given by the intersections $W_{\l{V} ,t} := \mathcal{F}_t \cap H_{\l{V}}$.

When all of the dihedral angles between walls are of the form $\pi/k$, where $k \ge 2$ is an integer,
the Poincar{\'e} lemma (Theorem \ref{poincare lem v2}) implies that $\mathcal{F}_t$ is a fundamental domain for $\rho_t( \Gamma_{22})$,  the group generated by reflections in the walls of $\mathcal{F}_t$.
The dihedral angles are $\pi/k$ when $t = t_{2k}$.  However,  our description of  the polytopes determined by the half-spaces corresponding to the space-like vectors in table \ref{qt table} is valid for \emph{all} values of $t$ strictly between $1$ and $\sqrt{3/5}$.

Because of the symmetries that are preserved by our deformation there will only be $3$ geometrically different types of walls, corresponding respectively to the members of the letter, negative, and positive walls.  Each  type of wall will be described both as a $3$-dimensional combinatorial object and as a $3$-dimensional hyperbolic polyhedron.

\subsection*{The geometry of $W_{\l{A},t}$} 

We first consider, for each $t$, the wall $W_{\l{A},t}$ of $\mathcal{F}_t$.  To simplify the notation, we will drop the $t$ subscript for the walls, it being understood that the unwritten $t$ takes a value strictly between $\sqrt{3/5}$ and $1$.  All the other letter walls are isometric to $W_{\l{A}}$.  Its faces are given by intersections with the $8$ walls $W_{\p{0}}$ through $W_{\m{3}}$.  These walls all intersect  $W_{\l{A}}$ orthogonally so the dihedral angle between them and $W_{\l{A}}$ is $\pi/2$.  The boundary of the hyperplane $H_{\l{A}}$ meets the sphere at infinity of $\Hf$ in a $2$-sphere which we denote by $\partial H_{\l{A}}$.  Its intersection pattern with  the boundaries of the hyperplanes $H_{\p{0}}$ through $H_{\m{3}}$ is represented by the black circles and lines of figure \ref{circles2 figure}, where the $xy$-plane represents $\partial H_{\l{A}}$.  

\begin{figure}[ht]
\labellist
\small\hair 2pt
\pinlabel $\m{0}$ [br] at 86 344
\pinlabel \textcolor{blue}{$\p{0}$} at 87 159
\pinlabel $\p{3}$ at 152 267
\pinlabel $\m{3}$ at 147 79
\pinlabel \textcolor{blue}{$\p{1}$} at 163 361
\pinlabel \textcolor{blue}{$\m{1}$} at 200 243
\pinlabel $\m{2}$ at 272 324
\pinlabel \textcolor{blue}{$\p{2}$} [tl] at 244 132 
\pinlabel {Front faces are labeled in black.} [tl] at 208 59
\pinlabel {Back faces are labeled in \textcolor{blue}{blue}.} [tl] at 208 40
\endlabellist
\centering
\includegraphics[scale=0.6]{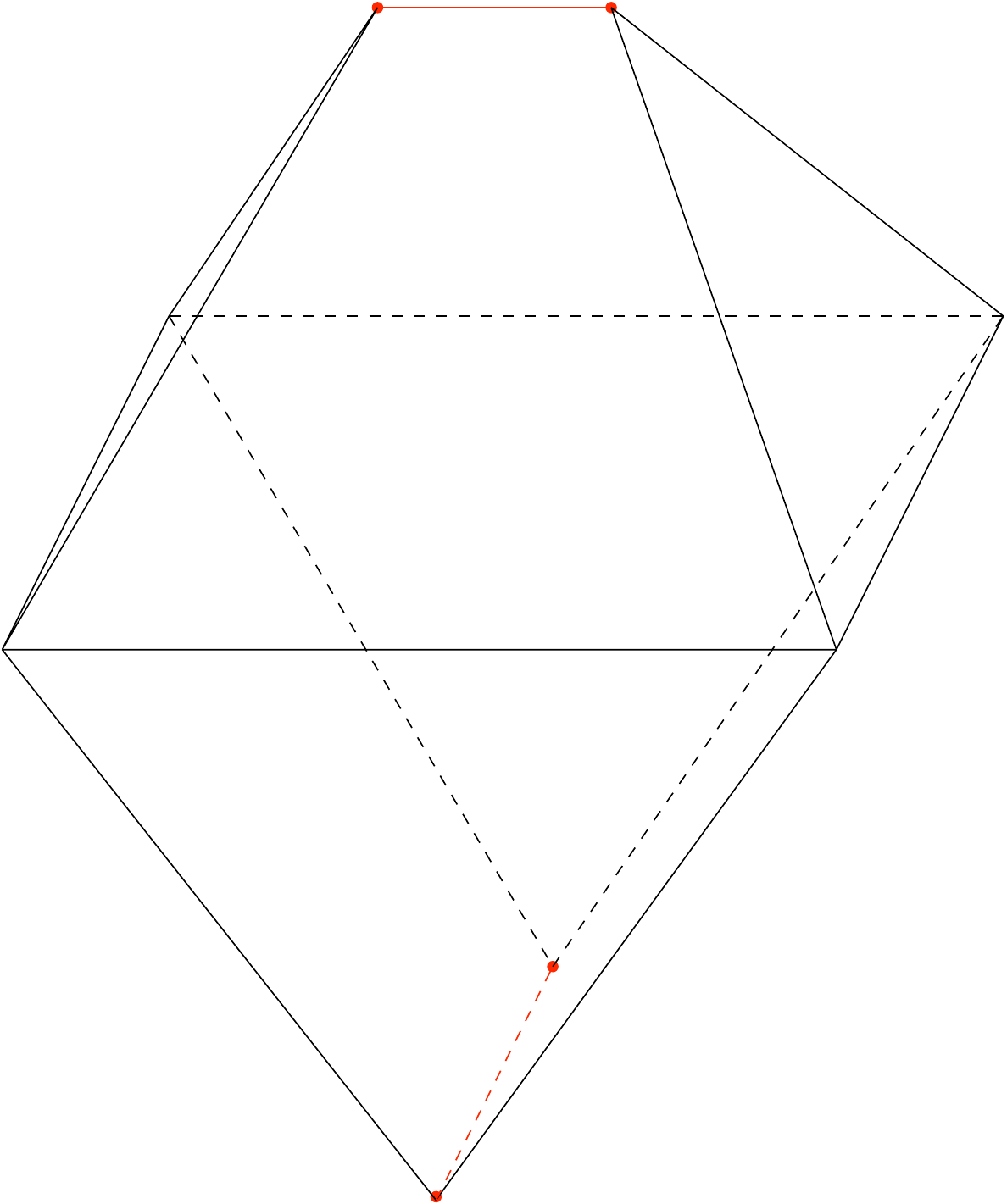}
\caption{The combinatorics of the wall $W_{\l{A}}$}
\label{deformedoctahedron figure}
\end{figure}

\begin{figure}[p!]
\includegraphics[scale=0.4,angle=-90]{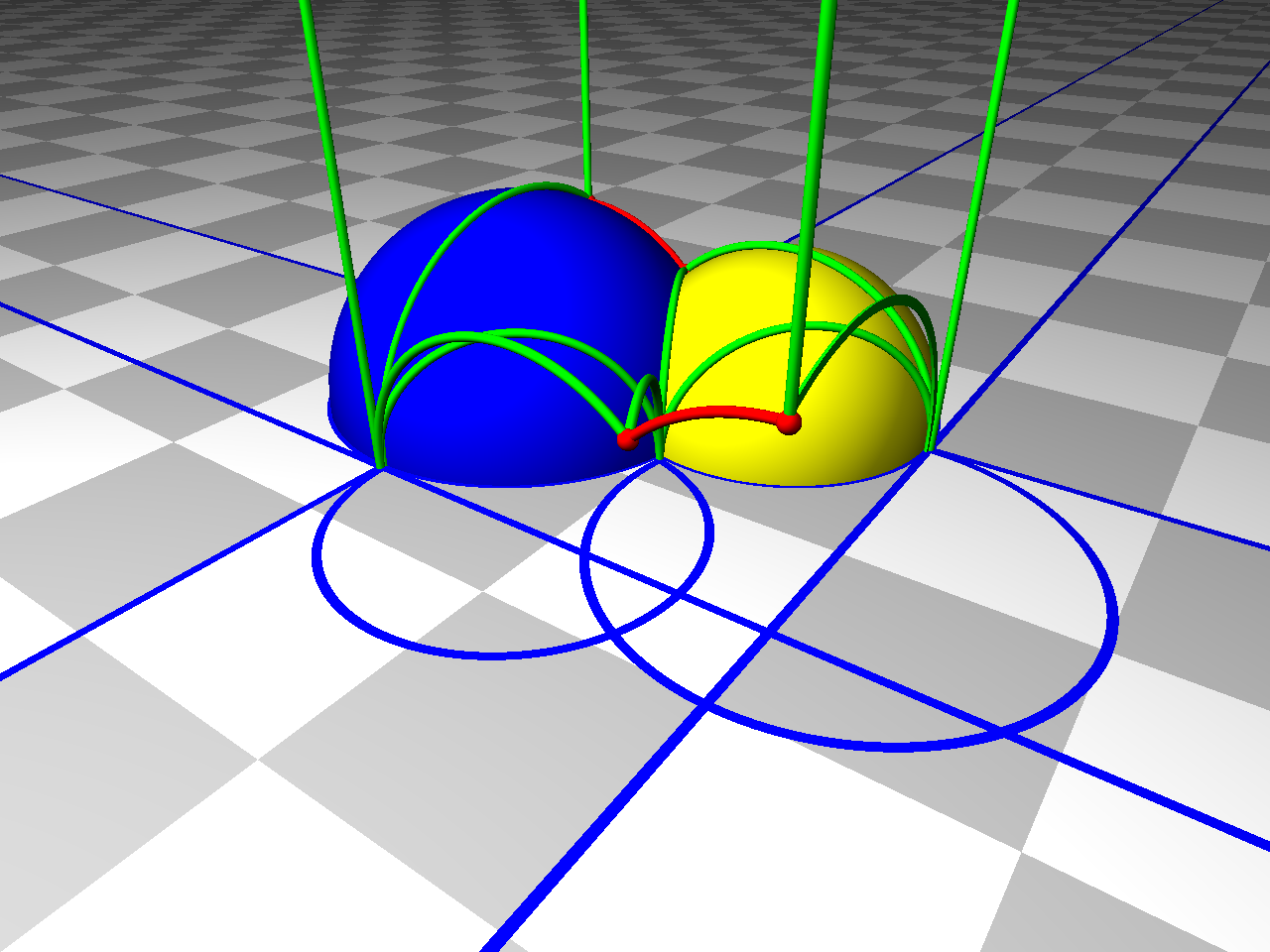}
\caption{The geometry of the wall $W_{\l{A}}$}
\label{octahedron3 figure}
\end{figure}

From this pattern on the sphere at infinity it is easy to construct the pattern of intersections of the hyperplanes  themselves and to see the geometry of  the wall $W_{\l{A}}$.  We view $\partial H_{\l{A}}$
as the boundary of the upper halfspace model of $\mathbb{H}^3$.  The intersections of the
hyperplanes $H_{\p{0}}$ through $H_{\m{3}}$ with $ H_{\l{A}}$ are hyperbolic planes that are 
represented by $2$-spheres perpendicular to $\partial H_{\l{A}}$.  Note that, because the hyperplanes intersect $W_{\l{A}}$ orthogonally, the dihedral angles between these hyperbolic planes, hence between the faces $W_{\l{A}}$ are the same as the dihedral angles between the hyperplanes that create those faces.   A geometric depiction of this configuration of geodesic planes is given by figure \ref{octahedron3 figure} 
while figure \ref{deformedoctahedron figure} provides a combinatorial model of the resulting polyhedron.

Let us explain the features of figure \ref{deformedoctahedron figure}.  It represents an eight sided polyhedron in $\mathbb{H}^3$.  The red edges are geodesics of length $\ell$, where $\ell$ is the function of $t$ given in Proposition \ref{angle prop}.  The red vertices are therefore points in $\mathbb{H}^3$.  The remaining vertices are points at infinity, and the black edges are infinite or semi-infinite geodesics.  The dihedral angle of the black edges is $\pi/2$.  The dihedral angle of the red edges is $\theta$, where $\theta$ is the function of $t$ given in Proposition \ref{angle prop}.  Each face of the polyhedron is labeled with a number.  This number indicates which wall intersects $W_{\l{A}}$ along that face.  For example, the face labeled $\p{3}$ represents the intersection $W_{\l{A}} \cap W_{\p{3}}$.  Faces on the front of the polyhedron are labeled with black numbers.  Faces on the back of the polyhedron are labeled with blue numbers.  (The back faces are the faces of figure \ref{deformedoctahedron figure} the reader could not see if the polyhedron were opaque.)  

Figure \ref{octahedron3 figure} provides a geometrically accurate picture of this polyhedron, viewed in an upper halfspace corresponding to $H_{\l{A}}$.  The blue circles and lines in the $xy$-plane show the intersections of $\partial H_{\p{0}}$ through $\partial H_{\m{3}}$ with $\partial H_{\l{A}}$.  The blue hemisphere represents the intersection $H_{\p{2}} \cap H_{\l{A}}$ and the yellow hemisphere represents $H_{\m{1}} \cap H_{\l{A}}$.  The edges of $W_{\l{A}}$ are green and red.  The green edges have dihedral angle $\pi/2$; the red edges have dihedral angle $\theta$.  In the figure the deformation parameter $t$ is $0.8$, which means $\theta \approx 55.8\,^{\circ}$.  The red balls indicate vertices of $W_{\l{A}}$ that lie in $H_{\l{A}}$, as opposed to being ideal vertices at infinity.    Notice that $W_{\l{A}}$ is not compact, but it does have finite volume.  

\subsection*{The geometry of $W_{\m{0}}$}

Next we will analyze the geometry of the wall $W_{\m{0}}$.  All of the other negative walls will be isometric to it. 

When $t=1$ this wall is an ideal octahedron with a Fuchsian end added to one face (see Section \ref{def:orbifold}), resulting from the removal of  the wall $W_{\l{G}}$  when forming $P_{22}$. The
remaining faces are formed by orthogonal intersections with seven other walls.  Thus, the dihedral angles between $W_{\m{0}}$ and each of these walls is $\pi/2$.  The pattern of intersection at infinity of the corresponding $2$-spheres is given in figure \ref{beginslicewith1.fig} for $t=1$.   The deformed intersection pattern is pictured in figure \ref{slicewith1.fig}.     Using this information, one can deduce the geometry of $W_{\m{0}}$.  As before, this is done by taking an upper half-space lying over the plane in figure \ref{slicewith1.fig} (which represents $\partial H_{\m{0}}$) and inserting seven geodesic planes whose intersections at infinity produce the pattern of that figure.   The polyhedron in 
$\mathbb{H}^3$ they bound is identified with $W_{\m{0}}$.  As before,  because the walls intersect $W_{\m{0}}$ orthogonally, the dihedral angles between the faces $W_{\m{0}}$ are the same as the dihedral angles between the hyperplanes that create those faces.   All of the latter angles are either 
$\pi/2$ or $\theta$ and can be read off from either  Proposition \ref{angle prop} or figure \ref{slicewith1.fig}.

\begin{figure}[ht]
\labellist
\small\hair 2pt
\pinlabel $\p{0}$ at 179 677
\pinlabel $\l{B}$ at 64 537 
\pinlabel $\l{C}$ at 294 537
\pinlabel $\p{5}$ at 179 518
\pinlabel $\textcolor{blue}{\l{A}}$ at 179 616
\pinlabel $\textcolor{blue}{\p{3}}$ at 93 380
\pinlabel $\textcolor{blue}{\p{1}}$ at 269 380
\pinlabel {Front faces are labeled in black.} [tl] at 10 294
\pinlabel {Back faces are labeled in \textcolor{blue}{blue}.} [tl] at 10 270
\endlabellist
\centering
\includegraphics[scale=0.4]{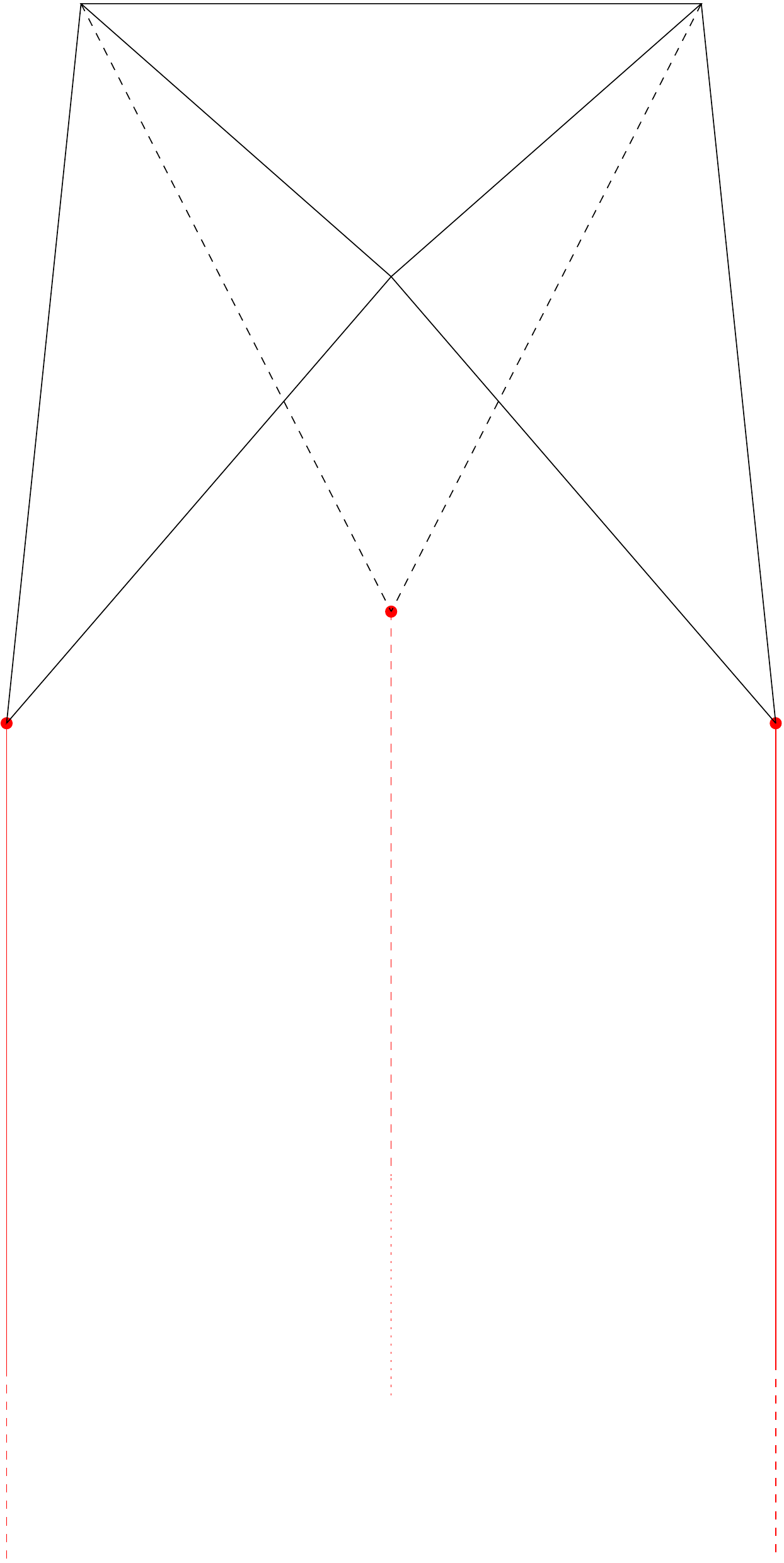}
\caption{The combinatorics of the wall $W_{\m{0}}$}
\label{combinatoricsof1 figure}
\end{figure}

\begin{figure}[p!]
\includegraphics[scale=0.4,angle=-90]{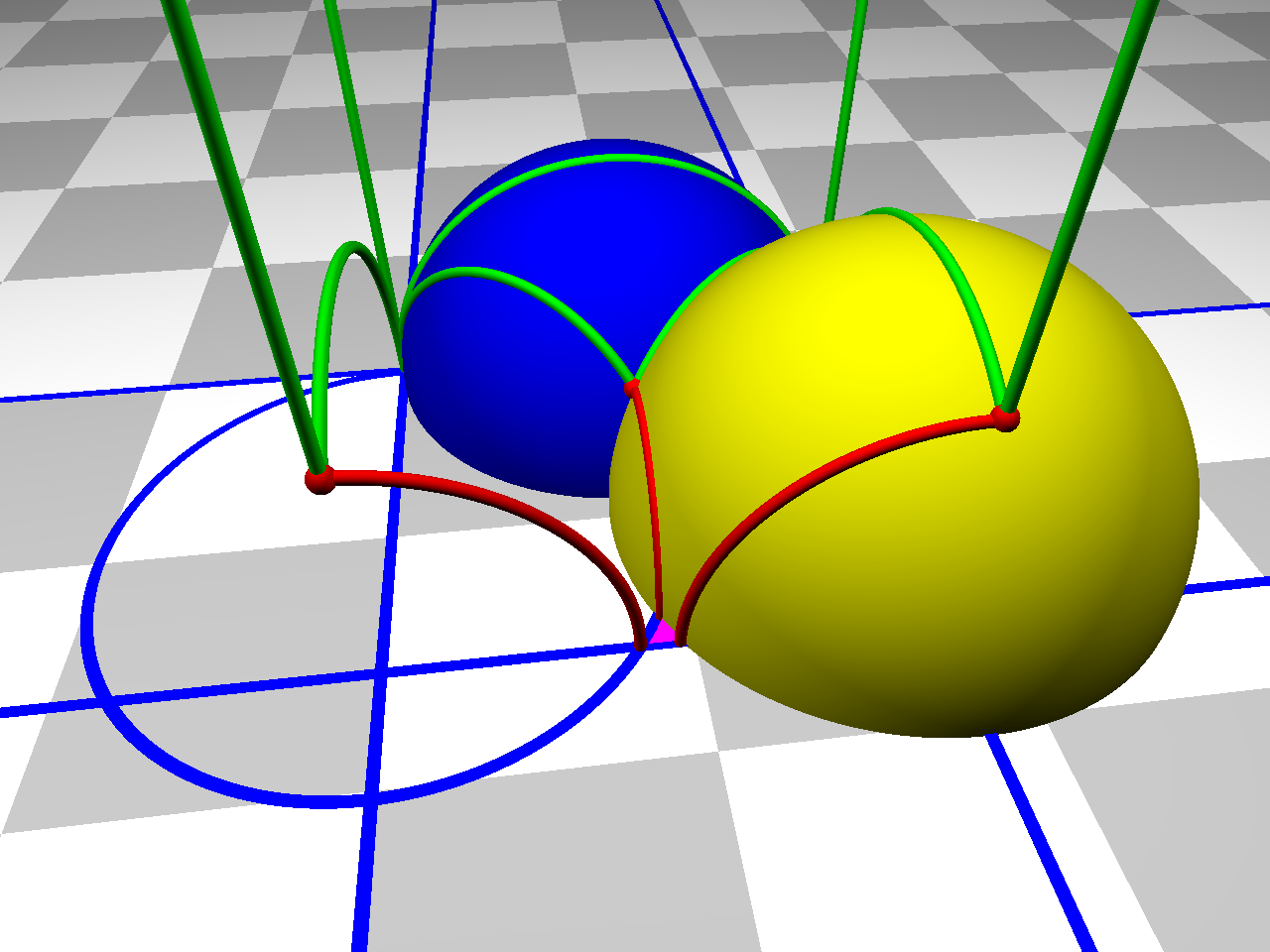}
\caption{The geometry of the wall $W_{\m{0}}$}
\label{octahedron4 figure}
\end{figure}

The combinatorics of $W_{\m{0}}$ are described by figure \ref{combinatoricsof1 figure}, while the geometry of $W_{\m{0}}$ is depicted in figure \ref{octahedron4 figure}.  

Figure \ref{combinatoricsof1 figure} shows a seven sided unbounded polyhedron.  Each side represents a face of $W_{\m{0}}$.  A face of $W_{\m{0}}$ corresponds to an intersection of $W_{\m{0}}$ with another wall of $\mathcal{F}_t$, and each face of figure \ref{combinatoricsof1 figure} is labeled with the wall it intersects.  Black edges correspond to edges of $W_{\m{0}}$ with dihedral angle $\pi/2$.  Red edges correspond to edges with dihedral angle $\theta$.  The fact that the polyhedron of figure \ref{combinatoricsof1 figure} is unbounded signifies that $W_{\m{0}}$ has infinite volume and nonempty intersection with the sphere at infinity $\partial H_{\m{0}}$, as was visible in figure \ref{slicewith1.fig}.

Figure \ref{octahedron4 figure} shows the edges of $W_{\m{0}}$ sitting inside the $3$-dimensional hyperbolic space $H_{\m{0}}$.  The checkerboard plane is the sphere at infinity of $H_{\m{0}}$.  The blue lines and circles in the plane represent the intersections of $H_{\p{0}}$, $H_{\p{1}}$, $H_{\p{3}}$, $H_{\p{5}}$, $H_{\l{A}}$, $H_{\l{B}}$, and $H_{\l{C}}$ with $\partial H_{\m{0}}$.  This planar pattern is a reproduction of figure \ref{slicewith1.fig}.  One should imagine $2$-dimensional hyperplanes in $H_{\m{0}}$ producing these intersection patterns at infinity.  Figure \ref{octahedron4 figure} shows two of the hyperplanes.  The blue hemisphere represents the intersection $H_{\l{C}} \cap H_{\m{0}}$ and the yellow hemisphere represents the intersection $H_{\p{5}} \cap H_{\m{0}}$.  The green edges represent edges of $W_{\m{0}}$ with dihedral angle $\pi/2$.  The red edges represent edges with dihedral angle $\theta$.  In the figure the deformation parameter $t$ is $0.8$, and thus $\theta \approx 55.8\,^{\circ}$.  The red balls indicate vertices of $W_{\m{0}}$ which lie in $H_{\m{0}}$ (as opposed to being ideal vertices at infinity).  Similar to figure \ref{slicewith1.fig}, the magenta triangle is again the intersection $W_{\m{0}} \cap \partial H_{\m{0}}$.

\subsection*{The geometry of $W_{\p{3}}$}

The next wall to examine is $W_{\p{3}}$.  The other positive walls are isometric to it.

As before, one can use the intersection pattern at infinity (figure \ref{slicewith3.fig}) and knowledge of the dihedral angles between hyperplanes (Proposition \ref{angle prop})  to deduce the geometry of $W_{\p{3}}$.   At $t=1$, $W_{\p{3}}$ has infinite volume and intersects the sphere at infinity $\partial H_{\p{3}}$ in a conformal triangle with interior angles zero.  Similarly, for the values of $t$ considered here,$W_{\p{3}}$ has infinite volume.  However, its intersection with infinity is more complicated than before.  The faces of $W_{\p{3}}$ corresponding to hyperplanes $\p{1}$, $\p{5}$, and $\p{7}$ do not exist when $t=1$.  These new intersections cause the intersection $W_{\p{3}} \cap \partial H_{\p{3}}$ to become a right angled hexagon.  This is visible in figure \ref{slicewith3.fig}, where we have labeled the magenta hexagon where $W_{\p{3}}$ intersects $\partial H_{\p{3}}$.

\begin{figure}[ht]
\labellist
\small\hair 2pt
\pinlabel $\m{3}$ at 259 519
\pinlabel $\textcolor{blue}{\l{D}}$ at 259 475 
\pinlabel $\textcolor{blue}{\p{7}}$ at 259 145
\pinlabel $\m{0}$ at 259 339
\pinlabel $\l{B}$ at 163 408
\pinlabel $\l{A}$ at 326 408
\pinlabel $\textcolor{blue}{\m{4}}$ at 89 429
\pinlabel $\textcolor{blue}{\m{2}}$ at 417 429
\pinlabel $\p{5}$ at 64 239
\pinlabel $\p{1}$ at 462 239
\pinlabel {Front faces are labeled in black.} [tl] at 151 100
\pinlabel {Back faces are labeled in \textcolor{blue}{blue}.} [tl] at 151 82
\endlabellist
\centering
\includegraphics[scale=0.5]{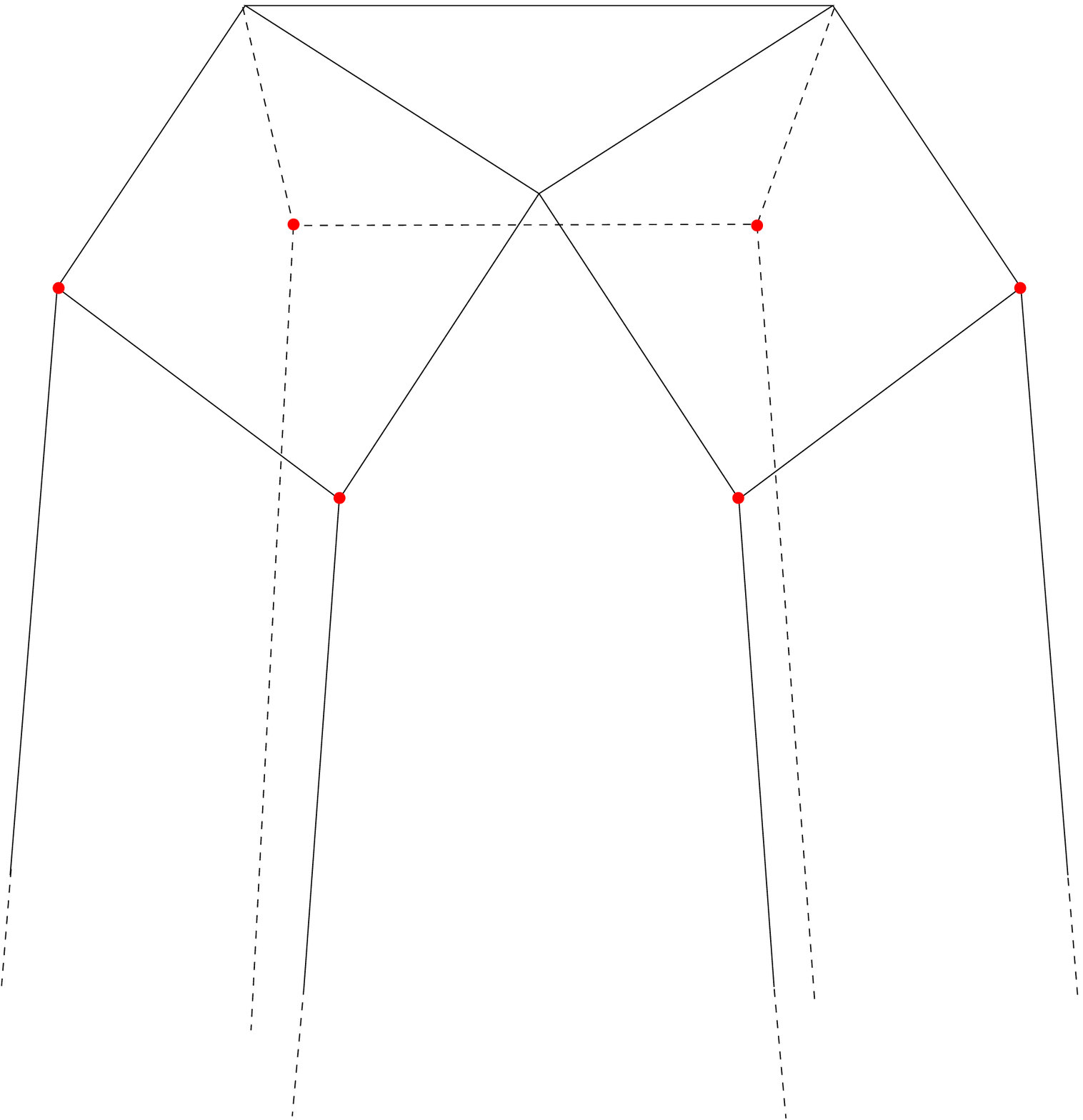}
\caption{The combinatorics of the wall $W_{\p{3}}$}
\label{combinatoricsof3 figure}
\end{figure}

\begin{figure}[p!]
\includegraphics[scale=0.4,angle=-90]{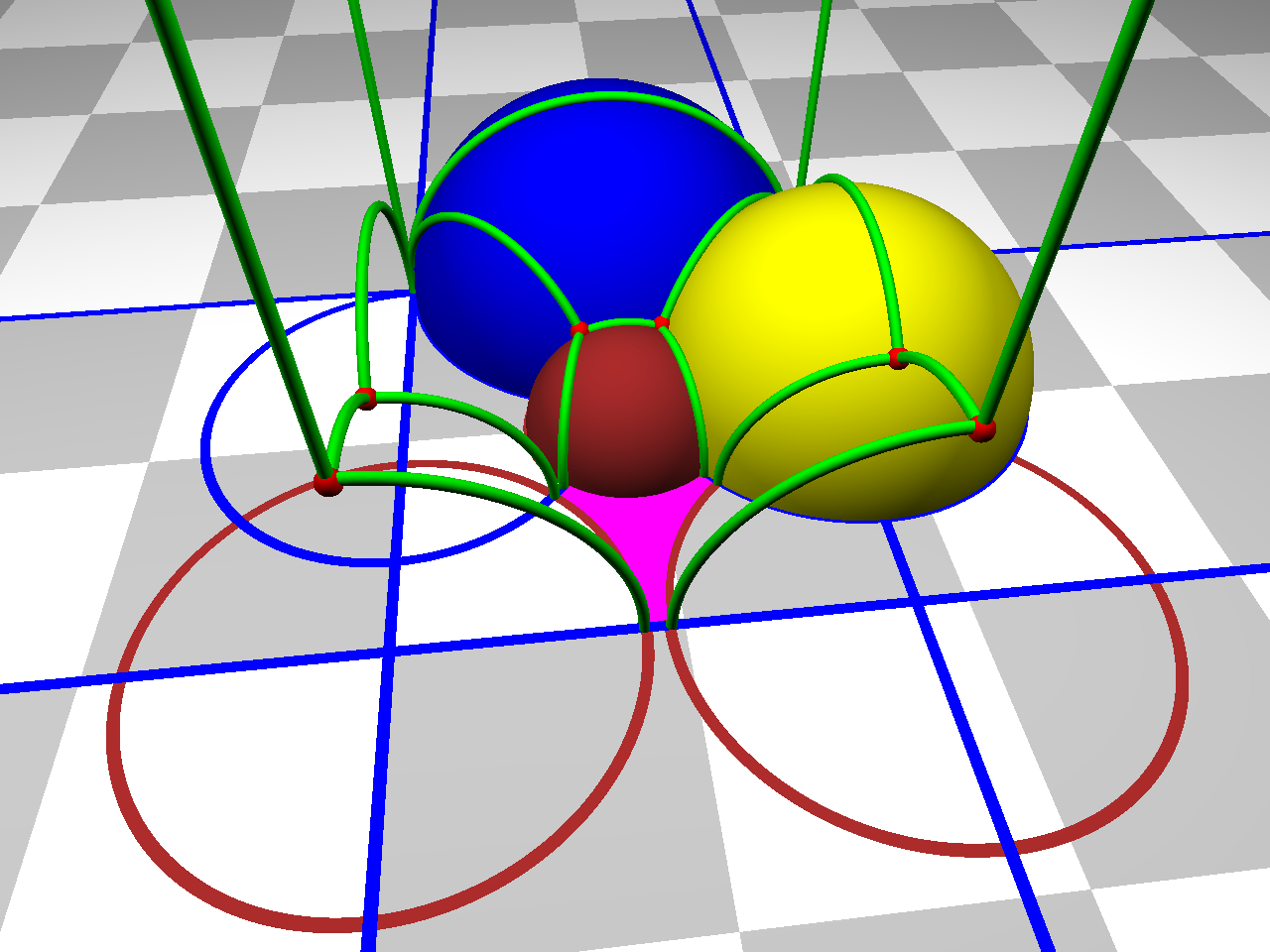}
\caption{The geometry of the wall $W_{\p{3}}$}
\label{octahedron5 figure}
\end{figure}

The combinatorics of $W_{\p{3}}$ are described by figure \ref{combinatoricsof3 figure}, while the geometry of $W_{\p{3}}$ is depicted in figure \ref{octahedron5 figure}.  Let us describe these figures.  

Figure \ref{combinatoricsof3 figure} shows a ten sided unbounded polyhedron.  Each side represents a face of $W_{\p{3}}$.  A face of $W_{\p{3}}$ corresponds to an intersection of $W_{\p{3}}$ with another wall of $\mathcal{F}_t$, and each face of figure \ref{combinatoricsof3 figure} is labeled with the wall it intersects.   The fact that the polyhedron of figure \ref{combinatoricsof3 figure} is unbounded signifies that $W_{\p{3}}$ has infinite volume and nonempty intersection with the sphere at infinity $\partial H_{\p{3}}$, as was visible in figure \ref{slicewith3.fig}.

Note that, in contrast with the previous cases, the dihedral angles between $W_{\p{3}}$ and the hyperplanes that create its faces are not all $\pi/2$.  In particular, the angles between $W_{\p{3}}$ and the hyperplanes $\p{1}$, $\p{5}$, and $\p{7}$ equal $\theta$,  the angle expressed as a function of
$t$ in Proposition \ref{angle prop}.  (The remaining intersections are orthogonal.)  Thus,  the dihedral angles for $W_{\p{3}}$ do not generally correspond to the dihedral angles between the corresponding
hyperplanes.  However, we can determine them as follows.  The edges of the face corresponding to the hyperplane $\p{1}$  are created by intersection with the hyperplanes $\l{A}$, $\m{0}$, and $\m{2}$ (and $W_{\p{3}}$), respectively.  These intersect $W_{\p{3}}$ and the $\p{1}$ hyperplane orthogonally.  It follows that the dihedral angle of these edges is $\pi/2$.  Similarly, the dihedral angles of the edges of the faces corresponding to the $\p{5}$ and $\p{7}$ hyperplanes are $\pi/2$.  The remaining edges are formed by the orthogonal intersection of $W_{\p{3}}$ with two hyperplanes that meet each other orthogonally so their dihedral angles are also $\pi/2$.  Hence all of the dihedral angles between the faces of  $W_{\p{3}}$ are $\pi/2$.

Analogous to the previous cases, figure \ref{octahedron5 figure} shows the edges of $W_{\p{3}}$ sitting inside the $3$-dimensional hyperbolic space $H_{\p{3}}$.  The checkerboard plane is again the sphere at infinity of $H_{\p{3}}$.  The blue lines and circles in the plane represent the orthogonal intersections of $H_{\m{0}}$, $H_{\m{3}}$, $H_{\m{2}}$, $H_{\m{4}}$, $H_{\l{A}}$, $H_{\l{B}}$, and $H_{\l{D}}$ with $\partial H_{\p{3}}$.  The brown circles represent the angle $\theta$ intersections of $H_{\p{1}}$, $H_{\p{5}}$, and $H_{\p{7}}$ with $\partial H_{\p{3}}$.  This planar pattern is a reproduction of figure \ref{slicewith1.fig}.  To help the imagination, figure \ref{octahedron5 figure} shows three of the bounding hyperplanes.  The blue hemisphere represents the intersection $H_{\l{D}} \cap H_{\p{3}}$,  the yellow hemisphere represents the intersection $H_{\m{4}} \cap H_{\p{3}}$, and the brown hemisphere represents the intersection $H_{\p{7}} \cap H_{\p{3}}$.  As discussed above, all the edges of $W_{\p{3}}$ have dihedral angle $\pi/2$.  The red balls indicate vertices of $W_{\p{3}}$ which lie in $H_{\p{3}}$.  Similar to figure \ref{slicewith3.fig}, the magenta hexagon is the intersection $W_{\p{3}} \cap \partial H_{\p{3}}$.

We can now use the symmetries of our deformation to deduce the geometry of the other walls.  In Section \ref{defining the deformation 2} it was shown that the deformation has a symmetry group which permutes the walls of $\mathcal{F}_t$.  This action has three orbits,  the letter walls, positive numbered walls, and the negative numbered walls.  We have examined one wall from each orbit.  All the walls in an orbit have the same geometry, so we have determined the geometry and combinatorics of all $22$ walls of $\mathcal{F}_t$.  Since we have also determined the angle of intersection between the walls, this provides a complete description of the deformed polytopes.



\section{Geometry of the ends} \label{miscellaneous section}

In this section we will study the geometry of the ends of the polytopes $\mathcal{F}_{t}$.  We will restrict ourselves to values of $t$ satisfying $1> t > \sqrt{3/5}$.  This corresponds to values of the angle $\theta$ less than $\pi/3$ and the analysis from Section \ref{fundamental domain} will
apply. Understanding the geometry of the ends will, in particular, imply certain properties of the groups $\Lambda_n = \rho_{t_n}(\Gamma_{22})$ for $n = 2 m$ where $m \ge 4$.    For these even values of $n$ the angle $\theta$ equals $2 \pi / n = \frac{\pi}{m}$, a submultiple of $\pi$.  We can therefore directly apply the Poincar{\'e} lemma to the convex region $\mathcal{F}_{t_n}$ to conclude that it is a fundamental domain for $\Lambda_n$.

The polytope $P_{22}$ constructed by removing the two walls $W_{\l{G}}$ and $W_{\l{H}}$
from the hyperbolic $24$-cell, $P_{24}$, has infinite volume with two ends whose intersections with the sphere at infinity of $\Hf$ have non-empty interior.  These ends are bounded by the walls that had previously intersected
the removed walls in $P_{24}$.  Specifically, the end corresponding to the wall  $W_{\l{G}}$ is
bounded by the odd, positive numbered walls and the even, negative numbered walls. We saw in the previous sections that for the values of $t$ we are considering, this end will continue to be bounded by the
same collection of walls.  However, its geometry will vary with $t$.  We will focus
on this end which we will denote by $E_{\l{G}}$.  The other end, corresponding to $W_{\l{H}}$, will be isometric to it via the roll symmetry $R$.

At the beginning, when $t=1$, the ends of $P_{22}$ are Fuchsian.  Recall that an end being Fuchsian 
means by definition that the reflections in the walls determining that end preserve a common 
$3$-dimensional hyperplane in $\Hf$.  A natural question is whether this property holds for
any of the deformed polytopes.  In fact, this is never the case when $1 > t > \sqrt{3/5}$.

To see this, first consider the even, negative numbered walls, $W_{\m{0}}$, $W_{\m{2}}$, $W_{\m{4}}$,
and $W_{\m{6}}$, which bound $E_{\l{G}}$.  The corresponding space-like vectors are
equal to $\sqrt{2} e_0 \pm e_1 \pm e_2 \pm e_3 - t e_4$, where there is an even number of negative
signs among the $\pm e_i, i=1,2,3$.  It is easy to check that, up to multiplication by a scalar,
the unique vector in $\mathbb{R}^{1,4}$ orthogonal to all $4$ of these vectors equals $(-t, 0, 0, 0, \sqrt{2})$.
This is a space-like vector which implies that there is a unique $3$-dimensional hyperplane orthogonal
to all $4$ of the corresponding walls.  Hence, there is a unique hyperplane preserved by reflection
in the even, negative numbered walls.  The odd, positive numbered walls,  $W_{\p{1}}$, $W_{\p{3}}$, $W_{\p{5}}$, and $W_{\p{7}}$, correspond to the space-like vectors  $\sqrt{2} e_0 \pm e_1 \pm e_2 \pm e_3 - \frac{1}{t} e_4$, where there is an odd number of negative
signs among the $\pm e_i, i=1,2,3$.  Up to scale, the only vector orthogonal to all $4$ of these
vectors is $(-\frac{1}{t}, 0, 0, 0, \sqrt{2})$.  This equals the vector orthogonal to the even, negative walls if and only if $t =1$.  Thus, there is no hyperplane preserved
by reflection in all $8$ walls bounding $E_{\l{G}}$ except at the beginning of the family of
polytopes.

While this analysis holds for all values of $t \neq 1$,
we will see in Sections \ref{disappear} and \ref{manual lattice} that the topology of the ends
of $\mathcal{F}_{t}$ changes at  $t = \sqrt{3/5}$.  At that point the infinite volume ends of the polytope 
will be bounded by fewer walls;  the negative walls will no longer be involved.  Our computation above
shows that there will be a hyperplane preserved by reflection in the odd, positive numbered walls
when $\sqrt{3/5} \geq t > \sqrt{1/2}$ and, hence, that the end $E_{\l{G}}$ will be Fuchsian.  
At $t = \sqrt{1/2}$
the vector $(-\frac{1}{t}, 0, 0, 0, \sqrt{2})$ becomes light-like, implying that the four walls will 
have a common point on the sphere at infinity.  We will also see in these later sections that
the polytope has finite volume at this stage and that these ends will no longer exist for all smaller
values $\sqrt{1/2} > t > 0$.  The four odd, positive
numbered walls will have a common point in the interior of $\Hf$.

\begin{figure}[ht]
\labellist
\small\hair 2pt
\pinlabel $\p{7}$ at 229 179
\pinlabel $\p{5}$ at 145 301
\pinlabel $\m{6}$ at 29 134
\pinlabel $\m{4}$ at 302 350
\pinlabel $\textcolor{blue}{\m{0}}$ at 87 319
\pinlabel $\textcolor{blue}{\m{2}}$ at 280 35
\pinlabel $\textcolor{blue}{\p{3}}$ at 249 284
\pinlabel $\textcolor{blue}{\p{1}}$ at 78 228
\pinlabel {Front faces are labeled in black.} [tr] at 134 55 
\pinlabel {Back faces are labeled in \textcolor{blue}{blue}.} [tr] at 134 37 
\endlabellist
\centering
\includegraphics[scale=0.6]{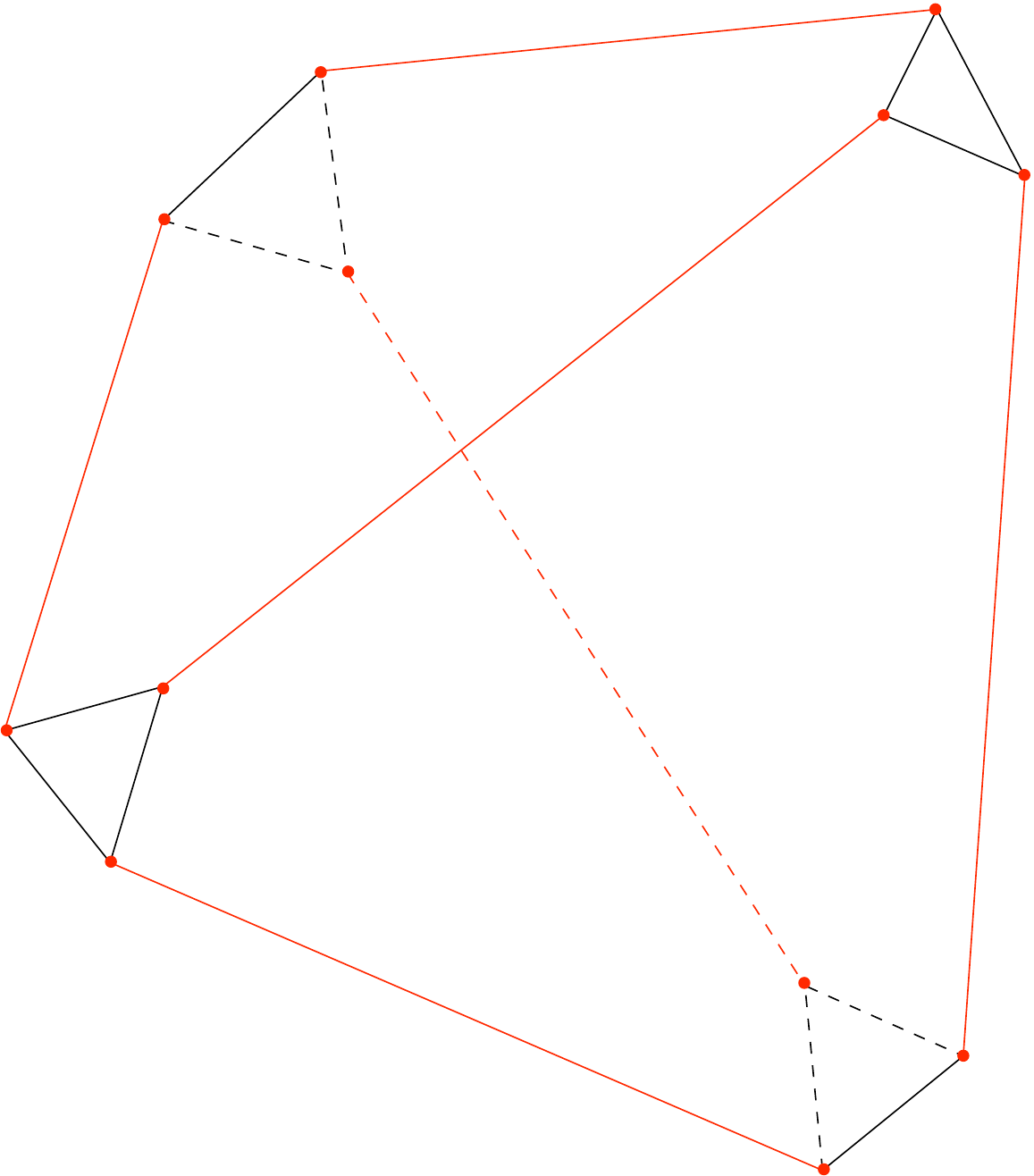}
\caption{The combinatorics of a component of the conformal boundary}
\label{dd1 figure}
\end{figure}
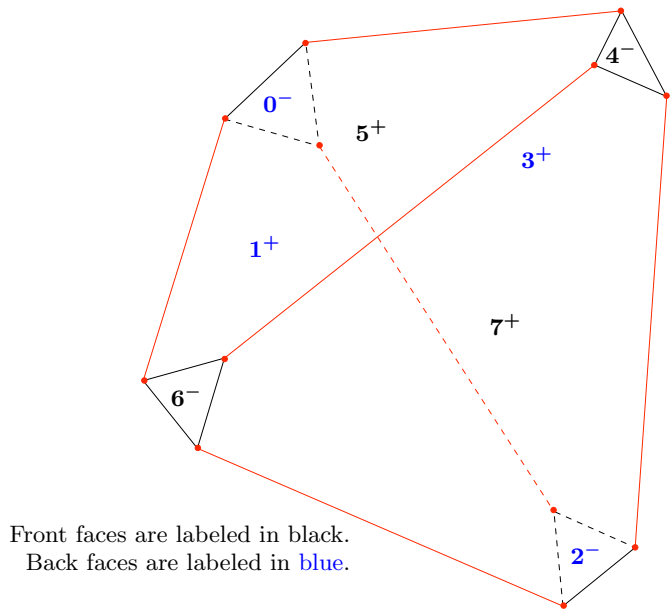

Returning to the case when there is a nontrivial end, we can further study its geometry at infinity.  Specifically, we can consider the closure, $\overline{\mathcal{F}}_{t}$, of $\mathcal{F}_{t}$ 
as a subset of the compactified topological space $\Hf \cup \mathbb{S}^3_\infty$.  
We will refer to the components of
$\overline{\mathcal{F}}_{t} \cap \mathbb{S}^3_\infty$ as the boundary components of
$\overline{\mathcal{F}}_{t}$. They inherit a conformally flat structure from the sphere at
infinity.
Recall figures \ref{octahedron4 figure} and \ref{octahedron5 figure}, depicting the geometry of $W_{\m{0}}$ and $W_{\p{3}}$ respectively.  These two walls each contribute one face to a component of  $\partial \overline{\mathcal{F}}_{t}$.  In the figures these faces are shown as magenta polygons in the complex plane.  Specifically, $W_{\m{0}}$ contributes a triangle and $W_{\p{3}}$ contributes a hexagon.  Recall that the walls indexed by $\p{5}$, $\p{7}$, and $\p{1}$ are all isometric to $W_{\p{3}}$, and the walls indexed by $\m{4}$, $\m{6}$, and $\m{2}$ are isometric to $W_{\m{0}}$.  Using these facts one can piece together the eight faces of the component $C$ of $\partial \overline{\mathcal{F}}_{t}$ corresponding to $E_{\l{G}}$, thereby determining its combinatorics.  These combinatorics  are depicted in figure \ref{dd1 figure}; the polyhedron is the familiar truncated tetrahedron.  As in previous similar figures, the front faces are labeled with black numbers and the back faces are labeled with blue numbers.  Faces are labeled with the number of the wall that intersects them.  The red edges have dihedral angle $\theta$.  The black edges have dihedral angle $\pi/2$.

\begin{figure}[p!]
\includegraphics[scale=0.4,angle=-90]{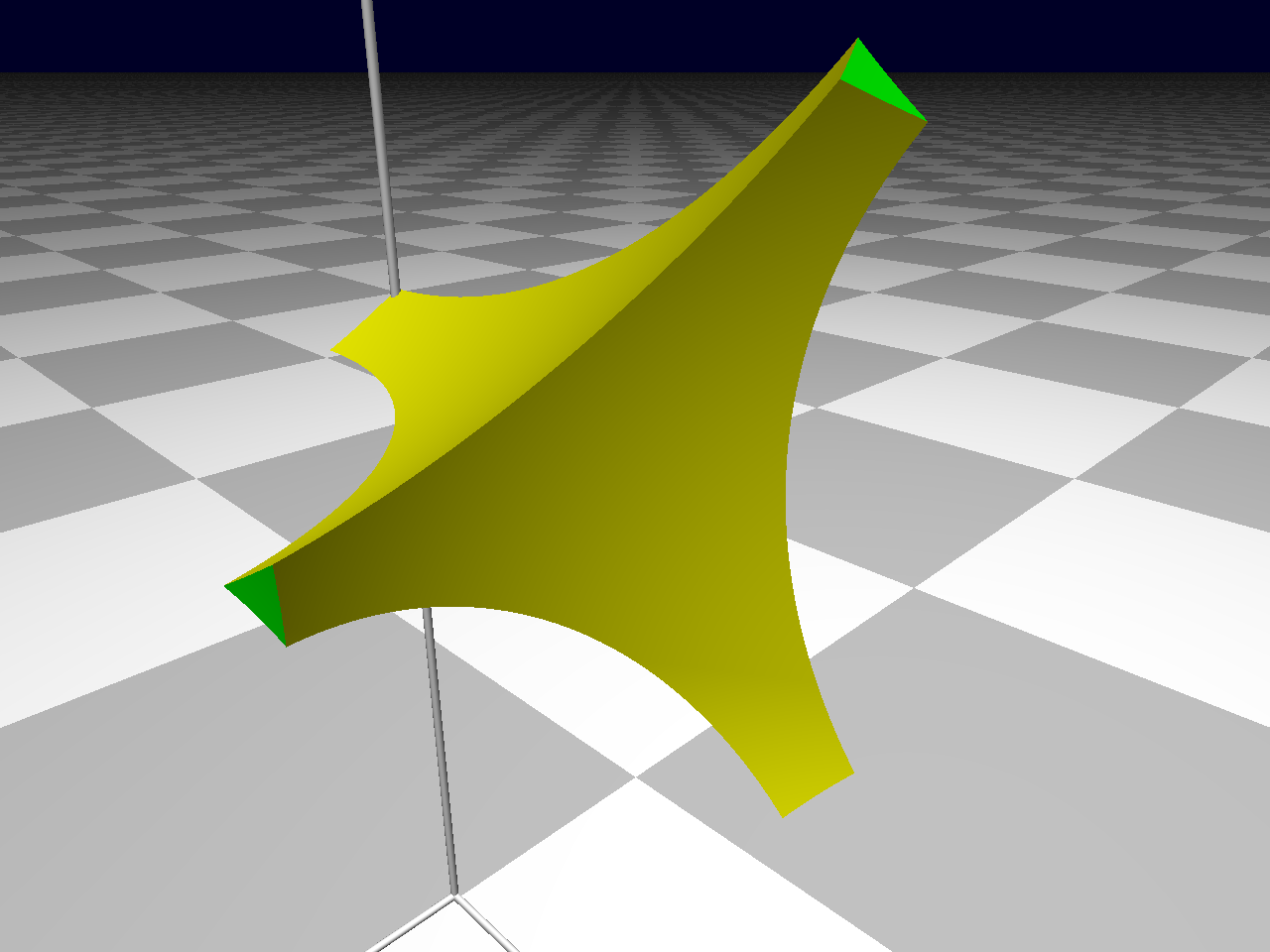}
\caption{The conformal geometry of a component of the conformal boundary}
\label{dd1 geom}
\end{figure}

Using a computer it is possible to create a picture showing how $C$ sits inside the sphere at infinity $\mathbb{S}^3_\infty$, thereby displaying its conformal geometry.  Such a picture is shown in figure \ref{dd1 geom}, where the parameter value $t = 0.8$.  Figures \ref{dd1 figure} and \ref{dd1 geom} have been arranged so that the face labels of figure \ref{dd1 figure} correspond to the face labels of figure \ref{dd1 geom}.

Let us consider the subgroup $B \subset \rho_t(\Gamma_{22})$ generated by reflections in the 
$\m{0}$, $\m{2}$,$\m{4}$, $\m{6}$,$\p{1}$, $\p{3}$, $\p{5}$, and $\p{7}$ walls.  These are the walls that determine the end $E_G$ for $1 > t > \sqrt{3/5}$.  On the sphere at infinity they determine the faces of
the corresponding component $C$ of $\partial \overline{\mathcal{F}}_{t}$.  For this reason we will refer to $B$ as a boundary subgroup.  

We can consider the inclusion map of the boundary subgroup $B$ into $\Isom (\Hf)$ as a point in the representation variety $\Hom (B, \Isom(\Hf))$ and ask whether or not it has nontrivial deformations.
In considering this question we will assume that the walls that intersect in angle $\theta$ as well as those
that intersect in angle $\pi/2$ continue to do so under deformation.  For angle $\pi/2$ this is forced by relations in
the abstract group $B$ but for angle $\theta$ there are corresponding relations only when the angle is
a rational multiple of $\pi$.  Assuming that the angles are preserved in all cases leads to a more uniform treatment for all values of $\theta$.

We can quickly find a $4$-dimensional space of local deformations of $B$ as follows:
Fix the odd, positive walls. This ensures that the angle $\theta$ is preserved.
Consider the vectors $(\frac{r+1}{\sqrt{2}}, r, r,r,-t)$.  When $r = 1$ this is the space-like
vector corresponding to the wall $W_{\m{0}}$ at time $t$.  For all values of $r$ it will remain
orthogonal to the space-like vectors corresponding to the $\p{1}$, $\p{3}$, $\p{5}$ walls and it
will still be a space-like vector for values of $r$ near $1$.  There are no further constraints on the relationships between the generators so this gives a $1$-dimensional
family of small deformations.  Similar, independent deformations are possible 
for the other even, negative walls.  This gives a $4$-dimensional space of nontrivial local deformations.
 
It is possible to compute the dimension of the Zariski tangent space of a scheme corresponding to
this representation variety (with the angle constraints) and see that it is
$14$-dimensional.  We will not do that here.  Note that $\Isom(\Hf)$ is $10$-dimensional.
Using the argument that will be provided
in the proof of Theorem \ref{infinitesimal thm} in Section \ref{computing infinitesimally}, this computation, together with the existence of a smooth $4$-dimensional family of nontrivial deformations,
implies that the representation variety is smooth at this point and that these are all the small deformations, modulo conjugation, near this point.

We now specialize to the values $t = t_n$ for $n = 2m$, $m$ an integer at least $4$.  We will
assume that $t$ and $n$ are of this form for the remainder of this section.  For these
values, the angle $\theta$ equals $\pi/m$ and all of the dihedral angles of the finite-sided
polytope $\mathcal{F}_{t_n}$ are of the form $\pi/k$ for some integer $k$.  Poincar{\'e}'s lemma 
(Theorem \ref{poincare lem v2}) implies that $\mathcal{F}_{t_n}$ is a fundamental domain for the discrete hyperbolic reflection group $\rho_{t_n} (\Gamma_{22}) = \Lambda_n$.  Since these polytopes have infinite volume, the quotient space $\Hf / \Lambda_n$ has infinite volume.  However, because it has a finite-sided fundamental domain, its convex core has finite volume.  (See \cite{Bow} for a general discussion of the relationship between fundamental domains and the volume of the convex core.  The required implication in our case is not so difficult because  all the cusps are of maximal rank, but the general situation is quite subtle.) 

The group acts properly discontinuously on an open region $\Omega$ in the sphere at infinity.  The convex core has two boundary components (orbifold-)homeomorphic to the boundary components of 
$(\Hf \cup \Omega)/\Lambda_n$, viewed as an orbifold with boundary. 
Neither boundary component is totally geodesic.  If one were totally geodesic, reflection in the walls determining the corresponding end of $\mathcal{F}_{t_n}$ would preserve a hyperplane in $\Hf$.  We saw at the beginning of
this section that no such hyperplane exists for the values of $t$ under consideration.

We state this as a theorem.

\begin{prop} \label{boundary never geodesic}
For even $n\ge 8$ the convex core of $\Lambda_n$ does not have totally geodesic boundary.
\end{prop}

We cannot say much more about the boundary of the convex hull of $\Lambda_n$; its geometry 
remains a mystery.

For the values $t = t_n$ the boundary subgroup $B$ is also a discrete hyperbolic reflection group.
It is generated by reflections in $8$ walls that bound a convex region which is a fundamental
domain for the group. It also is geometrically finite with a finite volume convex core.
Unlike the convex core of $\Lambda_n$, which has cusps, the convex core of $B$ is compact.
This follows from the fact that none of the walls bounding the fundamental domain are tangent.
They are either disjoint or intersect.  There are no parabolic elements in $B$.  Such a group
is called convex cocompact.  By compactness, this property holds under small deformations.

We have seen that there are always small deformations of the boundary group $B$ for the values of
$t$ under consideration.  The roll symmetry induces an automorphism of $\Lambda_n$ taking $B$ to 
the boundary subgroup of the other boundary component so that subgroup has the same properties. 
Again, we record these cases as a theorem.

\begin{prop} \label{boundary groups prop}
The discrete group $\Lambda_n < \Isom(\Hf)$ for $n \in \{ 8,10,12,\ldots\}$ has two (conjugacy classes of) boundary groups, one for each end of the orbifold $\mathbb{H}^4 / \Lambda_n$.  Each boundary group is convex cocompact, and has nontrivial convex cocompact deformations as a subgroup of $\Isom(\Hf)$.
\end{prop}

We will see in Section \ref{computing infinitesimally} that the groups $\Lambda_n$ are, nonetheless, infinitesimally rigid.  In other words, none of the deformations of the boundary groups extend over the larger group.



\section{The disappearing domain of discontinuity} \label{disappear} 

This section continues the analysis of the geometry of the polytope $\mathcal{F}_t$, focusing on the changes that occur in the conformal  boundary at infinity.  We first explain what the expected behavior would be, based on geometric, rather than conformal models.  Then we produce a sequence of pictures describing the intersection patterns at infinity analogous to those in Section \ref{viewed in the sphere at infinity}.  These confirm the expected behavior.  However, since they were produced by computer using floating point  computation, they should not be considered as a rigorous proof.  A rigorous computational  proof would be possible but lengthy;  primarily one must show that there are no new $k$-fold intersections between the hyperplanes determining the polytope other than those described.  However, by using the symmetry of the deformation and analyzing the reflection groups extended by these symmetries we will be able more easily to prove these results rigorously in Section \ref{manual lattice}.

Recall that  for $t$ values between $t_6$ and $1$ a component of the intersection $\mathcal{F}_t \cap \partial \Hf$ looks like the truncated tetrahedron of figures \ref{dd1 figure} and \ref{dd1 geom}.   This figure has a realization as a polyhedron in $\mathbb{H}^3$  with the same angles.  Let us analyze how such
polyhedra would change as the angles change, referring  to figure \ref{dd1 geom} as if it were a hyperbolic model.
 As $t$ decreases from $1$ to $t_6$, the dihedral angle along the red edges 
 increases to $\pi/3$.  Consider one of the triangular faces of figure \ref{dd1 figure}.  It is a hyperbolic  equilateral triangle with interior angle increasing to $\pi/3$.  A hyperbolic equilateral triangle must have interior angle strictly smaller than $\pi/3$. As $t$ decreases to $t_6$, the triangular faces shrink as their interior angle increases to $\pi/3$.  At $t_6$, the triangular faces have disappeared, leaving a regular  tetrahedron with dihedral angle $\pi/3$.

We will see the same behavior occur in the conformal boundary at infinity.  In particular, the triples of hyperplanes which form the $3$ edges  of a triangular face in figure \ref{dd1 figure} intersect in a single point at infinity when $t =  t_6 = \sqrt{3/5}$.  At this stage,  the analyses of Sections \ref{fundamental domain} and \ref{miscellaneous section} are no longer valid, which is why we restricted attention to values
of $t$ between $t_6$ and $1$.

Let us now focus on the expected geometry at infinity of the polytope $\mathcal{F}_t$ for $0.378 \approx \sqrt{1/7} = t_3 < t \le t_6$.  (We will explain the lower bound later.)  We will pay particular attention to the parameter values $t_6$ and $t_4$, where $\mathcal{F}_t$ is a fundamental domain for the discrete groups $\Lambda_6$ and $\Lambda_4$.  

As $t$ decreases from $t_6$, the dihedral angle of the tetrahedron increases from $\pi/3$.  There is a regular  hyperbolic tetrahedron with these angles; all of its vertices are finite.  As the angles increase,  the tetrahedron shrinks.  At $t = \sqrt{1/2}$ the dihedral angle is $\arccos(1/3)$, which is the dihedral angle of a regular Euclidean tetrahedron, and the tetrahedron will collapse to a point.   We will see that at this value of 
$t$ the conformal boundary components will disappear.   At $t = \sqrt{1/2}$ the two infinite volume ends of $\mathcal{F}_t$ will collapse to  cusps with cross-sections equal to a regular Euclidean tetrahedron.
In particular,  at $t=\sqrt{1/2}$ the polytope has finite volume!

The cusps will become  finite vertices of $\mathcal{F}_t$ as $t$ decreases past $\sqrt{1/2}$.  A small neighborhood of the vertices in $\mathcal{F}_t$ will be isometric to a hyperbolic cone on a spherical regular tetrahedron.  At $t_4$ this spherical tetrahedron will have all orthogonal dihedral angles.  The corresponding reflection group $\Lambda_4$ will generate a finite covolume lattice.  Finally, at $t_3$ this spherical tetrahedron has dihedral angles $2\pi/3$ and we will see in Section \ref{manual lattice} that the reflection group is again discrete with finite covolume.

We now present a series of computer-drawn pictures which provide strong evidence for the fact that the expected behavior actually occurs.   Recall that via the symmetries, the letter walls of $\mathcal{F}_t$ are all isometric.  Similarly, the walls of the positive octet are all isometric, as are the walls of the negative octet.  This is true for all $t>0$.  We therefore need only examine one wall from each octet.  As in Section \ref{fundamental domain}, we will examine the walls corresponding to $\l{A}$, $\m{0}$, and $\p{3}$.  For each we will show a series of figures similar to figures \ref{circles2 figure}, \ref{slicewith1.fig} and \ref{slicewith3.fig}.  Let us begin with $\l{A}$.

\begin{figure}[ht]
\centering
\includegraphics[scale=0.5]{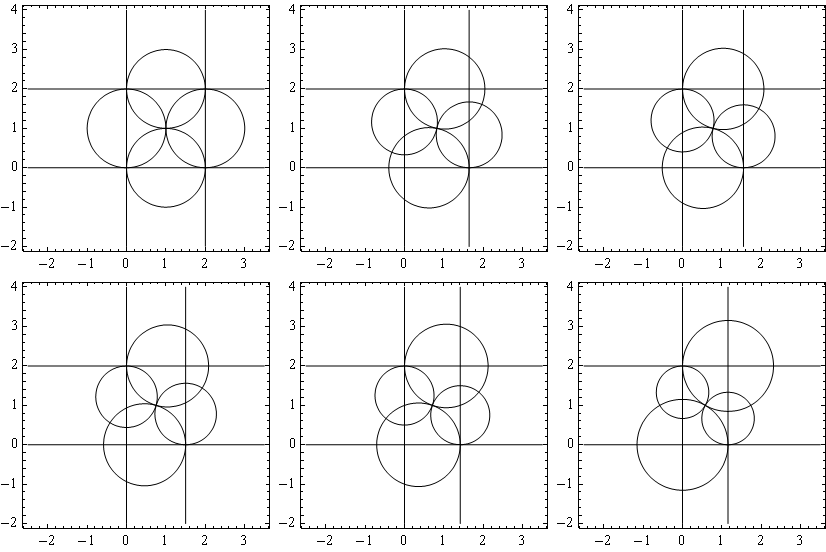}
\caption{Slicing along wall $\l{A}$}
\label{A_wall_figure}
\end{figure}

\begin{figure}[ht]
\centering
\includegraphics[scale=0.5]{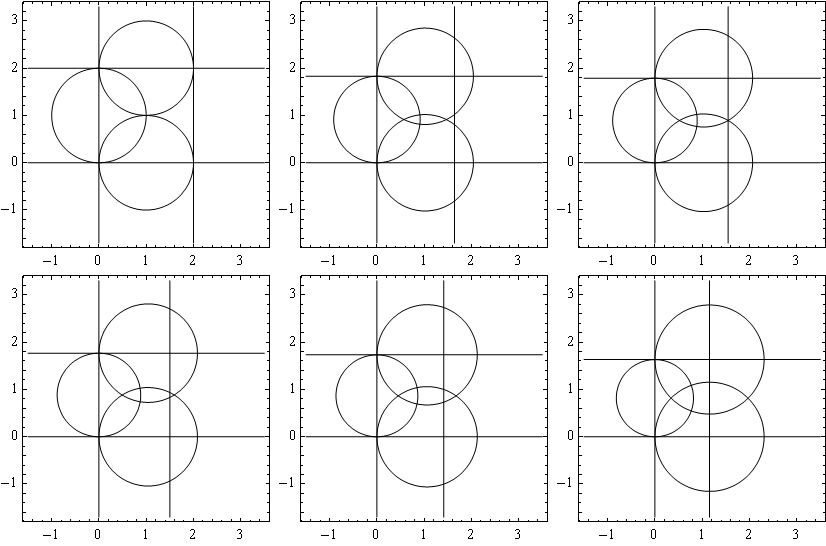}
\caption{Slicing along wall $\m{0}$}
\label{0-_wall_figure}
\end{figure}

\begin{figure}[ht]
\centering
\includegraphics[scale=0.5]{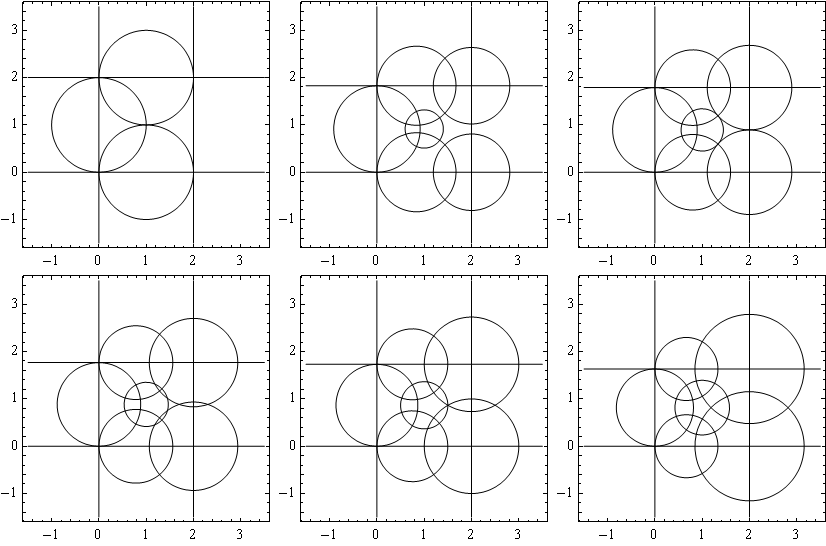}
\caption{Slicing along wall $\p{3}$}
\label{3+_wall_figure}
\end{figure}

Observe figure \ref{A_wall_figure}.  Its six pictures form a rough animation proceeding from left to right and then top down.  The first picture is a copy of figure \ref{circles1 figure}.  In it the plane represents $\partial H_{\l{A}}$.  (Recall $H_{\l{A}}$ is the hyperplane containing the wall $W_{\l{A}}$.)  The conformal circles represent the intersections of $\partial H_{\l{A}}$ with the boundaries of the other hyperplanes of the arrangement when $t=1$.  These circles (and lines) should be labeled as in figure \ref{circles1 figure}.  The other five pictures show the same intersections as the parameter $t$ decreases.  The six pictures show the parameter values $\{1, \, 0.82,\,  t_6 = \sqrt{3/5}, \, 0.75, \,\sqrt{1/2},t_4= \sqrt{1/3} \}$.  (These values were selected uniformly for figures \ref{A_wall_figure}, \ref{0-_wall_figure}, and \ref{3+_wall_figure} to highlight the interesting changes through the deformation.)  We remark that there are no other hyperplanes of the arrangement intersecting $H_{\l{A}}$ for this range of parameter values.  The third and sixth pictures correspond to the discrete groups $\Lambda_6$ and $\Lambda_4$ respectively.  There are no combinatorial changes in this range of parameter values, and we are ready to proceed to $\m{0}$.

Next observe figure \ref{0-_wall_figure}.  Its six pictures form another animation which should be read in the same way as figure \ref{A_wall_figure}.  The second picture is nearly a copy of figure \ref{slicewith1.fig}.  In it the plane represents $\partial H_{\m{0}}$ and the circles represent its intersections with the other hyperplanes when $t=0.82$.  (Figure \ref{slicewith1.fig} shows the same thing for the parameter value $t=0.8$.)  As before, the circles should be labeled as in figure \ref{slicewith1.fig}.  The other five pictures show the same intersections, and should be labeled in the natural way.  The pictures show the same parameter values $\{1, 0.82, t_6 = \sqrt{3/5}, 0.75, \sqrt{1/2}, t_4= \sqrt{1/3} \}$.  There are no other hyperplanes of the arrangement intersecting $H_{\m{0}}$ for this range of parameter values.

Unlike the previous case, a combinatorial change is visible in the series of figure \ref{0-_wall_figure}.  Recall the magenta triangle visible in figure \ref{slicewith1.fig}.  This triangle indicated an intersection of $\mathcal{F}_t$ with the sphere at infinity.  This triangle is visible in the second picture of figure \ref{0-_wall_figure}.  In the first picture it is expanded to a cusped triangle.  More interestingly, in the third picture this triangle has collapsed to a point.  For $t_6 < t \le 1$ the hyperplanes corresponding to $\m{0}$, $\p{1}$, $\p{3}$, and $\p{5}$ do not have a common intersection at infinity.  Then at $t_6$ they come together to form a new ideal vertex of $\mathcal{F}_{t_6}$.  For smaller values of $t$ these $4$ hyperplanes intersect in a finite point that is a new finite vertex of $\mathcal{F}_t$.    The wall $W_{\m{0}}$ intersects the sphere at infinity only at points corresponding to cusps which have existed througout the deformation.
This confirms the first transition suggested by our geometrical analysis.

Finally we consider the pictures of figure \ref{3+_wall_figure}.  As before, its six pictures form an animation which should be read left to right and then top down.  The second picture is nearly a copy of figure \ref{slicewith3.fig}.  The plane represents $\partial H_{\p{3}}$ and the circles represent its intersections with the other hyperplanes.  The circles of the second picture should be labeled as in \ref{slicewith3.fig}, and these labelings should be extended to the other pictures in the natural way.  The same parameter values are shown as in the previous two figures.  There are no other hyperplanes of the arrangement intersecting $H_{\p{3}}$ for this range of parameter values.

The combinatorial changes are most drastic in this case.  For $t_6 < t < 1$, the wall $W_{\p{3}}$ intersects the sphere at infinity in a hexagon, visible in the second picture and colored magenta in figure \ref{slicewith3.fig}.  At $t_6$, three edges of this hexagon collapse to points, and the intersection becomes a cusped triangle.  These collapsing edges each correspond to an edge of one of  the disappearing triangles  at infinity for the negative walls.  The disappearing triangle in $W_{\m{0}}$ was described above and each of the negative walls undergoes the same process.  For $\sqrt{1/2} < t < t_6$ the triangle $W_{\p{3}} \cap \partial \Hf$ shrinks as its interior angle increases.  Finally at $t= \sqrt{1/2}$ the triangle collapses to a point.  This is visible in the fifth picture of figure \ref{3+_wall_figure}.  The hyperplanes corresponding to $\p{1}$, $\p{3}$, $\p{5}$, and $\p{7}$ have come together to intersect in single point at infinity, creating a new ideal vertex of $\mathcal{F}_{\sqrt{1/2}}$.  For $t < \sqrt{1/2}$ this vertex becomes a finite vertex and  the wall $W_{\p{3}}$ intersects the sphere at infinity only in points corresponding to cusps which have existed throughout the deformation.  This implies that for  $t \leq \sqrt{1/2}$, $\mathcal{F}_t$  will be a \emph{finite volume}  polytope in $\Hf$.

In particular, if this analysis  were rigorous, it would imply  that $\Lambda_3$ and $\Lambda_4$
are lattices.  A method for proving this rigorously (without a computer) will be explained in Section \ref{manual lattice}. 

 Beginning at  $t_3 = \sqrt{1/7}$ there are new pairwise intersections between the hyperplanes determined by the space-like vectors in table \ref{qt table}.  Although these intersections do not actually affect the geometry of the polytope determined by these vectors, this is difficult to see using the analogous pictures for those values of $t$.  Thus, we will postpone the analysis of these values until Section \ref{manual lattice}.



\section{Infinitesimal deformations of the letters} \label{inf letter} 

At this point of the paper we mentally back up a few steps and re-examine,  at the infinitesimal level, the deformation given in table \ref{qt table}.  Recall in Section \ref{defining the deformation 2} we showed that the deformation of $\Gamma_{22}$ given in table \ref{qt table} is unique up to conjugation \emph{among deformations preserving the roll symmetry}.  The goal of this section and the next is to show this symmetry hypothesis is not necessary; the deformation is simply unique up to conjugation.  

We begin by returning to the following set of six walls from the letter octet:

\[
\begin{array}{lll}
\l{A} = \left( 1, \sqrt{2},0,0,0 \right) & \qquad & \l{B} = \left(1,0,\sqrt{2},0,0\right) \\
\l{C} = \left( 1,0,0,\sqrt{2},0 \right) & \qquad & \l{D} = \left(1,0,0,-\sqrt{2},0\right) \\
\l{E} = \left( 1,0,-\sqrt{2},0,0\right) & \qquad & \l{F} = \left(1,-\sqrt{2},0,0,0\right)
\end{array}
\]

As in Section \ref{defining the deformation 2}, these six vectors describe walls $W_{\l{A}}$ through $W_{\l{F}}$ in $\mathbb{H}^4$.  Let $V$ be the copy of $\mathbb{H}^3$ orthogonal to the space-like vector $e_4 = (0,0,0,0,1)$.  Note the six letter walls $W_{\l{A}}$ through $W_{\l{F}}$ all intersect $V$ orthogonally.  On the sphere at infinity  the intersections $W_{\l{A}} \cap \partial V$ through $W_{\l{F}} \cap \partial V$ will form the conformal arrangement of figure \ref{diag_slice}.  

Suppose $\l{A}(t)$ through $\l{F}(t)$ are paths of space-like vectors defining one-parameter deformations of walls $W_{\l{A}}$ through $W_{\l{F}}$.  Assume that $\l{A}(1)=\l{A}$, $\l{B}(1) = \l{B}$, etc..  Moreover, assume that the deformation preserves the tangencies of the original walls, since this will be implied by Lemma \ref{cube deformations}  applied to deformations of $\Gamma_{22}$.  Finally, we define 
\[ \ld{A} = (a_0, a_1,a_2,a_3,a_4) = \l{A}'(1), \quad \ld{B} = (b_0, b_1,b_2,b_3,b_4) = \l{B}'(1), \quad \text{etc..} \]

As we are only interested in the deformations of the walls $W_{\l{A}}$ through $W_{\l{F}}$, the norms of the space-like vectors $\l{A}(t)$ through $\l{F}(t)$ are irrelevant, and we are free to assume without loss of generality that 
\[ \| \l{A}(t) \| = \ldots = \| \l{F}(t) \| = 1 \]
for all $t$.  (Note that the initial vectors $\l{A}(1)$ through $\l{F}(1)$ are unit vectors.)  
 We will also conjugate our paths $\l{A}(t)$ through $ \l{F}(t)$ into a simple form by a path of Minkowski isometries.  Thus we are choosing a particular slice of the set of conjugacy classes in the deformation space of the letter walls (and of $\Gamma_{22}$).  This slice will be chosen to be a real algebraic variety.

This is easiest to describe in $\partial \mathbb{H}^4$, where the $1$-parameter families $W_{\l{A}(t)} - W_{\l{F}(t)}$ become the $1$-parameter families of spheres $S_{\l{A}(t)} - S_{\l{F}(t)}$, and the conjugating Minkowski isometries become conformal automorphisms of $\partial \mathbb{H}^4$ which we will identify with $\mathbb{R}^3 \cup \infty$ via stereographic projection.   First, we may assume the point of tangency between $S_{\l{A}(t)}$ and $S_{\l{B}(t)}$ is fixed and equal to the point at infinity for all $t$.    Second, conjugating by similarities of $\mathbb{R}^3$, we may assume that $S_{\l{A}(t)} = S_{\l{A}(1)}$ and $S_{\l{B}(t)} = S_{\l{B}(1)}$ are equal to the planes $z = 0$ and $z=2$, respectively, for all $t$.  This determines the configuration up 
to a Euclidean isometry of the $xy$-plane.  Since the $\l{C}$ wall is tangent to the $\l{A}$ and $\l{B}$ walls for all $t$,  we can arrange, by a translation, that 
$S_{\l{C}(t)} = S_{\l{C}(1)}$ for all $t$.   Finally, by a rotation of the $xy$-plane fixing the point of tangency of $S_{\l{C}(t)}$ with 
$ S_{\l{A}(t)}$ we can assume that the point of tangency of $S_{\l{D}(t)}$ and $S_{\l{A}(t)}$ lie on the same line for all $t$.  Specifically, we may assume that $S_{\l{D}(t)}$, hence the entire quadruple of spheres $S_{\l{A}(t)} - S_{\l{D}(t)}$, remains orthogonal to $\partial V \subset \partial \mathbb{H}^4$ for all $t$. Pictorially, these normalizations tell us the the arrangement of the three spheres $S_{\l{A}(t)} - S_{\l{C}(t)}$ remains as shown in figure \ref{diag_slice} for all $t$, and the sphere $S_{\l{D}(t)}$ moves only horizontally.

Since we have fixed the space-like vectors determining the $\l{A}, \l{B}, \l{C}$ walls, the resulting
variety is determined by the variation of the remaining space-like vectors.  The conditions placed
on these involve tangency and orthogonality relations, along with the norm-squared normalization.
The final normalization condition on $S_{\l{D}(t)}$ is equivalent to the condition that the vector
determining the $\l{D}$ wall is orthogonal to $e_4$.  These conditions are equivalent to various 
Minkowski dot products remaining constant,  so all of the conditions are polynomial
and the resulting variety is real algebraic.  We will now compute the Zariski tangent space 
of the scheme given by these equations.

Taking  the time derivatives of the fact that the $\l{A}, \l{B}, \l{C}$ walls are fixed yields the relations
\[ \ld{A} = 0, \quad \ld{B} = 0, \quad \ld{C} = 0 \]
and the norm-squared normalizations yield the relations
\[ \left( \l{D}, \ld{D} \right) = \left( \l{E}, \ld{E} \right)  = \left( \l{F}, \ld{F} \right) = 0. \]
The fact that the $\l{D}$ vector is always perpendicular to $e_4$ implies that
\[ d_4 = \left( \ld{D}, e_4 \right) = 0. \]

The sphere $S_{\l{D}(t)}$ is tangent to $S_{\l{A}(t)}$ and $S_{\l{B}(t)}$.  From these tangencies and the fact that $\ld{A}=\ld{B}=0$, we conclude that
\[ \left( \l{A},\ld{D} \right) = \left( \l{B},\ld{D}\right) = 0. \]
Together with  $\left( \l{D}, \ld{D} \right) = 0$ and $d_4=0$, these facts imply that
\[ \ld{D} = \left( \sqrt{2} \delta, \delta,\delta,- \delta,0 \right) \quad \text{for some } \delta.\]
Performing a similar analysis for $\ld{E}$, we use the relations 
\[ \left( \l{A},\ld{E} \right) = \left( \l{C}, \ld{E} \right) = \left( \l{E},\ld{E} \right) = 0\]
to conclude that
\[ \ld{E} = \left( \sqrt{2} \varepsilon',\varepsilon',-\varepsilon',\varepsilon',\varepsilon \right) \quad
    \text{for some } \varepsilon, \varepsilon'.\]
Repeating this procedure for $\ld{F}$ shows that
\[ \ld{F} = \left( \sqrt{2} \phi', - \phi',\phi',\phi',\phi \right) \quad \text{for some } \phi, \phi'. \]

We have not yet used the tangency relation
\[ \left( \ld{E}, \l{D} \right) + \left( \l{E}, \ld{D} \right) = 0.\]
Using the above expressions for $\ld{D}, \l{D}, \ld{E}$, and $\l{E}$ this now expands to the equation
\[ - \sqrt{2} \varepsilon' - \sqrt{2} \varepsilon' - \sqrt{2} \delta - \sqrt{2} \delta = 0,\]
implying that $\varepsilon'=-\delta$.  By replacing $\l{E}$ with $\l{F}$ we can similarly conclude that $\phi' = - \delta$.
With this the tangency relation
\[ \left( \ld{E},\l{F} \right) + \left( \l{E}, \ld{F} \right) = 0 \]
implies $\delta = 0$.  We have therefore shown that 
\[
\begin{array}{c}
\ld{A} = \ld{B} = \ld{C} = \ld{D} = 0, \\ 
 \ld{E} = \left( 0,0,0,0,\varepsilon\right), \quad \text{and  } \ld{F} = \left( 0,0,0,0,\phi \right) .
\end{array}
\]

Let us restate this analysis as a proposition.  

\begin{prop} \label{inf letter walls prop}
Let $\ld{A}, \ld{B}, \ldots, \ld{F}$ be an infinitesimal deformation of the space-like vectors $\l{A}$ through $\l{F}$.  Moreover, let us assume that this infinitesimal deformation preserves the tangencies between the walls $W_{\l{A}}$ through $W_{\l{F}}$ and is infinitesimally norm preserving.  It is possible to conjugate this infinitesimal deformation so that 
\[ \ld{A} = \ld{B} = \ld{C}  = \ld{D} = 0 , \]
and the final pair of vectors satisfy the equations
\[ \ld{E} = \left(0,0,0,0,\varepsilon\right) \quad \text{and} \quad \ld{F} =\left(0,0,0,0,\phi\right).\]
\end{prop}

One can show that this variety is, in fact,  $2$-dimensional and smooth at this initial point, but
since we will not need this fact we will not prove it here.    Qualitatively, this analysis shows that infinitesimal deformations of the arrangement in figure \ref{diag_slice} only involve the spheres $S_{\l{F}}$ and $S_{\l{E}}$ moving in the direction orthogonal to the page.  However,  we note that we have \emph{not} shown (and it is not true) that $\l{D}(t) = \l{D}(0)$ for all $t$, merely that $\l{D}(t)$ is constant to first order.


\section{computing infinitesimally} \label{computing infinitesimally}

This section will build on the computations of Section \ref{inf letter}.  The goal is to prove that the deformation of $\Gamma_{22}$ described in table \ref{qt table} is unique at the infinitesimal level.  Specifically, we will show that at any point $\rho_t$ of the deformation described in Section \ref{defining the deformation 2}, the Zariski tangent space of the deformation space is exactly $1$ dimensional after modding out by conjugation.  From this we can conclude that there are no new infinitesimal deformations of $\rho_t$ other than those coming from infinitesimally varying the $t$ parameter.   An application of the implicit function theorem will imply that the deformation variety of $\Gamma_{22}$,  
where we take the slice of conjugacy classes described in the previous section, is $1$-dimesional and smooth at all of the points $\rho_t, t >0$.  We emphasize that no symmetry assumptions are placed on the deformations in this section.  

Fix a value $t_0 > 0$ and consider the representation $\rho_{t_0}$ given by table \ref{qt table}.  Given an infinitesimal deformation of $\rho_{t_0}$,  denote the infinitesimal deformations of the generating space-like vectors by $\pd{0}, \md{0}, \ldots, \md{7}, \ld{A}, \ldots , \ld{F}$.  The goal is to prove that this deformation is equivalent to an infinitesimal variation of the parameter $t$.  (By equivalent we mean up to infinitesimal scaling and conjugation in $O(1,4)$.  Infinitesimal scaling means, for example, adding multiples of $\p{0}$  to $\pd{0}$.)

We would like to apply Proposition \ref{inf letter walls prop} to $\ld{A}$ through $\ld{F}$.  We recall that,  for any value of $t_0>0$,  the rank $3$ cusps of $\Gamma_{22}$ remain rank $3$ cusps in $\rho_{t_0}(\Gamma_{22})$ and the space-like vectors corresponding to the letter walls are unchanged.    Furthermore, by Lemma \ref {cube deformations},  under any small deformation of 
$\rho_{t_0}(\Gamma_{22})$ the rank $3$ cusps remain rank $3$.   This
implies, in particular, that any infinitesimal deformation will infinitesimally preserve the tangencies between the lettered walls.  Thus we can apply Proposition \ref{inf letter walls prop}.  Making the constant norm and conjugacy  normalizations of the previous section  we conclude that  the infinitesimal deformations of the letters have the form
\[ \ld{A} = \ld{B} = \ld{C}  = \ld{D} = 0 , \]
and
\[ \ld{E} = \left(0,0,0,0,\varepsilon\right) \quad \text{and} \quad \ld{F} =\left(0,0,0,0,\phi\right).\]

After this normalization we have no freedom to conjugate our remaining generators further.  However, we may assume  that the numbered vectors $\pd{0}, \md{0}, \ldots, \md{7}$ describe a norm preserving infinitesimal deformation.  In other words
\[ \left( \p{0}, \pd{0} \right) = \left( \m{0}, \md{0} \right) = \ldots = \left( \m{7}, \md{7} \right) = 0.\]
  
In Section \ref{defining the deformation 2} the only use of the assumption that the deformation preserved the roll symmetry was to conclude that all of the letter walls stayed constant.  Once we established this property, the proof that our family of deformations is the unique one fixing the letter walls was straight forward.  The main step in this section is to prove that the letter walls remain infinitesimally fixed, without the symmetry assumption.  It remains to show  that $\ld{E}$ and $\ld{F}$ are also $0$.  This requires a few intermediate computational lemmas.

\begin{lem} \label{computation lem 1}
There exists a constant $a$ such that
\begin{eqnarray*}
\pd{0} = \frac{-a}{t_0^2} \cdot \left( \sqrt{2}, 1, 1, 1, t_0 \right) & \quad & 
\md{0} = a \cdot \left( \sqrt{2}, 1, 1, 1, -1/t_0 \right) \\
\pd{3} = \frac{-a}{t_0^2} \cdot \left( \sqrt{2}, 1, 1, -1, -t_0 \right) & \quad & 
\md{3} = a \cdot \left( \sqrt{2}, 1,1,-1, 1/t_0 \right)
\end{eqnarray*}
\end{lem}
\begin{pf}
Begin by analysing $\pd{0}$.  We know $\left( \pd{0}, \p{0} \right) = 0$.  As $\ld{A} = \ld{B} = \ld{C} = 0$, derivatives of the relations
\[ \left( \p{0}, \l{A} \right) = \left( \p{0}, \l{B} \right) = \left( \p{0}, \l{C} \right) = 0 \]
yield
\[ \left( \pd{0}, \l{A} \right) = \left( \pd{0},  \l{B} \right) = \left( \pd{0}, \l{C} \right) = 0 .\]
From this we conclude that
\[ \pd{0} = c_0^+ \cdot \left( \sqrt{2}, 1, 1, 1, t_0 \right).\]

Using the relation $\ld{D} = 0$ and arguing similarly for $\m{0}$, $\p{3}$, and $\m{3}$ yields the equations
\begin{eqnarray*}
\pd{0} = c_0^+ \cdot \left( \sqrt{2}, 1, 1, 1, t_0 \right) & \quad & 
\md{0} = c_0^-  \cdot \left( \sqrt{2}, 1, 1, 1, -1/t_0 \right) \\
\pd{3} = c_3^+ \cdot \left( \sqrt{2}, 1, 1, -1, -t_0 \right) & \quad & 
\md{3} = c_3^-  \cdot \left( \sqrt{2}, 1,1,-1, 1/t_0 \right).
\end{eqnarray*}

It remains to relate the constants $c_0^+$, $c_0^-$, $c_3^+$, and $c_3^-$.  This we do using derivatives of the orthogonality relations
\[ \left( \m{0}, \p{0} \right) = \left( \m{0}, \p{3} \right) = \left( \p{3}, \m{3} \right) = 0.\]
These show that
\[ c_0^+ = -\, \frac{c_0^-}{t_0^2}, \qquad c_3^+ = -\, \frac{c_0^-}{t_0^2}, \qquad c_3^- = -\, c_3^+ \cdot t_0^2 = c_0^-.\]
Setting $a = c_0^-$ completes the proof of the lemma.
\end{pf}

The main fact driving the proof of Lemma \ref{computation lem 1} is that $\ld{A}$, $\ld{B}$, $\ld{C}$, and $\ld{D}$ are all $0$.  Attempting to emulate the above computation for the quartet $\p{1}$, $\m{1}$, $\p{2}$, and $\m{2}$ will run into complications because we do not know that $\ld{E}=0$.  However, repeating the above steps will prove the following slightly weaker lemma.

\begin{lem} \label{computation lem 2}
Let $a$ be the constant from Lemma \ref{computation lem 1}.  There exist constants $c_1^+$, $c_1^-$, $c_2^+$, and $c_2^-$ such that
\begin{eqnarray}
\pd{1} &=& \left( \sqrt{2} c_1^+, c_1^+, a/t_0^2, c_1^+, a/t_0 \right)  \label{cl1} \\
\md{1} &=& \left( \sqrt{2} c_1^-, c_1^-, - a, c_1^-, a/t_0 \right)	\label{cl2} \\
\pd{2} &=& \left( \sqrt{2} c_2^+, c_2^+, a/t_0^2, -c_2^+, -a / t_0 \right) \label{cl3} \\
\md{2} &=& \left( \sqrt{2} c_2^-, c_2^-, -a, -c_2^-, -a/t_0 \right)	\label{cl4}
\end{eqnarray}
\end{lem}

\begin{pf}
Begin with $\pd{1}$.  Apply the relations
\[ \left( \pd{1}, \p{1} \right) = \left( \pd{1}, \l{A} \right) = \left( \pd{1}, \l{C} \right) = 0 \]
to conclude
\[ \pd{1} = \left( \sqrt{2}, c_1^+, c_1^+, d_1^+, c_1^+, d_1^+ \, t_0 \right), \]
for some constant $d_1^+$.  The orthogonality condition $\left( \p{1}, \m{0} \right) = 0$ yields
\[ \left( \pd{1}, \m{0} \right) + \left( \p{1}, \md{0} \right) = 0.\]
Expanding this relation yields $d_1^+ = a/t_0^2$, as desired.  This establishes equation \ref{cl1}.

Arguing similarly for $\m{1}$, use the equtions
\[ \left( \md{1}, \l{A} \right) = \left( \md{1}, \l{C} \right) = \left( \md{1}, \m{1} \right) = \left( \md{1}, \p{0} \right) + \left( \m{1}, \pd{0} \right) = 0 \]
to deduce equation \ref{cl2}.  Equation \ref{cl3} follows from 
\[ \left( \pd{2}, \l{A} \right) = \left( \pd{2}, \l{D} \right) = \left( \pd{2}, \p{2} \right) = \left( \pd{2}, \m{3} \right) + \left( \p{2}, \md{3} \right) = 0. \]
Finally, equation \ref{cl4} follows from 
\[ \left( \md{2}, \l{A} \right) = \left( \md{2}, \l{D} \right) = \left( \md{2}, \m{2} \right) = \left( \md{2}, \p{3} \right) + \left( \m{2}, \pd{3} \right) = 0.\]
\end{pf}

Armed with equations \ref{cl1} - \ref{cl4}, only a couple additional computations are needed to show that $\ld{E}=0$.  Specifically,
\[ \left( \md{1}, \p{2} \right) + \left( \m{1}, \pd{2} \right) = 0 \]
implies
\[ a \cdot \left( 1 - 1/t_0^2 \right) = c_1^- + c_2^+, \]
and
\begin{eqnarray*}
0 &=& \left( \m{1}, \ld{E} \right) + \left( \md{1}, \l{E} \right) + \left( \p{2}, \ld{E} \right) + \left( \pd{2}, \l{E} \right) \\
&=& - t_0 \, e - \sqrt{2} c_1^- + \sqrt{2} \, a - e / t_0 - \sqrt{2} \, c_2^+ - \sqrt{2} \, c/t_0^2 \\
&=& -\left( 1+1/t_0 \right) \, e - \sqrt{2} \left( c_1^- + c_2^+ \right) + \sqrt{2} \left( 1 - 1/t_0^2 \right) \, a \\
&=& - \left( 1 + 1/t_0 \right) e.
\end{eqnarray*}
As $t_0 >0$, this implies $e=0$, and thus $\ld{E}=0$.

Now that we know $\ld{E}=0$, we can argue as in Lemma \ref{computation lem 1} to conclude that
\[ c_1^+ = -a / t_0^2, \quad c_1^- = a, \quad c_2^+ = - a/t_0^2, \quad c_2^- = a. \]
This determines $\pd{0}$, $\md{0}$, $\ldots$, $\md{3}$ up to a single scalar $a$.

This is half of our goal.  It remains to show $\ld{F}=0$.  Notice that all we did was begin with the four numbered generators $\p{0}$, $\m{0}$, $\p{3}$, $\m{3}$ of the front-right portion of the cube of figure \ref{cube figure}, move left, and apply known relations until $\ld{E}=0$ popped out the calculations.  The next step is to proceed right around the cube to the quartet $\p{4}$, $\m{4}$, $\p{5}$, and $\m{5}$.  Argue exactly as we did for the quartet $\p{1}$, $\m{1}$, $\p{2}$, and $\m{2}$ to show that $\ld{F}=0$.

Once we know that $\ld{F}=0$ one can compute
\begin{eqnarray*}
\pd{4} = -a/t_0^2 \cdot \left( \sqrt{2}, -1,1,-1,t_0 \right) & \quad &
    \md{4} = a \cdot \left( \sqrt{2}, -1, 1, -1, -1/t_0 \right) \\
\pd{5} = -a/t_0^2 \cdot \left( \sqrt{2}, -1, 1, 1, -t_0 \right) & \quad &
    \md{5} = a \cdot \left( \sqrt{2}, -1, 1,1,1/t_0 \right)
\end{eqnarray*}
Finally, repeating this procedure again for the quartet $\p{6}$, $\m{6}$, $\p{7}$, and $\m{7}$ yields
\begin{eqnarray*}
\pd{6} = -a/t_0^2 \cdot \left( \sqrt{2}, -1,-1,1,t_0 \right) & \quad &
    \md{6} = a \cdot \left( \sqrt{2}, -1,-1,1,-1/t_0 \right) \\
\pd{7} = -a/t_0^2 \cdot \left( \sqrt{2}, -1,-1,-1,-t_0 \right) & \quad &
    \md{7} = a \cdot \left( \sqrt{2}, -1,-1,-1,1/t_0 \right)
\end{eqnarray*}

Because we have not imposed all of the required orthogonality conditions, these calculations demonstrate only that the Zariski tangent space at $\rho_{t_0}$ is at most $1$-dimensional,
corresponding to the choice of the constant $a$ above.  As we will see in the proof of 
Theorem \ref{infinitesimal thm} below, because of the existence of our path of representations, 
this is sufficient to show that the Zariski tangent space is exactly $1$-dimensional and that 
the variety is smooth.  But first we note that, for a suitable choice of $a$, the infinitesimal 
variation in the space-like vectors given above equals that obtained by simply taking the $t$-derivative of the deformation given in table \ref{qt table} after rescaling.  To see this 
we must scale the vectors of table \ref{qt table} so that their norms remain constant.  
Namely, multiply the $+$ vectors by the function
\[ \sqrt{ \frac{1 + 1/t_0^2}{1 + 1/t^2} }, \]
and the $-$ vectors by
\[ \sqrt{ \frac{1 + t_0^2}{1 +t^2} }. \]
(Note that $t_0$ is constant.)  The resulting rescaled family will have $t$-derivatives at $t_0$ equal to $\pd{0}$, $\md{0}$, $\ldots$, $\md{7}$ above for $a=-t_0 / (1+t_0^2)$.  As always, this rescaling has no effect on the $1$-parameter family of representations it describes.

An application of the implicit function theorem now implies the following theorem: 

\begin{thm} \label{infinitesimal thm}
For any $t_0>0$ consider the representation $\rho_{t_0} \in \Hom (\Gamma_{22},G)$, where $G=\Isom (\Hf)$.  Let 
$\Def(\Gamma_{22},G)$ denote the slice of conjugacy classes in $\Hom(\Gamma_{22},G)$ determined 
by the conditions that the space-like vectors corresponding to the $A$, $B$, and $C$ walls are 
constant and the vector corresponding to the $D$ wall is orthogonal to the vector $e_4$.  Then the point $\rho_{t_0}$ has a neighborhood in $\Def(\Gamma_{22},G)$ that is a smooth $1$-dimensional manifold.  The Zariski tangent space of $\Def (\Gamma_{22},G)$ at $\rho_{t_0}$ is equal to the tangent space.    Near $\rho_{t_0}$ this deformation space is parametrized by the $1$-parameter family of representations $\rho_t$ described in table \ref{qt table}.
\end{thm}

\begin{pf}
A point in $\Def(\Gamma_{22},G)$ near $\rho_{t_0}$ is determined by space-like 
vectors $\{q_i\}$ in $\mathbb{R}^{1,4}$  associated with reflections in the $22$ walls of $P_{22}$.
The coordinates of all these vectors determine a point in $\mathbb{R}^k$, for some large value $k$.
(If we suppress the coordinates
of the $A$, $B$, and $C$ vectors and the last coordinate of the $D$ vector, all of which are
assumed fixed in our slice, then $k = 94$ coordinates suffice to uniquely determine the vectors, but this is not important in our argument.)   We denote by $\rho(t_0)$ the point in $\mathbb{R}^k$
corresponding to the representation $\rho_{t_0}$.
Vectors $\{q_i\}$ correspond to a point in 
$\Def(\Gamma_{22},G)$ if and only if certain pairs are orthogonal.  This means that their 
dot products are $0$, a quadratic condition.  Since a reflection depends only on a space-like vector up to scale we assume that the norms of the vectors are the same as those corresponding to $\rho(t_0)$, again a quadratic condition arising from dot products.

Define a map $r: \mathbb{R}^k \to \mathbb{R}^{\ell}$ by taking a point in $\mathbb{R}^k$ to this 
collection of dot products of the associated vectors.  Then, near $\rho(t_0)$, $\Def(\Gamma_{22},G)  = r^{-1}(p)$, where $p = r(\rho(t_0))$.  We have thus identified $\Def(\Gamma_{22},G)$ near 
$\rho_{t_0}$ with a subset of $\mathbb{R}^k$ containing the point $\rho(t_0)$, in fact as the preimage of a single point of a polynomial map. 

The computations of the last two sections imply that
the kernel of the derivative $dr$ at $\rho(t_0)$ is at most $1$-dimensional. A priori we can 
only conclude that there is an inequality because we did not impose all of the orthogonality conditions.  On the other hand, consider the map 
$\rho: \mathbb{R}_{+} \to \mathbb{R}^k$ sending $t \in \mathbb{R}_{+}$ to the coordinates
of the space-like vectors in table \ref{qt table} at time $t$, but rescaled as above.  
The derivative of this map is injective. Since all of $\rho(\mathbb{R}_{+}) \subset \Def(\Gamma_{22},G)$ maps to 
$p= r(\rho(t_0))$ under $r$, the kernel of $dr$ is
at least $1$-dimensional, hence exactly $1$-dimensional.

Now consider the sequence of maps
\[ \mathbb{R}_{+} \rightarrow \mathbb{R}^k \rightarrow \mathbb{R}^{\ell} \]
given by the composition $r \circ \rho$.  This composition is constant with image equal to
$p$.  Furthermore, we have just shown that the corresponding sequence 
$dr \circ d\rho$ of maps of tangent spaces
is exact at the central term.

These are the hypotheses of Lemma 6.8 of \cite{Rag}, which is an application
of the implicit function theorem (and is used there to prove Weil's lemma).
The conclusion of this Lemma is that, near the point where the derivative is computed, the
preimage $r^{-1}(p)$ is a smooth manifold whose tangent space equals the kernel of $dr$, 
which equals the image of $d \rho$.   Theorem \ref{infinitesimal thm} now follows.

\end{pf}

Near $\rho_{t_0}$, the orbit of the action of $G$ on 
$\Hom(\Gamma_{22},G)$ by conjugation is a smooth $10$-dimensional manifold, the same dimension as 
$G$.  Since the slice described above is transverse to this orbit,  
Theorem \ref{infinitesimal thm} could equally well be phrased as saying that $\Hom(\Gamma_{22},G)$
is a smooth $11$-dimensional manifold near $\rho_{t_0}$.

Theorem \ref{infinitesimal thm} can be used to study the quotient groups $\Lambda_n = \rho_{t_n}(\Gamma_{22})$ of $\Gamma_{22}$, where $n \in \{3,4, \ldots\}$ and $t_n$ is defined by equation \ref{discrete values}.  Recall that each $\Lambda_n$ is a discrete subgroup of $G$.  Our analysis of the infinitesimal deformations of $\Gamma_{22}$ implies:

\begin{cor} \label{infinitesimal rigidity cor}
The inclusion map of $\Lambda_n$ into $G$ is an infinitesimally rigid representation of $\Hom(\Lambda_n,G)$.  Thus, at this representation, the variety $\Hom(\Lambda_n,G)$ is 
$10$-dimensional and smooth, corresponding locally to the $G$-action by conjugation.
\end{cor}

\begin{pf}
The group $\Lambda_n$ contains additional torsion relations not found in the presentation of $\Gamma_{22}$.  For example, walls $\p{1}$ and $\p{3}$ intersect at an angle $2 \pi/n$ in the arrangement corresponding to $\Lambda_n$.  Therefore the group element obtained by concatenating these two reflections will be a rotation around their plane of intersection by angle 
$4 \pi/n$ and hence will have finite order.  Any infinitesimal deformation of 
$\Lambda_n$ must infinitesimally preserve this relation and thus the angle between these walls.  Since any such infinitesimal deformation is also an infinitesimal deformation of $\Gamma_{22}$,  it corresponds to varying the $t$ parameter of $\rho_{t}$ infinitesimally.  But, as computed in Proposition \ref{angle prop}, this changes the angle of intersection to the first order. Thus, there are no nontrivial infinitesimal deformations of $\Lambda_n$.

The final statement is Weil's lemma \cite{We,Rag} applied to the $10$-dimensional Lie group $G$, using the fact that these representations have trivial centralizer in $G$.
\end{pf}

Corollary \ref{infinitesimal rigidity cor} and Proposition \ref{boundary groups prop} combine to give explicit examples of $4$-dimensional geometrically finite groups that are infinitesimally rigid even though their boundary subgroups have nontrivial deformations.  This demonstrates the failure in dimension four of the $3$-dimensional Ahlfors-Bers deformation theory of geometrically finite Kleinian groups \cite{And}.  That theory implies that the deformations of a
$3$-dimensional geometrically finite group of hyperbolic isometries are locally parametrized by the
conformal structures of the boundary surfaces at infinity.  In particular, if the boundary
subgroups have nontrivial deformations, so will the groups themselves.

Combining this with Proposition \ref{boundary never geodesic} we can conclude the following:

\begin{cor} \label{infinitesimal rigidity cor 2}
For $n = 2m$, where $m \ge 4$ is an integer, the geometrically finite, infinite covolume group $\Lambda_n < G$ is infinitesimally rigid, Zariski dense, and its convex hull does not have totally geodesic boundary.  The boundary subgroups of $\Lambda_n$ all have nontrivial deformations in
$G$.
\end{cor}

\section{Lattices in the deformation} \label{manual lattice}

In Section \ref{defining the deformation 2} we showed that  the deformed reflection group 
$\U{\Lambda}_n := \rho_{t_n} (\U{\Gamma}_{22}) < \text{Isom}(\Hf)$, is discrete for $n$ an integer at least $3$.  Figure \ref{coxeter_diagram2} shows the Coxeter diagram for this group, where $\nu=n$.  Using this diagram one can prove the following proposition.

\begin{prop} \label{penny prop}
The discrete group $\U{\Lambda}_n := \rho_{t_n} (\U{\Gamma}_{22}) < \text{Isom}(\Hf)$ is a nonuniform lattice when $n \in \{3,4,5\}$.
\end{prop}

\noindent This proposition immediately implies the finite index subgroup $\Lambda_n = \rho_{t_n}(\Gamma_{22})$ is also a lattice for $n \in \{3,4,5 \}$.  A floating point computation suggesting this was explained in Section \ref{disappear}.  Note that Sections \ref{viewed in the sphere at infinity} and \ref{fundamental domain} built fundamental domains only for $n>6$; other methods are necessary to analyse the remaining cases.  We now explain how to prove Proposition \ref{penny prop}.

The proof is a lengthy but manageable application of a criterion, due to Vinberg \cite{Vin1},  for determining whether or not a hyperbolic reflection group is a lattice.  This criterion is given in terms of the group's Coxeter diagram, using the diagram to determine the combinatorics of the polytope which forms a fundamental region for the Coxeter group.  For simplicity we focus on the case of $\Hf$.  The idea of Vinberg's criterion is to look for the polytope's edges in the diagram, and check whether or not each end of each edge terminates in either a vertex inside $\Hf$ or a cusp of rank $3$.  This is done by examining subdiagrams of the Coxeter diagram.  At this point the reader might find it useful to obtain a table of the connected elliptic and parabolic Coxeter diagrams.  For example, see \cite[pgs.202-203]{Vin2}.

The link of  a finite vertex determines a $3$-dimensional spherical Coxeter group; this corresponds to an elliptic  subdiagram of the original Coxeter diagram.  Every $n$-dimensional spherical Coxeter group is obtained by reflection in some $n$-dimensional spherical simplex, so for a $4$-dimensional hyperbolic Coxeter group,   the subdiagrams corresponding to a finite vertex will have exactly $4$ vertices.  The disjoint union of an elliptic diagram and a vertex is again elliptic, corresponding to the spherical cone on a lower dimensional spherical simplex.  All elliptic diagrams can be built up from connected ones in this way.  Similarly,  by intersecting a neighborhood of an ideal vertex with a horosphere based there, one obtains a Euclidean
Coxeter group; this corresponds to a parabolic subdiagram.  Any $n$-dimensional Euclidean Coxeter
group determined by a \emph{connected} parabolic diagram arises as a reflection group in some 
$n$-dimensional Euclidean simplex.  The disjoint union of parabolic diagrams is again parabolic, corresponding to
the direct product of the associated reflection groups.  All parabolic diagrams can be built
from connected ones in this way, and one can easily compute the rank from the ranks of the connected diagrams.

In dimension $4$ an edge is seen in the diagram as a triple of 
vertices forming an elliptic subdiagram.  The corresponding Coxeter subgroup is the isotropy group of the edge in the associated polytope.  It is a $2$-dimensional spherical reflection group whose fundamental polygon is determined
by the intersection of  a hyperplane normal to the edge with a small neighborhood of the edge.

Vinberg's criterion in dimension $4$ is then given algorithmically by first finding all the triples
of vertices in the original Coxeter diagram that determine elliptic subdiagrams and then checking
that each one extends in exactly two ways: either to an elliptic subdiagram with $4$ vertices
or to a parabolic subdiagram of rank $3$.

We outline this process in our case; some details are left to the reader.  We refer constantly to the Coxeter diagram in figure \ref{coxeter_diagram2} where we will assume for the moment that the
value $\nu > 2$ equals an integer $n$.   The fundamental polytope for the group $\U{\Lambda}_n$ will be denoted by $\U{\mathcal{F}}_{t_n}$.

Since there are $8$ vertices
there are $56$ possible triples.  An obvious necessary condition for a triple to determine an edge is
that  each pair of the triple should correspond to walls which intersect. This rules out $28$ possibilities immediately.  
For example,
the triple $\{ \l{A}, \l{L}, \l{M} \}$ does not form an edge because $W_{\l{A}}$ and $W_{\l{L}}$ are tangent.

\begin{figure}[ht]
\labellist
\small\hair 2pt
\pinlabel $\p{0}$ [b] at 25 650
\pinlabel $\p{3}$ [t] at 25 4
\pinlabel $\l{L}$ [r] at 435 327
\pinlabel $\l{A}$ [b] at 685 360
\pinlabel $\m{0}$ [b] at 880 650
\pinlabel $\m{3}$ [t] at 880 4
\endlabellist
\centering
\includegraphics[scale=0.1]{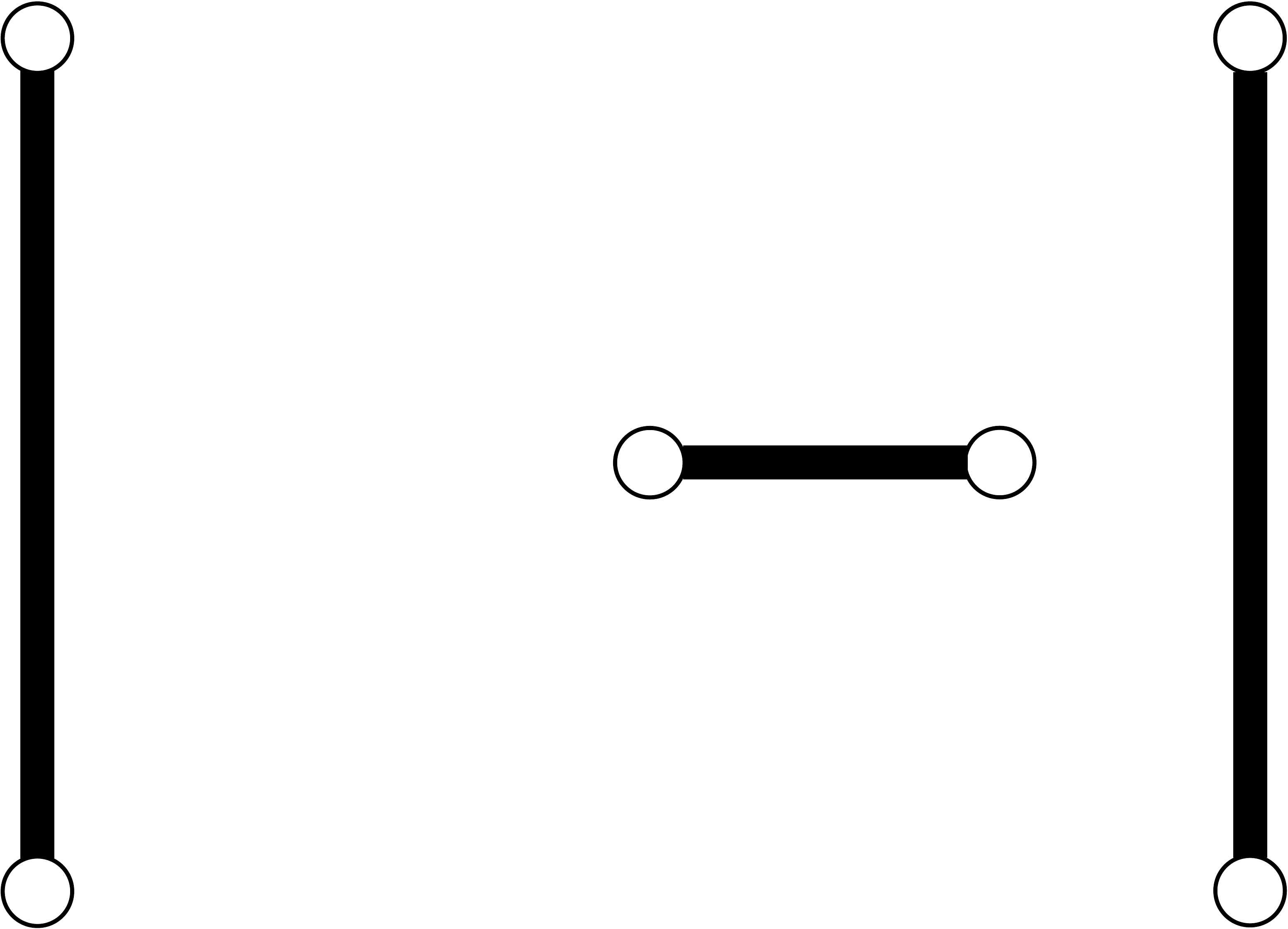}
\caption{The Coxeter diagram of a cubical reflection group}
\label{coxeter_diagram3}
\end{figure}
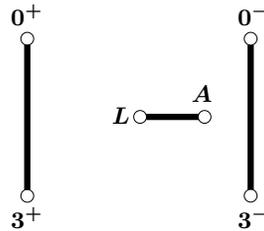

We note that there is exactly one parabolic subdiagram of rank $3$ in this Coxeter diagram.
It is shown in figure \ref{coxeter_diagram3};  it  corresponds to the group generated by reflections in the faces of a cube, which is a horocyclic cross-section of the unique cusp of $\U{\mathcal{F}}_{t_n}$.
There are
$8 = 2^3$ edges coming from the different ways of choosing one vertex from each of the connected components in this parabolic subdiagram.  The fact that these triples enlarge to a parabolic subdiagram indicates that one end of each edge terminates at infinity at the fixed point of a rank $3$ cusp group.  The $8$ edges correspond to the $8$ vertices of the cube  cross-section.  As there is only one cusp, we must show that the other end of each of these edges terminates at a finite vertex.


There are $20$ remaining triples that could possibly correspond to edges and, when $n < 6$, they do, indeed, define elliptic subdiagrams.   Rather than show this directly, we make the following argument.
Any edge that corresponds to one of these remaining triples must be a finite edge because
we have used up all of the ends of edges going out the cubical cusp.  Thus we should have $40$ finite ends of edges from these triples.  Together with the $8$ finite ends of the infinite edges,
there should be a total of  $48$ finite ends of edges.  Any finite vertex has valence $4$ since its
spherical link is a combinatorial tetrahedron.  Thus it suffices to describe $12$ finite vertices.

We will describe $6$ vertices; the other $6$ are the images of these under the roll symmetry
(which is the reflection in horizontal line through $\l{L}$ and $\l{A}$ in figure \ref{coxeter_diagram2}).
The first two vertices correspond to the quartets $\{ \p{3}, \m{3}, \l{N}, \l{A} \}$ and 
$\{ \p{3}, \m{3}, \l{N}, \l{L } \}$  whose polyhedra are the spherical cone on the $(2,2,2)$ and 
$(2,2,3)$ spherical triangles, respectively.  These are clearly spherical and independent of $n$.  The second two
vertices correspond to $\{ \p{3}, \l{M}, \l{N}, \l{A} \}$ and $\{ \p{3}, \l{M}, \m{0}, \l{A} \}$.  Their associated
polyhedra are the same, equal to the cone on the $(2,2,n)$ spherical triangle.  Although the
geometry depends on $n$ it is always spherical.  Finally, the last two quartets are only
elliptic when $n < 6$.  They are $\{ \p{3}, \l{M}, \l{L}, \m{0} \}$ and $\{ \p{3}, \l{M}, \l{L}, \l{N} \}$.
The first is the cone on the $(2,3,n)$ triangle which is spherical when $n < 6$ but not otherwise.
The second determines a tetrahedron which is spherical for integer values $n \leq 5$.

This completes  the proof of Proposition \ref{penny prop}.

\subsection{Transitional combinatorics}

In this section we will use the analysis from the previous section to understand the transitional 
behavior of the combinatorics of the polytopes $\U{\mathcal{F}}_{t}$.  Rather than considering only values of $t$ when the angle 
between $W_{\p{3}}$ and $W_{\l{M}}$ and between $W_{\p{0}}$ and $W_{\l{N}}$ is of
the form $\pi/n$, we allow $t$ and hence the angle to vary continuously.  This amounts to letting the
value $\nu$  in the Coxeter diagram in figure \ref{coxeter_diagram2} vary continuously.  We always
assume that $\nu >2$.

Our previous analysis amounted to considering the edges and vertices of the polytopes $\U{\mathcal{F}}_{t_n}$ via the Coxeter diagram.  This makes sense for any value of $t$  
and the analog of Vinberg's criterion is the same.  A finite vertex in the polytope corresponds 
to a quartet of vertices in the Coxeter diagram that determines a spherical tetrahedron and an
edge in the polytope corresponds to a triple of vertices in the Coxeter diagram that determines
a spherical triangle.  The polytope has finite volume if and only if every edge has two vertices
each of which is either finite or ends at a rank $3$ cusp.

There are two transitional values for $1 > t > 0$  (or for angles between $0$ and $\pi/2$), one at 
$t = t_6$, when the angle equals $\pi /3$, and the
other at $t = \sqrt{1/2}$, which corresponds to the angle 
$\frac {\text{cos}^{-1} (1/3)} {2}$.
In section \ref{disappear} we showed numerically that the combinatorics of the polytopes
$\mathcal{F}_{t}$ for the groups $\rho_t (\Gamma_{22})$ change at these values.  We will now be able to prove rigorously that, for the polytopes $\U{\mathcal{F}}_{t}$ of the representations of the  extended group,  these changes
occur and that these are the only combinatorial changes that occur.  This follows from the
fact that any combinatorial change will create or destroy vertices or edges and we can
now analyze all such occurrences.  We will then check that these changes in the polytopes 
$\U{\mathcal{F}}_{t}$ correspond to the numerically observed changes in the polytopes
$\mathcal{F}_{t}$.

In the previous section we showed that all $28$ triples of vertices in the Coxeter diagram that can correspond to edges for some value of $\nu$ do actually determine an edge when 
$\nu = n \leq 5$.  For these values, there are
$12$ quartets that correspond to finite vertices.   These statements are not true for values
$\nu \geq 6$.

\begin{figure}[ht]
\labellist
\small\hair 2pt
\pinlabel $\l{M}$ [l] at 54 27
\pinlabel $\l{L}$ [l] at 54 328
\pinlabel $\l{N}$ [l] at 54 623
\endlabellist
\centering
\includegraphics[scale=0.1]{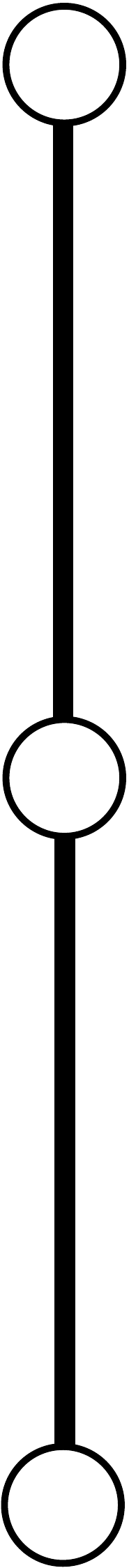}
\caption{A problematic edge}
\label{coxeter_diagram4}
\end{figure}

\begin{figure}[ht]
\labellist
\small\hair 2pt
\pinlabel $\p{0}$ [t] at 29 600 
\pinlabel $\l{N}$ [l] at 473 623
\pinlabel $\l{L}$ [l] at 473 328
\pinlabel $\l{M}$ [r] at 421 27
\pinlabel $\l{N}$ [l] at 1190 622
\pinlabel $\l{L}$ [l] at 1190 327
\pinlabel $\l{M}$ [l] at 1190 28
\pinlabel $\p{3}$ [b] at 740 54
\pinlabel $\nu$ [b] at 231 630
\pinlabel $\nu$ [b] at 956 34
\endlabellist
\centering
\includegraphics[scale=0.1]{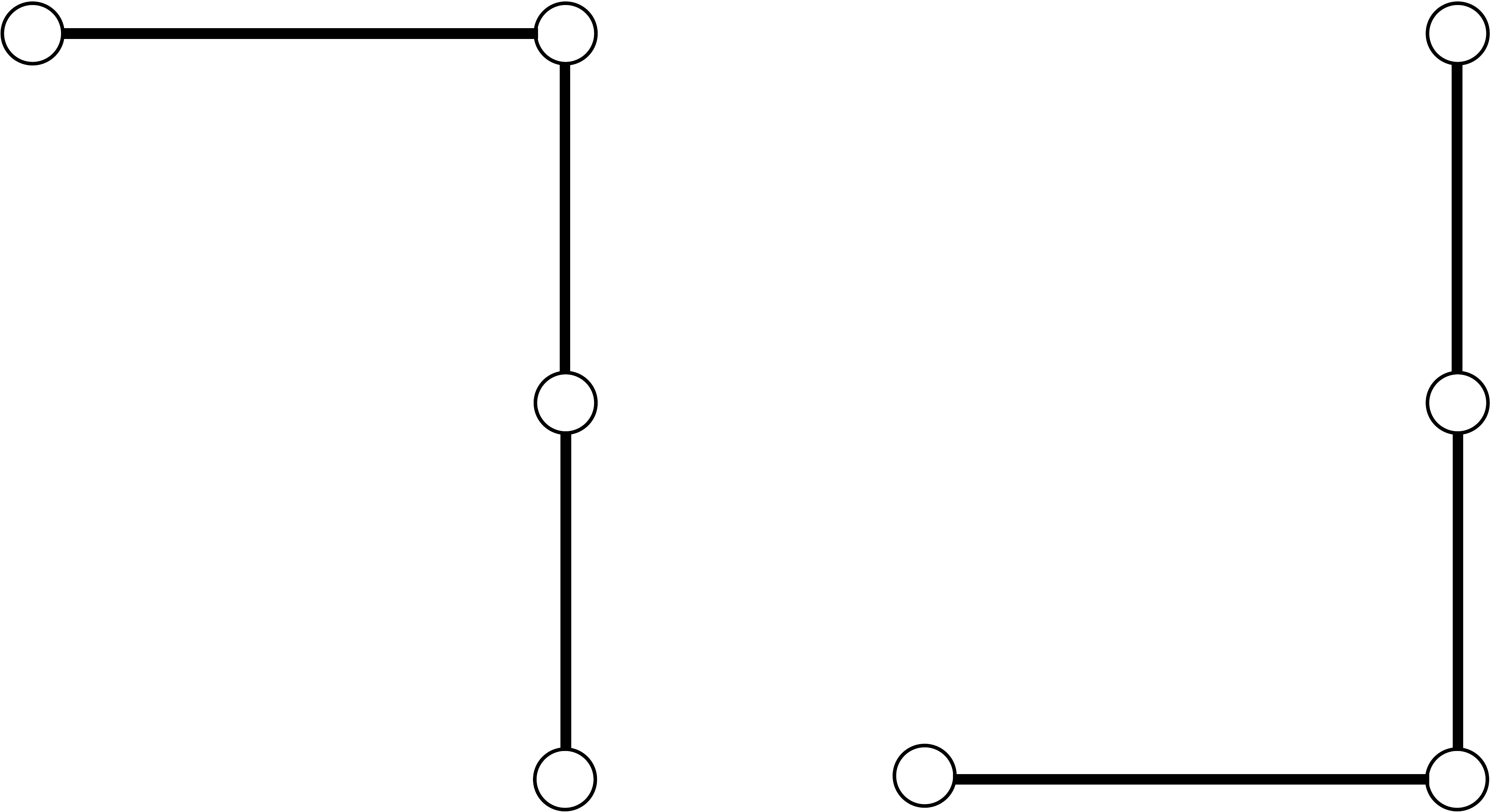}
\caption{Vertices of a problematic edge}
\label{coxeter_diagram5}
\end{figure}

To understand more precisely how things change for  $\nu >5$, consider the subdiagram of figure \ref{coxeter_diagram2} shown in figure \ref{coxeter_diagram4}.  This subdiagram corresponds to an edge of $\U{\mathcal{F}}_{t}$,  independent of $t$.  For $\nu \leq  5$, this edge has finite length and terminates in vertices corresponding to the subdiagrams in figure \ref{coxeter_diagram5}.

\begin{figure}[ht]
\labellist
\small\hair 2pt
\pinlabel $2$ [t] at 301 37
\pinlabel $3$ [l] at 596 226
\pinlabel $2$ [bl] at 446 390
\pinlabel $2$ [br] at 140 250
\pinlabel $3$ [b] at 334 197
\pinlabel $\nu$ [l] at 385 308
\endlabellist
\centering
\includegraphics[scale=0.2]{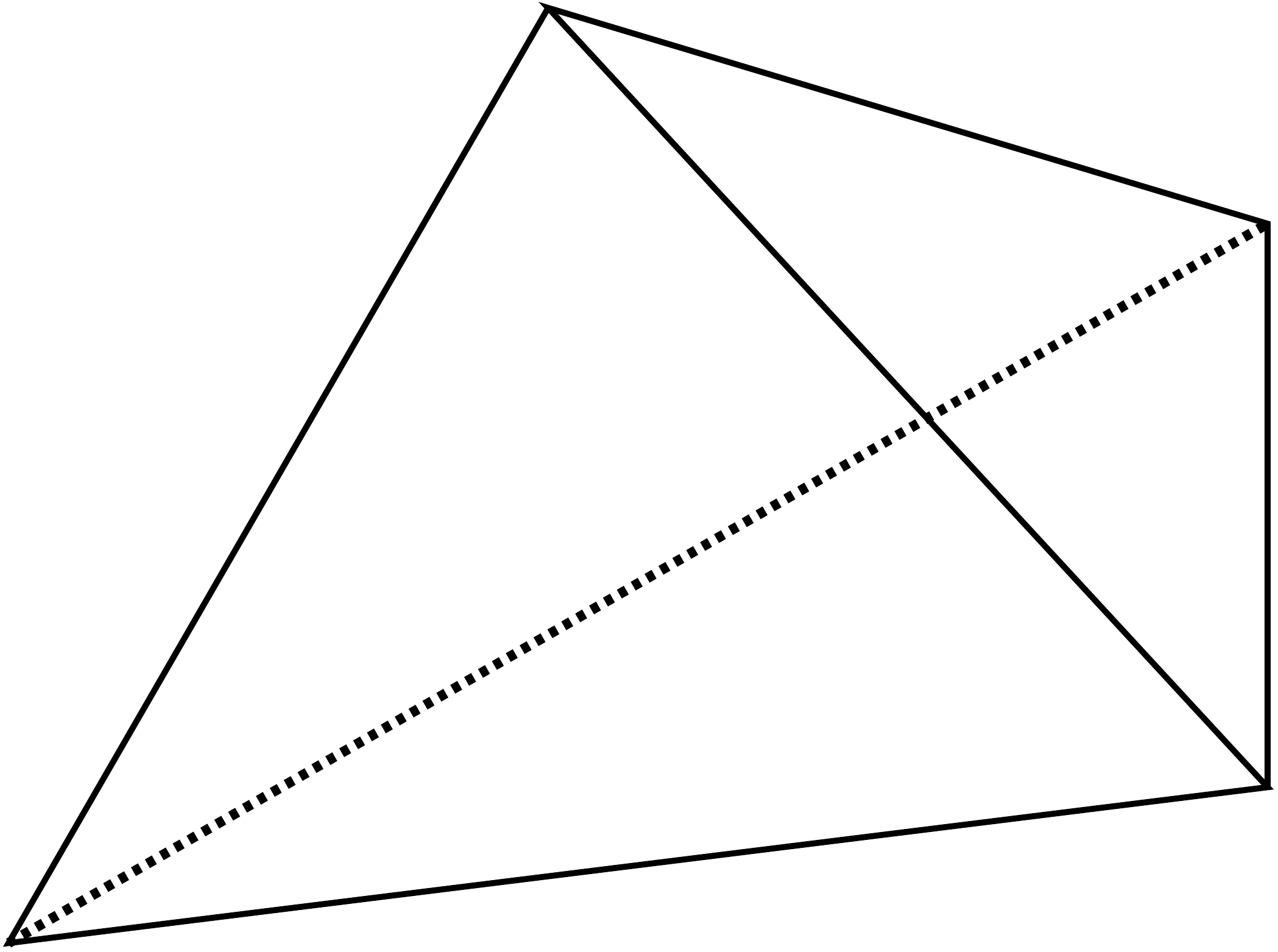}
\caption{The link of the disappearing vertices}
\label{tetrahedron1}
\end{figure}

When $n \leq 5$, the subdiagrams in figure \ref{coxeter_diagram5} correspond to an elliptic Coxeter group with polyhedron given by the spherical tetrahedron shown in figure \ref{tetrahedron1}, where  the dihedral angle at an edge equals $\pi$ divided by its label.   For any value of $\nu$, one obtains a tetrahedron with these angles by taking  a regular tetrahedron with dihedral angle  $\frac {2\pi}{\nu}$ 
and dividing out by its symmetries.   When the dihedral angle at the edge labeled $\nu$ equals 
$\frac {\text{cos}^{-1} (1/3)} {2}$ (at $\nu \approx 5.104$),  the tetrahedron in figure \ref{tetrahedron1} is no longer spherical.  It is Euclidean and the edge ends at two ideal vertices.  (The dihedral angles of a regular Euclidean
tetrahedron equal  $\text{cos}^{-1} (1/3)$.)   At this value,  the volume of the polytope is still finite.  However, for any larger value of $\nu$, the tetrahedron is not Euclidean or spherical (it is actually hyperbolic) and the corresponding vertex has disappeared from $\U{\mathcal{F}}_{t}$.  The edge corresponding to the subdiagram in figure \ref{coxeter_diagram4} has no finite or cusped vertices 
so the Vinberg criterion fails, and the polytope has infinite volume.

There is a second combinatorial change that occurs when the dihedral angle passes through $\pi/6$
(when $\nu$ passes through value $6$).
At this value, two vertices and two edges are lost.   The two edges that are lost are associated to the triple $ \{ \p{3}, \l{M}, \l{L}\}$ and its image $\{ \p{0}, \l{N}, \l{L}\}$ under the roll symmetry.   These
subdiagrams correspond to the group of reflections in the $(2,3, \nu)$  triangle ( with angles $\pi/2$,  $\pi/3$, and  $\frac{\pi} {\nu}$), which is spherical if and only if $\nu < 6$.  The two vertices
lost are associated to the quartet $ \{ \p{3}, \l{M}, \l{L}, \m{0}\}$ and its image $\{ \p{0}, \l{N}, \l{L},\m{3}\}$ under the roll symmetry.  At $\nu = 6$ they are parabolic subdiagrams but of rank $2$ and, 
for larger values
of $\nu$,  the corresponding four walls have no common intersection, even at infinity.

These are the only possible combinatorial changes.  To see this, note that the arguments in the previous section show that all $28$ triples of vertices in the Coxeter diagram that could correspond edges for some value of $\nu$ did occur for $\nu  \leq 5$.  Thus,  no new edges
can occur as $\nu$ is increased.  Also, for 
 $\nu \leq 5$,  every edge has at least one finite vertex, so all possible edge triples occur as a 
subset of at least one of the
$12$ quartets that correspond to finite vertices.  A quick check shows that,  except for 
$ \{ \p{3}, \l{M}, \l{L}\}$ and $\{ \p{0}, \l{N}, \l{L}\}$,  all such triples correspond to 
spherical triangles for all values of $\nu$.  Thus these are the only ones that can be lost and this occurs only at $\nu = 6$.  
Similarly, there can be no new finite vertices and vertices can only be lost when and in the manner
described above.

We now confirm that the combinatorial changes in $\U{\mathcal{F}}_{t}$ at these values correspond to the changes in the polytopes $\mathcal{F}_{t}$ for the unextended groups $\rho_t (\Gamma_{22})$ 
that were observed numerically in Section \ref{disappear}.
As in that section we consider decreasing the value of $t$ from slightly less than $1$ where we
know that the description of the polytopes in Section \ref{fundamental domain} is correct.
For values of $t$ larger than $t_6$  (hence for angles smaller than $\pi/3$) the finite vertex corresponding to the quartet $ \{ \p{3}, \l{M}, \l{L}, \m{0}\}$ does not exist while it does exist for values of $t$ less than $t_6$.  At $t = t_6$ the quartet corresponds to an ideal vertex where the associated
$4$ walls intersect at infinity.  

To understand how the polytope $\mathcal{F}_{t}$ changes at
$t_6$ we note that it is the union of the $24$ copies of the polytope $\U{\mathcal{F}}_{t}$ which
are its orbit under the group generated reflections in the $\l{L}$, $\l{M}$ and $\l{N}$ walls.  A point of intersection of the walls in the quartet $ \{ \p{3}, \l{M}, \l{L}, \m{0}\}$ provides a point of intersection in $\mathcal{F}_{t}$ between the walls in the orbit of $ \{ \p{3}, \m{0}\}$ under the subgroup
generated by reflection in the $\l{L}$ and $\l{M}$ walls.  This subgroup preserves the wall $W_{\m{0}}$ so  this equals  $W_{\m{0}}$ together with the orbit of $W_{\p{3}}$.
This orbit equals the $ \{ \p{3}, \p{1}, \p{5}\}$ walls.  The intersection of these $3$ walls with $W_{\m{0}}$,
first at infinity when $t = t_6$ then at a finite vertex for slightly smaller values of $t$ is apparent in
the third and fourth pictures of figure \ref{0-_wall_figure}.  This corresponds to the triangular
region at infinity of $W_{\m{0}}$ first collapsing to a point and then disappearing.  The same process occurs on the other even negative walls, which are the orbit of $W_{\m{0}}$ under the full group
of symmetries.  Applying the roll symmetry  one sees it  occur also on the odd negative walls.

The transition at $t = t_6$ can also be seen in the third and fourth pictures in figure \ref{3+_wall_figure} which shows the pattern of intersection for $W_{\p{3}}$ at infinity.  
Here one sees the hexagonal region at infinity for 
$W_{\p{3}}$ collapsing first to a triangle with internal angles $0$ and then changing to a triangle with small internal angles.   One of the vertices of the triangle with internal angles $0$ is the ideal vertex
corresponding to the subdiagram $ \{ \p{3}, \l{M}, \l{L}, \m{0}\}$ which is parabolic at $t_6$.  The
fact that it is not of rank $3$ can be seen from the fact that it is on the boundary of an open region
at infinity for $W_{\p{3}}$.  The other vertices of this triangle are the orbit of this one under the
symmetries preserving $W_{\p{3}}$.  

In figure \ref{3+_wall_figure} one can also see the new edges that
appear for values $t > t_6$.  In $\U{\mathcal{F}}_{t}$ the triple $ \{ \p{3}, \l{M}, \l{L}\}$ determines an edge for these values of $t$.  By considering the orbit of $\p{3}$ under the group generated by reflection in the $\l{L}$ and $\l{M}$ walls, this provides an edge which equals the intersection of the 
$ \{ \p{3}, \p{1}, \p{5}\}$ walls.  One endpoint of this edge is one of the vertices of the triangle
in the fourth picture and the triple intersection is apparent from the fact that the circles corresponding
to the $ \p{1}$ and $\p{5}$ walls overlap.   This edge has one finite vertex corresponding 
to the quartet $ \{ \p{3}, \p{1}, \p{5}, \m{0}\}$ which defines a spherical polyhedron for $t > t_6$
but the fact that it has no other finite or cusped vertex implies that the polytope has infinite
volume.

One can do a similar analysis at $t = \sqrt{1/2}$ which corresponds to $1/2$ the dihedral
angle of a regular Euclidean tetrahedron.  At this value the quartet $\{ \p{3}, \l{M}, \l{L}, \l{N} \}$
determines a Euclidean tetrahedron, hence a quadruple intersection at infinity, and for
larger values of $t$ it determines a spherical tetrahedron, hence a finite vertex.  In the
polytope $\mathcal{F}_{t}$ this provides similar quadruple intersections between the
$ \{ \p{3}, \p{1}, \p{5}, \p{7}\}$ walls, where these represent the orbit of $\p{3}$ under the full
symmetry group.  The intersection at infinity case can be seen in the fifth picture of figure \ref{3+_wall_figure}
where the triangular region at infinity for $W_{\p{3}}$ has collapsed to a point and the finite
vertex case occurs in the sixth picture when this region has disappeared completely.

\subsection{A bonus lattice}

\begin{figure}[ht]
\labellist
\small\hair 2pt
\pinlabel $\p{0}$ [b] at 25 648
\pinlabel $\l{N}$ [b] at 451 648
\pinlabel $\m{0}$ [b] at 868 648
\pinlabel $\m{3}$ [t] at 25 2
\pinlabel $\l{M}$ [t] at 451 2
\pinlabel $\p{3}$ [t] at 868 2
\pinlabel $6$ [b] at 226 648
\pinlabel $6$ [t] at 226 2
\pinlabel $\l{L}$ [r] at 423 328
\pinlabel $\l{A}$ [b] at 590 347
\pinlabel $\l{G}$ [b] at 693 493
\pinlabel $\l{H}$ [t] at 693 161
\endlabellist
\centering
\includegraphics[scale=0.15]{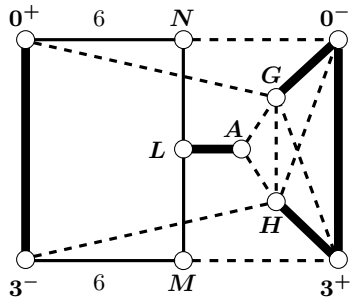}
\caption{The Coxeter diagram for $\Lambda_6$}
\label{coxeter_diagram6}
\end{figure}

Finally, we can use the polytopes in our family to construct one more lattice.  We have seen that
the group $\U{\Lambda}_6$ is  not a lattice.  However, we can extend it to a lattice by ``putting back in'' the reflections corresponding to $\l{H}$ and $\l{G}$.  Precisely, we redefine the letters $\l{H}$ and $\l{G}$ as
\[ \l{G}= \left( 1,0,0,0,\sqrt{\frac{6}{5}} \right) \quad \l{H} = \left(1,0,0,0,- \sqrt{\frac{6}{5}} \right), \]
and consider the extension defined as
\[ L_6 := \langle \U{\Lambda}_6 , \l{G}, \l{H} \rangle.\]
The Coxeter diagram for $L_6$ is the slightly messy figure \ref{coxeter_diagram6}.

\begin{figure}[ht]
\labellist
\small\hair 2pt
\pinlabel $6$ [b] at 243 33
\endlabellist
\centering
\includegraphics[scale=0.1]{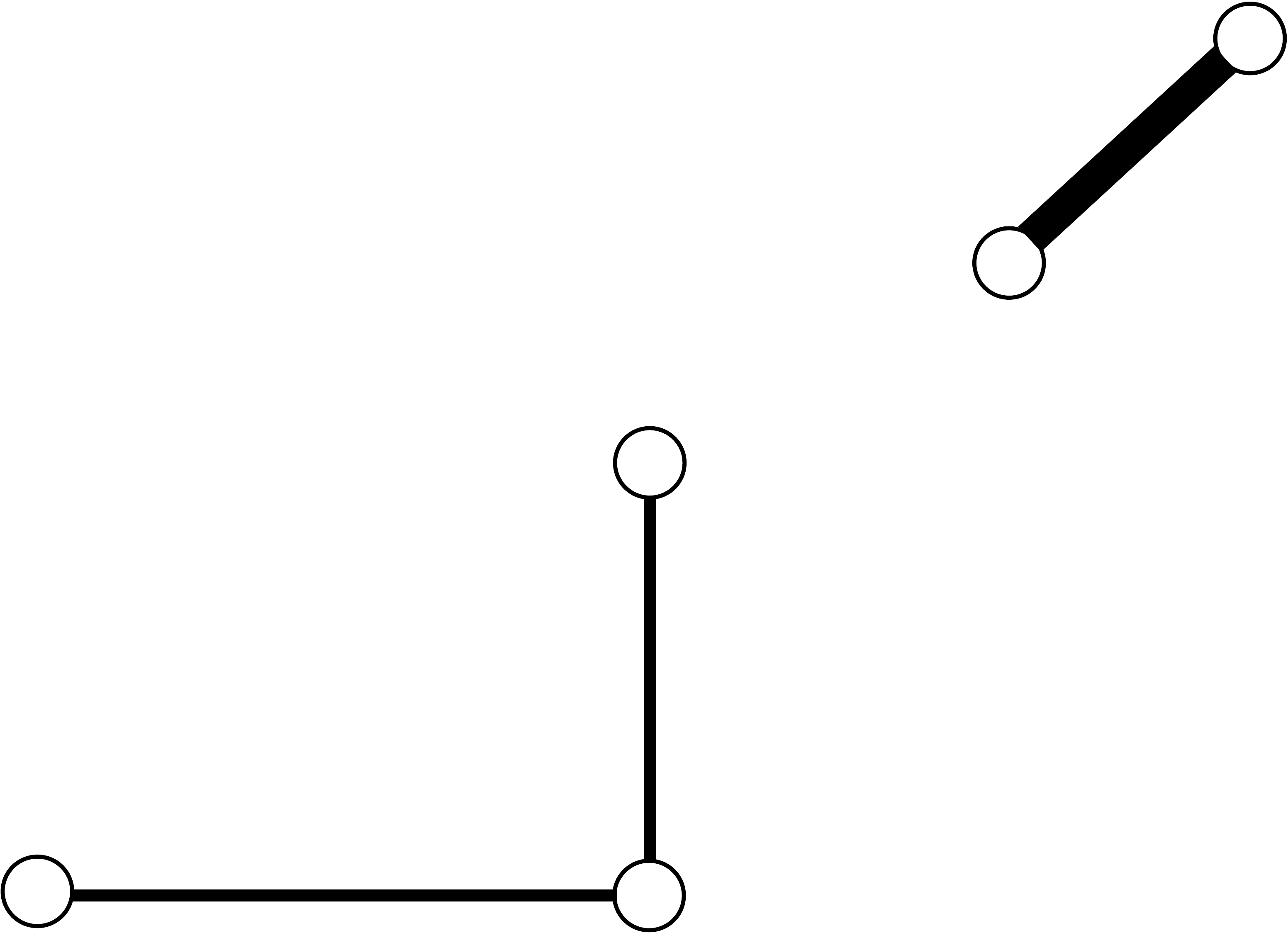}
\caption{A cusp group of $L_6$}
\label{coxeter_diagram7}
\end{figure}

By running through the procedure desribed above one can verify that the Coxeter diagram of figure \ref{coxeter_diagram6} corresponds to a finite volume polytope, and $L_6$ is therefore a lattice in $\text{Isom}(\Hf)$.  The only complication in this verification is the appearance of a new type of cusp with Coxeter diagram shown in figure \ref{coxeter_diagram7}.  (The figure is arranged to correspond with figure \ref{coxeter_diagram6}.)  This cusp has cross-section equal to the product of an interval with
the Euclidean $(2,3,6)$ triangle so it is Euclidean with rank $3$.

Let $\U{\mathcal{F}}_{t_6}$ be the Coxeter polytope of $\U{\Lambda}_6$.  The nonempty intersections with the new walls corresponding to $\l{G}$ and $\l{H}$ in $\U{\mathcal{F}}_{t_6}$ are all orthogonal.  This means $\l{G}$ and $\l{H}$ each cut off an infinite volume end of $\U{\mathcal{F}}_{t_6}$ with orthogonal intersections.  From this we can conclude that the Coxeter polytope of $L_6$ sits inside $\U{\mathcal{F}}_{t_6}$ as a convex core, and the convex hull of $\U{\Lambda}_6$ therefore has totally geodesic boundary.  The analogous fact for the polytope $\mathcal{F}_{t_6}$ of the non-extended group was nonrigorously observed in Section \ref{disappear}.

\subsection{Arithmeticity}

We next turn our attention to the question of when the group $\U{\Lambda}_n$ is arithmetic.  We will have no need of the general definition of arithmeticity for lattices in $\text{Isom}(\Hf)$.  For this we refer the reader to \cite{Borel1,Morris}.  In our setting we need only a criterion due to Vinberg for arithmeticity applicable to reflection groups in $\text{Isom}(\Hf)$ with noncompact finite volume Coxeter polytopes \cite[Ch.5,Thm.3.1]{Vin2}.

To describe this criterion, let $R < \text{Isom}(\mathbb{H}^n)$ be a reflection group with a noncompact finite volume Coxeter polytope.  (There is another criterion in the compact case, but it is more complicated.)  Given the Coxeter diagram of $R$, we relabel its edges as follows:  For any pair of walls intersecting in an angle $\pi/n$ for $n>2$, label the corresponding edge with $2 \cos (\pi / n)$.  For any pair of walls tangent at infinity label the corresponding thickened edge with a $2$.  Finally, for any pair of walls separated by a distance $d>0$, label the corresonding edge with $2 \cosh (d)$.  In all three cases, the label is simply $-2$ times the Minkowski pairing of the respective unit space-like vectors.  Notice that 
for orthogonal pairs (which have no edge between them) this prescription would result in a label of 
$0$.  Since Vinberg's criterion involves taking products of such Minkowski inner products there is no need for a label between such pairs.
 
Now make a list of the closed cycles of the Coxeter diagram that have no repeated vertices except for the beginning and ending vertex.  (Going back and forth along a single edge is allowed.)  For each cycle of the list, form the product of the labels of the edges traversed.  This forms a list of numbers, one for each cycle.  Vinberg's criterion states that $R$ is an arithmetic lattice if and only if each number of this list is an integer.  (The general condition involves such products being algebraic integers in an algebraic
number field associated to $R$.  When the polytope is noncompact this field is always $\mathbb{Q}$.)

\begin{figure}[ht]
\labellist
\small\hair 2pt
\pinlabel $\p{0}$ [b] at 25 650
\pinlabel $\p{3}$ [t] at 25 4
\pinlabel $\l{M}$ [t] at 451 4
\pinlabel $\l{L}$ [r] at 435 327
\pinlabel $\l{A}$ [b] at 685 360
\pinlabel $\l{N}$ [b] at 451 650
\pinlabel $\m{0}$ [b] at 880 650
\pinlabel $\m{3}$ [t] at 880 4
\pinlabel $2\,t\,s$ [b] at 235 629
\pinlabel $2\,t\,s$ [t] at 235 21
\pinlabel $2\,s$ [b] at 669 629
\pinlabel $2\,s$ [t] at 669 21
\pinlabel $2$ [b] at 577 342
\pinlabel $1$ [r] at 443 482
\pinlabel $1$ [r] at 443 181
\pinlabel $2$ [l] at 881 330
\pinlabel $2$ [r] at 10 321
\endlabellist
\centering
\includegraphics[scale=0.15]{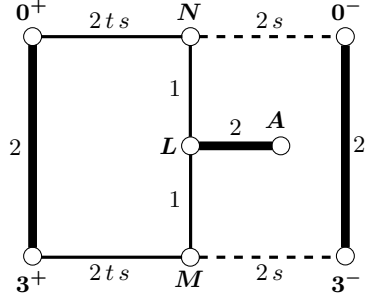}
\caption{The relabeled Coxeter diagram}
\label{coxeter_diagram8}
\end{figure}

Let us apply this criterion to $\U{\Lambda}_n$.  Figure \ref{coxeter_diagram8} shows its relabeled Coxeter graph, which is figure \ref{coxeter_diagram2} with edges marked as discussed above.  In the figure $t = t_n$ and 
\[ s = \sqrt{ \frac{2}{1 + t^2} }.\]  
A key observation is that any cycle of the Coxeter graph traversing a dotted edge must traverse a dotted edge twice.  Similarly, any cycle traversing an edge labeled $2\,t\,s$ must traverse such an edge twice.  To test Vinberg's criterion it therefore suffices to determine when
\[ \frac{8t^2}{1+t^2} \quad \text{and} \quad \frac{8}{1+t^2}\]
are integers.  This occurs if and only if $t \in \{ t_3 = \sqrt{\frac{1}{7}}, t_4 = \sqrt{\frac{1}{3}}, t_6 = \sqrt{\frac{3}{5}} \}$.  From this we may conclude that the lattices $\U{\Lambda}_3$ and $\U{\Lambda}_4$ are arithmetic, and the lattice $\U{\Lambda}_5$ is not arithmetic.  Moreover, we have shown that for $n>6$ no group $\U{\Lambda}_n$ 
is ever a subgroup of an arithmetic lattice.

Recall that arithmeticity is a property of a lattice's commensurability class \cite{Borel1,Morris}.  Using this we may restate the above results in terms of our original family $\Lambda_n$ as follows.

\begin{thm}\label{arithmeticity thm}
The lattices $\Lambda_3$ and $\Lambda_4$ are arithmetic.  The lattice $\Lambda_5$ is not arithmetic.
\end{thm}

Using these techniques, we leave it to the reader to prove that the extension $L_6$ of $\U{\Lambda}_6$ is also an arithmetic lattice.

\subsection{Retracts and Quotients}\label{retracts and quotients subsection}

Recall that a right-angled Coxeter group is a Coxeter group
\[ R = \langle r_1, r_2, \ldots, r_M \, | \, r_1^2 = \ldots = r_M^2 = (r_{i_1} r_{j_1})^{m_1} = \ldots 
    = (r_{i_N} r_{j_N} )^{m_N} = 1 \rangle \]
such that all the integers $m_i$ are equal to two.  This implies that all the generator pairs $r_{i_k}$, $r_{j_k}$ commute.  The lattice $\Gamma_{24}$ is an example of a right-angled Coxeter group.  These special groups have the following interesting property which we will use to construct nontrivial lattice quotients of $\Gamma_{24}$.

Pick a subset $S$ of the generators $\{r_1, \ldots , r_M \}$ of $R$.  Define the free Coxeter group
\[ C_M := \langle \rho_1, \ldots, \rho_M \, | \, \rho_1^2 = \ldots = \rho_M^2 = 1 \rangle, \]
and a homomorphism
\begin{eqnarray*}
f : C_M & \longrightarrow & \langle S \rangle < R \\
    \rho_i & \longmapsto & \begin{cases}
			    r_i & \text{if $r_i \in S$}\\
			    1  & \text{if $r_i \notin S$} \end{cases}
\end{eqnarray*}

We claim the homomorphism $f$ sends the relations of $R$ to the identity, and thus descends to a surjective homormorphism $\bar{f}:R \rightarrow \langle S \rangle$.  To see this, consider $r_i$ and $r_j$, a commuting pair of generators of $R$.  If $r_i$ and $r_j$ are in $S$ then
\[ f( \rho_i \rho_j \rho_i \rho_j ) = r_i r_j r_i r_j = 1. \]
If neither $r_i$ nor $r_j$ are in $S$ then $f(\rho_i \rho_j \rho_i \rho_j)$ is clearly the identity.  Finally, if $r_i \in S$ and $r_j \notin S$ then
\[ f(\rho_i \rho_j \rho_i \rho_j) = r_i \cdot 1 \cdot r_i \cdot 1 = 1. \]
This shows that $f$ descends to a surjection $\bar{f}: R \rightarrow \langle S \rangle$.  In the literature this homomorphism is called a retract of $R$ because the composition
\[ \langle S \rangle \stackrel{\text{incl}}{\longrightarrow} R \stackrel{\bar{f}}{\longrightarrow} 
    \langle S \rangle \]
is the identity map.

Let us apply this to $\Gamma_{22}$ and $\Gamma_{24}$.  By the above there exists a surjective homomorphism $\Gamma_{24} \rightarrow \Gamma_{22}$.  Recall that $\Lambda_n = \rho_n (\Gamma_{22})$ for $n \in \{3,4,5\}$ is a lattice in $\text{Isom}(\Hf)$.  The composition
\[ \Gamma_{24} \longrightarrow \Gamma_{22} \stackrel{\rho_n}{\longrightarrow} \Lambda_n \]
is a surjection between lattices of $\text{Isom}(\Hf)$ which is not an isomorphism.  In other words, the composition is a nontrivial quotient map.  We state this as a proposition.

\begin{prop}\label{quotient thm}
Inside $\text{Isom}(\Hf)$, the arithmetic lattices $\Lambda_3$ and $\Lambda_4$ and the nonarithmetic lattice $\Lambda_5$ are all nontrivial quotients of the arithmetic lattice $\Gamma_{24}$.  Moreover, by passing to finite index subgroups, an infinite number of similar lattice quotients can be constructed inside $\text{Isom}(\Hf)$.
\end{prop}

\section{To the cuboctahedron} \label{cubeoctahedron section}  

Finally, we consider the limiting representation of $\rho_t(\Gamma_{22})$ as $t \to 0$.
Referring to table \ref{qt table} we see that some of the entries of the positive numbered space-like vectors 
go to infinity.  However, the corresponding hyperplanes, and hence the corresponding reflections, depend
only on the space-like vectors up to scale.  If we rescale those vectors by multiplying by $t$, then they
all limit to $\pm e_4$; thus all the corresponding hyperplanes converge to a single hyperplane.  Meanwhile,
the $8$ negative numbered space-like vectors converge to the vectors $\sqrt{2} e_0 \pm e_1 \pm e_2 \pm e_3$.
All of these vectors are perpendicular to $e_4$ as are the letter space-like vectors (which are independent
of $t$).

Thus the limit representation, $\rho_0 (\Gamma_{22})$, preserves the $3$-dimensional hyperbolic
space $V \subset \Hf$ orthogonal to $e_4$.   On the sphere at infinity, the boundaries of hyperplanes
determined by the letter and negative numbered space-like vectors intersect $\partial V$ in the pattern
of circles of figure \ref{constant core intersection pattern}.  As we noted there, the entire collection of $2$-spheres at infinity corresponding
to all the letter and numbered walls intersect $\partial V$ in this same pattern for all values of
$t$.  However, the angle of intersection for the numbered walls varies with $t$ and we are now seeing
the pattern with the limit angles of $0$ and $\pi/2$ for the positive and negative numbered walls,
respectively.

\begin{figure}[ht]
\labellist
\small\hair 2pt
\pinlabel $\l{A}$ at 424 413
\pinlabel $\textcolor{blue}{\l{B}}$ at 655 107
\pinlabel $\textcolor{blue}{\l{C}}$ at 160 529
\pinlabel $\l{D}$ at 1020 533
\pinlabel $\l{E}$ at 585 948
\pinlabel $\textcolor{blue}{\l{F}}$ at 739 625
\pinlabel $\m{2}$ at 776 793
\pinlabel $\m{1}$ at 185 738
\pinlabel $\m{0}$ at 215 145
\pinlabel $\m{3}$ at 810 213
\pinlabel $\textcolor{blue}{\m{7}}$ at 951 892
\pinlabel $\textcolor{blue}{\m{6}}$ at 416 898
\pinlabel $\textcolor{blue}{\m{5}}$ at 396 266
\pinlabel $\textcolor{blue}{\m{4}}$ at 975 283
\pinlabel {Front faces are labeled in black.} [tl] at 825 48 
\pinlabel {Back faces are labeled in \textcolor{blue}{blue}.} [tl] at 825 0 
\endlabellist
\centering
\includegraphics[scale=0.2]{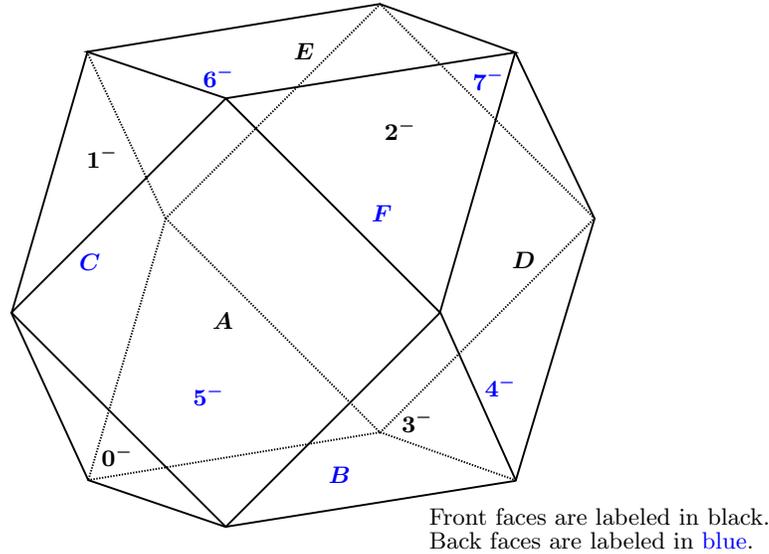}
\caption{The cuboctahedron}
\label{cubeoct2 fig}
\end{figure}

The group generated by reflections in the letter walls and the negative numbered walls preserves
$V$ and can thus be viewed as a discrete group of isometries of $\mathbb{H}^3$.  It equals the group of
reflections in a right-angled ideal \emph{cuboctahedron}.  Figure \ref{cubeoct2 fig} depicts a  combinatorial model of this $3$-dimensional
polyhedron, which we denote by $P_{co}$.  It derives its name from the fact that it
can obtained from either the cube or the octahedron by truncating their vertices by drawing lines
from the midpoints of the edges.  It can be realized as a polyhedron in $\mathbb{H}^3$ with all its
vertices at infinity and all its dihedral angles equal to $\pi/2$.  In our realization the rectangular faces come from the letter walls while the triangular faces come from the negative numbered walls.

Let $\Gamma_{co}$ denote the group generated by reflections in the faces of this hyperbolic
cuboctahedron.  It is a right-angled Coxeter group with $14$ order $2$ generators $r_i$ 
corresponding to the faces of the polyhedron and $12$ relations of the form $(r_j r_k)^2 = e$
corresponding to pairs of faces intersecting along an edge.  It is isomorphic to the subgroup
of $\Gamma_{22}$ generated by reflections in the letter and negative numbered walls.  Analogous 
to the discussion in the previous section, there is a retraction of 
$\Gamma_{22}$ onto $\Gamma_{co}$ obtained by
sending the generators corresponding to the positive walls to the identity.  Combined with
the retraction of $\Gamma_{24}$ onto $\Gamma_{22}$, we obtain a retraction of the group of reflections
in the walls of the $24$-cell onto the group of reflections in the faces of the cuboctahedron.

The image $\rho_0(\Gamma_{22})$ is isomorphic to $\mathbb{Z}_2 \times \Gamma_{co}$ where the
$\mathbb{Z}_2$ factor comes from reflection in the hyperplane $V$ which is the limit of all of
the positive numbered walls.  In order to better understand the geometry of this
limiting process, we first note that the positive even space-like vectors limit to
$+e_4$ while the positive odd space-like vectors limit to $-e_4$.  This implies that
the polytopes $\mathcal{F}_t$ are being trapped between two hyperplanes, one slightly
above and one slightly below $V$ (which is determined by $\pm e_4$).  The hyperplanes associated
to the even and odd positive space-like vectors, respectively, lie almost parallel to these two hyperplanes
which converge to $V$, forcing the polytope to flatten onto the $3$-dimensional subspace $V$.

However, one can change the metric on a neighborhood of $V$ so that the distance between the 
two squeezing hyperplanes remains constant and the limiting metric on the region between the hyperplanes 
is the product metric on $\mathbb{H}^3 \times I$, $3$-dimensional hyperbolic space crossed with the unit
interval.  It can be arranged that the geometric limit of the polytopes $\mathcal{F}_t$ with 
respect to this sequence of
normal rescalings will be the product of the cuboctahedron and an interval.  The top and bottom walls of this
polytope will be the geometric limits of the even and odd numbered walls, respectively.  Reflection in the walls of this polytope generates a discrete subgroup of the group of isometries of $\mathbb{H}^3 \times I$.

Providing the details of this geometric construction would require more space than perhaps is merited here.  On the other hand, the algebraic manifestation of this limiting process is fairly immediate.
Consider the Coxeter diagram in figure \ref{coxeter_diagram2} for the continuously varying representations
$\rho_t(\U{\Gamma}_{22})$ of the extended group.  As $t\to 0$ the label $\nu$ goes to $2$, signifying
the $\p{0}$ and $\l{N}$ walls and the $\p{3}$ and $\l{M}$ walls are becoming orthogonal.  When
$\nu = 2$ one removes the edge between the orthogonal walls and the Coxeter diagram becomes disconnected.
This Coxeter diagram is pictured in figure \ref{coxeter_diagram9} below.
The corresponding Coxeter group is the direct product of the group $\mathbb{Z}_2 * \mathbb{Z}_2$
generated by reflection in the $\p{0}$ and $\p{3}$ walls and the group generated by the remaining
walls.  This latter group is just the subgroup of $\U{\Gamma}_{22}$ corresponding to the subgroup of $\Gamma_{22}$ generated by the letter and negative numbered walls; this is isomorphic to the group of
reflections in the polyhedron $\U{P}_{co}$ obtained from the cuboctahedron by dividing out by its group of symmetries.
This group of symmetries is realized by the group generated by reflection in the $\l{L}$, $\l{M}$,
and $\l{N}$ walls acting on the hyperplane $V$.

\begin{figure}[ht]
\labellist
\small\hair 2pt
\pinlabel $\p{0}$ [b] at 25 650
\pinlabel $\p{3}$ [t] at 25 4
\pinlabel $\l{M}$ [t] at 451 4
\pinlabel $\l{L}$ [r] at 435 327
\pinlabel $\l{A}$ [b] at 685 360
\pinlabel $\l{N}$ [b] at 451 650
\pinlabel $\m{0}$ [b] at 880 650
\pinlabel $\m{3}$ [t] at 880 4
\endlabellist
\centering
\includegraphics[scale=0.1]{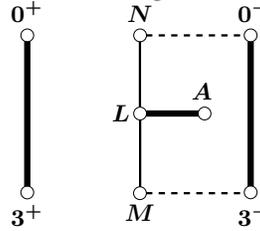}
\caption{A Coxeter diagram for the deformation $\widetilde{\rho}_0 ( \widetilde{\Gamma}_{22} ) $ }
\label{coxeter_diagram9}
\end{figure}

Thus the Coxeter group obtained when $t=0$ (and $\nu=2$) corresponds to the polytope equal to the
product of $\U{P}_{co}$ and an interval. It can be interpreted as the orbifold fundamental group of
this polytope viewed as an orbifold by mirroring its walls.  This group is the homomorphic image of 
the extended group $\U{\Gamma}_{22}$. Restricting to the subgroup corresponding to $\Gamma_{22}$, one
obtains a group isomorphic to the product of $\mathbb{Z}_2 * \mathbb{Z}_2$ with $\Gamma_{co}$ which is isomorphic to the group of reflections in the product polytope  $P_{co} \times I$.  Algebraically, this group is obtained
from $\Gamma_{22}$ by identifying the $4$ generators corresponding to the odd positive numbered walls
to one another and then doing the same to those corresponding to the $4$ even positive numbered walls.
One can see that this group can't be realized as a group of isometries of $\Hf$ because the centralizer 
of various elements is too large.  Hyperbolic geometry forces the identification of all $8$
elements corresponding to the positive numbered walls.  This results in a group isomorphic to
$\mathbb{Z}_2 \times \Gamma_{co}$, which we have seen is isomorphic to $\rho_0(\Gamma_{22})$.

Finally, we note that the phenomenon of a sequence of $n$-dimensional hyperbolic structures,
typically incomplete or complete but with orbifold or cone manifold singularities, collapsing 
to an $(n-1)$-dimensional
hyperbolic space is quite familiar from $3$-dimensional hyperbolic geometry.  For example, in 
\cite{Th} Thurston constructs several examples where a sequence of $3$-dimensional hyperbolic structures on the figure $8$ knot complement collapse to a $2$-dimensional image.  The limit corresponds to
a Dehn filling on the figure $8$ complement that is Seifert fibered with hyperbolic base.  The Seifert 
fibered space has a geometric structure modeled on $\U{SL}_2$, the universal covering of the group of
isometries of $\mathbb{H}^2$, and the limiting representations correspond to the holonomy group 
of a $\U{SL}_2$ structure, acting on $\mathbb{H}^2$.  A similar phenomenon occurs in the proof
of the Orbifold Theorem \cite{CHK,BLP} where certain Seifert fibered orbifolds with 
$\U{SL}_2$ structures
are approximated by $3$-dimensional hyperbolic cone manifold structures.

In our situation, we can view the polytope  $P_{co} \times I$ as a fibered orbifold with fiber
an interval (with mirrored endpoints) and hyperbolic base equal to $P_{co}$, viewed as a 
$3$-dimensional mirrored orbifold.  We have seen that this fibered geometric structure can
be approximated by a sequence of $4$-dimensional hyperbolic structures with (cone manifold) singularities.
Once again, three dimensional hyperbolic geometry pokes its head into dimension four.



\bibliography{Kerckhoff-Storm_2008-05-28}

\begin{thebibliography}{10}

\bibitem{And}
James~W. Anderson.
\newblock A brief survey of the deformation theory of {K}leinian groups.
\newblock In {\em The Epstein birthday schrift}, volume~1 of {\em Geom. Topol.
  Monogr.}, pages 23--49 (electronic). Geom. Topol. Publ., Coventry, 1998.

\bibitem{Andreev}
E.~Andreev.
\newblock Convex polyhedra of finite volume in {Lobachevsky} space.
\newblock {\em Mat. Sb. (N.S.)}, 83 (125):256--260, 1970.

\bibitem{BLP}
M.~Boileau, B.~Leeb, and J.~Porti.
\newblock Geometrization of $3$-dimensional orbifolds.
\newblock {\em Ann. of Math. (2)}, 162(1):195--290, 2005.

\bibitem{Borel1}
A.~Borel.
\newblock {\em Introduction aux groupes arithm\'etiques}.
\newblock Number 1341 in Publications de l'Institut de Math\'ematique de
  l'Universit\'e de Strasbourg, XV. Actualit\'es Scientifiques et
  Industrielles. Hermann, Paris, 1969.

\bibitem{Bow}
B.~Bowditch.
\newblock Geometrical finiteness for hyperbolic groups.
\newblock {\em J. Funct. Anal.}, 113(2):245--317, 1993.

\bibitem{C}
R.~Canary.
\newblock On the {Laplacian} and the geometry of hyperbolic 3-manifolds.
\newblock {\em J. Diff. Geom.}, 36:349--367, 1992.

\bibitem{CHK}
D.~Cooper, C.~Hodgson, and S.~Kerckhoff.
\newblock {\em Three-dimensional Orbifolds and Cone-Manifolds}, volume~5 of
  {\em MSJ Memoirs}.
\newblock Mathematical Society of Japan, 2000.

\bibitem{Cox1}
H.S.M. Coxeter.
\newblock {\em Twelve Geometric Essays}.
\newblock Southern Illinois University Press, Carbondale, Ill., 1968.

\bibitem{Cox2}
H.S.M. Coxeter.
\newblock {\em Regular Complex Polytopes. Second edition}.
\newblock Cambridge University Press, Cambridge, 1991.

\bibitem{dlH}
P.~de~la Harpe.
\newblock An invitation to {C}oxeter groups.
\newblock In {\em Group theory from a geometrical viewpoint ({T}rieste, 1990)},
  pages 193--253. World Sci. Publ., River Edge, NJ, 1991.

\bibitem{GR}
H.~Garland and M.~Raghunathan.
\newblock Fundamental domains for lattices in {$\mathbb{R}$}-rank $1$
  semisimple {Lie} groups.
\newblock {\em Ann. of Math. (2)}, 92:279--326, 1970.

\bibitem{Morris}
D.~Witte Morris.
\newblock Introduction to arithmetic groups.
\newblock Preprint available at \texttt{http://people.uleth.ca/\urltilda
  dave.morris}, 2003.

\bibitem{Mos}
G.~Mostow.
\newblock {\em Strong rigidity of locally symmetric spaces}.
\newblock Number~78 in Annals of Math. Studies. Princeton University Press,
  Princeton, N.J., 1973.

\bibitem{Pr}
G.~Prasad.
\newblock Strong rigidity of {$\mathbb{Q}$}-rank $1$ lattices.
\newblock {\em Invent. Math.}, 21:255--286, 1973.

\bibitem{Rag}
M.~S. Raghunathan.
\newblock {\em Discrete subgroups of {Lie} groups}, volume~68 of {\em
  Ergebnisse der Mathematik und ihrer Grenzgebiete}.
\newblock Springer-Verlag, New York-Heidelberg, 1972.

\bibitem{Sel}
A.~Selberg.
\newblock On discontinuous groups in higher-dimensional symmetric spaces.
\newblock In {\em Contributions to function theory (internat. Colloq. Function
  Theory, Bombay, 1960)}, pages 147--164. Tata Institute of Fundamental
  Research, Bombay, 1960.

\bibitem{Th}
W.~Thurston.
\newblock The topology and geometry of 3-manifolds.
\newblock Available from the MSRI website \texttt{www.msri.org}, 1976--1979.
\newblock Princeton Univ. lecture notes.

\bibitem{Vin1}
E.B. Vinberg.
\newblock On groups of unit elements of certain quadratic forms.
\newblock {\em Math. USSR Sbornik}, 16(1):17--35, 1972.
\newblock (translated by D.E. Brown).

\bibitem{Vin2}
E.B. Vinberg and O.V. Shvartsman.
\newblock Discrete groups of motions of spaces of constant curvature.
\newblock In {\em Geometry. II.}, volume~29 of {\em Encyclop{\ae}dia of
  Mathematical Sciences}, pages 139--248. Springer, Berlin, 1993.

\bibitem{Wa}
Hsien~Chung Wang.
\newblock Topics on totally discontinuous groups.
\newblock In {\em Symmetric spaces (Short Courses, Washington Univ., St. Louis,
  Mo., 1969--1970)}, volume~8 of {\em Pure and Appl. Math.}, pages 459--487.
  Dekker, New York, 1972.

\bibitem{We}
A.~Weil.
\newblock On discrete subgroups of {Lie} groups. {II}.
\newblock {\em Ann. of Math. (2)}, 75:578--602, 1962.

\end{thebibliography}
\end{document}